\documentclass[preprint,12pt,a4paper,review]{elsarticle}

\graphicspath{{Figures/}}
\DeclareGraphicsExtensions{.pdf,.png,.jpg}
\usepackage{amssymb}
\usepackage{amsthm}
\usepackage{graphicx}
\usepackage{amsmath,amsthm,amsbsy}
\usepackage{multirow}
\usepackage{cases}
\usepackage{pgfplots}
\pgfplotsset{compat=1.18}
\usepackage{subcaption}
\usepackage{tikz}
\usepackage{tkz-euclide}
\usetikzlibrary[patterns]
\usepackage{hyperref}
\usepackage{dsfont}
\usepackage{caption}
\usepackage{epstopdf}
\usepackage{epsfig}
\usepackage{array}
\newcolumntype{P}[1]{>{\centering\arraybackslash}p{#1}}
\usepackage{enumerate}
\usepackage{float}
\usepackage[T1]{fontenc}
\usepackage[latin9]{inputenc}
\usepgfplotslibrary{fillbetween}
\usetikzlibrary{patterns}
\usepackage{tikz-dimline}
\usepackage{color}
\usepackage{stmaryrd}
\usepackage{setspace}
\usepackage{amsthm}
\usepackage{textcomp}
\usepackage{soul,xcolor}
\usepackage{multirow}
\usepackage[inline]{enumitem}
\usepackage{pgfplotstable}
\usepackage{booktabs}
\usepackage{algorithmicx}

\usepackage{bm}

\biboptions{sort&compress}
\usepackage{dirtytalk}
\usepackage{csvsimple}

\pgfplotsset{
/pgfplots/bar cycle list/.style={/pgfplots/cycle list={%
    {blue,fill=blue!30!white,mark=none},%
    {red,fill=red!30!white,mark=none},%
    {brown!60!black,fill=brown!30!white,mark=none},%
    {black,fill=gray,mark=none},%
    {violet!80!black,fill=violet,mark=none},%
    {orange,fill=orange!50!white,mark=none}%
    }
},
}
\usepackage[noend]{algpseudocode}
\usepackage{algorithm}

\usetikzlibrary{shapes.geometric}
\usepackage[export]{adjustbox}

\usepackage{hyperref}
\hypersetup{
    colorlinks,
    linkcolor={magenta},
    citecolor={blue},
    urlcolor={blue!80!black},
    breaklinks=true,
	plainpages=true
}
\usepackage{lineno}

\usepackage{bm} 
\usepackage{mathrsfs}

\definecolor{mplblue}{HTML}{1F77B4}
\definecolor{mplorange}{HTML}{FF7F0E}
\newcommand{\ifcomment}{\iffalse}
\newcommand{\bs}[1]{\boldsymbol{#1}}

\def\abs#1{\left|#1\right|}

\newcommand{\G}{\Gamma}

\newcommand{\Om}{\Omega}

\newdefinition{rem}{Remark}

\newcommand{\cref}[2]{\hyperref[#2]{#1~\ref*{#2}}}
\newcommand{\colref}[2]{\hyperref[#2]{#1~\ref*{#2}}}
\newcommand{\eqnref}[1]{\colref{Eq.}{#1}}
\newcommand{\figref}[1]{\colref{Figure}{#1}}

\newcommand{\secref}[1]{\colref{Section}{#1}}

\newcommand{\Algref}[1]{\hyperref[#1]{Algorithm~\ref*{#1}}}

\newcolumntype{M}[1]{>{\centering\arraybackslash}m{#1}}

\definecolor{ActiveElement}{RGB}{147,194,74}
\definecolor{InterceptedElement}{RGB}{255,208,48}
\definecolor{FalseInterceptedElement}{RGB}{92,91,255}

\definecolor{cpu1}{HTML}{4CAF50}
\definecolor{cpu2}{HTML}{FFC107}
\definecolor{cpu3}{HTML}{F44336}
\definecolor{cpu4}{HTML}{2196F3}
\definecolor{cpu5}{HTML}{9932CC}

\definecolor{gpu1}{HTML}{A5D6A7}
\definecolor{gpu2}{HTML}{FFE082}
\definecolor{gpu3}{HTML}{EF9A9A}
\definecolor{gpu4}{HTML}{90CAF9}

\definecolor{sq_b1}{RGB}{37,52,148}
\definecolor{sq_b2}{RGB}{44,127,184}
\definecolor{sq_b3}{RGB}{65,182,196}
\definecolor{sq_b4}{RGB}{127,205,187}
\definecolor{sq_b5}{RGB}{199,233,180}
\definecolor{sq_b6}{RGB}{255,255,204}

\definecolor{sq_r1}{RGB}{189,0,38}
\definecolor{sq_r2}{RGB}{240,59,32}
\definecolor{sq_r3}{RGB}{253,141,60}
\definecolor{sq_r4}{RGB}{254,178,76}
\definecolor{sq_r5}{RGB}{254,217,118}
\definecolor{sq_r6}{RGB}{255,255,178}

\definecolor{sq_g1}{RGB}{0,104,55}
\definecolor{sq_g2}{RGB}{49,163,84}
\definecolor{sq_g3}{RGB}{120,198,121}
\definecolor{sq_g4}{RGB}{173,221,142}
\definecolor{sq_g5}{RGB}{217,240,163}
\definecolor{sq_g6}{RGB}{255,255,204}

\definecolor{div_c1}{RGB}{230,171,2}
\definecolor{div_c2}{RGB}{102,166,30}
\definecolor{div_c3}{RGB}{231,41,138}
\definecolor{div_c4}{RGB}{117,112,179}
\definecolor{div_c5}{RGB}{217,95,2}
\definecolor{div_c6}{RGB}{27,158,119}
\definecolor{div_c7}{RGB}{215,48,39}

\definecolor{div_d1}{RGB}{215,25,28}
\definecolor{div_d2}{RGB}{253,174,97}
\definecolor{div_d3}{RGB}{255,255,191}
\definecolor{div_d4}{RGB}{171,217,233}
\definecolor{div_d5}{RGB}{44,123,182}

\definecolor{lineclr}{RGB}{0,0,0}
\definecolor{utorange}{RGB}{0,0,255}
\definecolor{utsecblue}{RGB}{255,255,0}
\definecolor{utsecgreen}{RGB}{255,0,0}
\definecolor{red!15}{RGB}{0,255,255}
\definecolor{fillclr5}{RGB}{0,255,0}
\definecolor{fillclr6}{RGB}{255,0,255}
\definecolor{fillclr7}{RGB}{255,255,255}
\definecolor{fillclr8}{RGB}{0,0,0}

\definecolor{armygreen}{rgb}{0.29, 0.33, 0.13}
\definecolor{aurometalsaurus}{rgb}{0.43, 0.5, 0.5}
\definecolor{applegreen}{rgb}{0.55, 0.71, 0.0}
\definecolor{darkgreen}{rgb}{0.0, 0.4, 0.25}

\definecolor{c1}{RGB}{49,54,149}    
\definecolor{c2}{RGB}{69,117,180}
\definecolor{c3}{RGB}{116,173,209}
\definecolor{c4}{RGB}{171,217,233}
\definecolor{c5}{RGB}{253,174,97}
\definecolor{c6}{RGB}{244,109,67}
\definecolor{c7}{RGB}{215,48,39}    

\definecolor{sczmA}{HTML}{F28E2B}
\definecolor{ifmA}{HTML}{B35400}

\definecolor{sczmB}{HTML}{4CC9F0}
\definecolor{ifmB}{HTML}{1F4E79}

\def\drawcubeI(#1,#2,#3,#4,#5){ 
\coordinate (O) at (#1,#2,#3);
\coordinate (A) at (#1,#2+#4,#3);
\coordinate (B) at (#1,#2+#4,#3+#4);
\coordinate (C) at (#1,#2,#3+#4);
\coordinate (D) at (#1+#4,#2,#3);
\coordinate (E) at (#1+#4,#2+#4,#3);
\coordinate (F) at (#1+#4,#2+#4,#3+#4);
\coordinate (G) at (#1+#4,#2,#3+#4);
\draw[#5] (O) -- (C) -- (G) -- (D) -- cycle;
\draw[#5] (O) -- (A) -- (E) -- (D) -- cycle;
\draw[#5] (O) -- (A) -- (B) -- (C) -- cycle;
\draw[#5] (D) -- (E) -- (F) -- (G) -- cycle;
\draw[#5] (C) -- (B) -- (F) -- (G) -- cycle;
\draw[#5] (A) -- (B) -- (F) -- (E) -- cycle;
}

\def\drawcubeII(#1,#2,#3,#4,#5,#6,#7){ 
\coordinate (O) at (#1,#2,#3);
\coordinate (A) at (#1,#2+#4,#3);
\coordinate (B) at (#1,#2+#4,#3+#4);
\coordinate (C) at (#1,#2,#3+#4);
\coordinate (D) at (#1+#4,#2,#3);
\coordinate (E) at (#1+#4,#2+#4,#3);
\coordinate (F) at (#1+#4,#2+#4,#3+#4);
\coordinate (G) at (#1+#4,#2,#3+#4);
\draw[#5,fill=#6,opacity=#7] (O) -- (C) -- (G) -- (D) -- cycle;
\draw[#5,fill=#6,opacity=#7] (O) -- (A) -- (E) -- (D) -- cycle;
\draw[#5,fill=#6,opacity=#7] (O) -- (A) -- (B) -- (C) -- cycle;
\draw[#5,fill=#6,opacity=#7] (D) -- (E) -- (F) -- (G) -- cycle;
\draw[#5,fill=#6,opacity=#7] (C) -- (B) -- (F) -- (G) -- cycle;
\draw[#5,fill=#6,opacity=#7] (A) -- (B) -- (F) -- (E) -- cycle;
}

\def\drawNodes(#1,#2,#3,#4,#5,#6,#7){ 
\foreach \x in {#1,#7,...,#2}{
	\foreach \y in {#3,#7,...,#4}{
		\foreach \z in {#5,#7,...,#6}{
				\draw[fill=red!60] (\x,\y,\z) circle (0.15);
				}
			}
	}

}

\pgfplotsset{
  log x ticks with fixed point/.style={
      xticklabel={
        \pgfkeys{/pgf/fpu=true}
        \pgfmathparse{exp(\tick)}%
        \pgfmathprintnumber[fixed relative, precision=3]{\pgfmathresult}
        \pgfkeys{/pgf/fpu=false}
      }
  },
  log y ticks with fixed point/.style={
      yticklabel={
        \pgfkeys{/pgf/fpu=true}
        \pgfmathparse{exp(\tick)}%
        \pgfmathprintnumber[fixed relative, precision=3]{\pgfmathresult}
        \pgfkeys{/pgf/fpu=false}
      }
  }
}

\makeatletter
\newcommand\resetstackedplots{
\makeatletter
\pgfplots@stacked@isfirstplottrue
\makeatother
\addplot [forget plot,draw=none] coordinates{(48,0) (96,0) (192,0) (384,0) (768,0) (1536,0) (3072,0) (6144,0)};
}
\makeatother

\makeatletter
\newcommand\resetstackedplotsOne{
\makeatletter
\pgfplots@stacked@isfirstplottrue
\makeatother
\addplot [forget plot,draw=none] coordinates{(384,0) (768,0) (1536,0) (3072,0) (6144,0)};
}
\makeatother

\makeatletter
\newcommand\resetstackedplotsTwo{
\makeatletter
\pgfplots@stacked@isfirstplottrue
\makeatother
\addplot [forget plot,draw=none] coordinates{(16,0) (32,0) (64,0) (128,0) (256,0) (512,0) (1024,0) (2048,0) (4096,0) (8192,0) (16384,0) (32768,0)};
}
\makeatother

\makeatletter
\newcommand\resetstackedplotsThree{
\makeatletter
\pgfplots@stacked@isfirstplottrue
\makeatother
\addplot [forget plot,draw=none] coordinates{(2,0) (4,0) (8,0) (16,0) (32,0) (64,0)};
}
\makeatother

\makeatletter
\newcommand\resetstackedplotsFour{
\makeatletter
\pgfplots@stacked@isfirstplottrue
\makeatother
\addplot [forget plot,draw=none] coordinates{(4,0) (8,0) (16,0) (32,0) (64,0)};
}
\makeatother

\makeatletter
\newcommand\resetstackedplotsFive{
\makeatletter
\pgfplots@stacked@isfirstplottrue
\makeatother
\addplot [forget plot,draw=none] coordinates{(1,0) (2,0) (4,0) (8,0) (16,0) (32,0) (64,0) (128,0)};
}
\makeatother

\makeatletter
\newcommand\resetstackedplotsSix{
\makeatletter
\pgfplots@stacked@isfirstplottrue
\makeatother
\addplot [forget plot,draw=none] coordinates{(2,0) (4,0) (8,0) (16,0) (32,0) (64,0) (128,0)};
}
\makeatother

\newcommand{\jumpscalar}[1]{\left\llbracket #1 \right\rrbracket}

\newcommand{\neml}{\href{https://applied-material-modeling.github.io/neml2/}{\texttt{NEML2}}}
\newcommand{\moose}{\href{https://mooseframework.inl.gov/}{\texttt{MOOSE}}}

\newcommand{\sculpt}{\href{https://cubit.sandia.gov/files/cubit/16.10/help_manual/WebHelp/mesh_generation/meshing_schemes/parallel/sculpt.htm}{\textsc{Sculpt}}}
\newcommand{\gmsh}{\href{https://gmsh.info/}{\textsc{Gmsh}}}

\newcommand{\SBMElem}{\textsc{SBMElem}}

\allowdisplaybreaks

\journal{CMAME}

\begin{document}

\begin{frontmatter}

\title{A Shifted Cohesive-Zone Method for Non-Interface-Fitted Meshes with Applications to Crystal Plasticity}
\author[ANL]{Cheng-Hau Yang\texorpdfstring{\corref{cor}}{}}
\ead{chenghau.yang@anl.gov}
\author[ANL]{Mark C. Messner}
\ead{messner@anl.gov}
\author[ANL]{Tianchen Hu}
\ead{thu@anl.gov}
\cortext[cor]{Corresponding authors}
\address[ANL]{Argonne National Laboratory, 9700 S Cass Ave, Lemont, IL 60439, USA}

\begin{abstract}
The accurate simulation of interface-dominated solid mechanics problems on complex microstructures remains challenging, particularly when interface-fitted quadrilateral or hexahedral meshes are difficult to generate. We extend the shifted boundary method (SBM) to cohesive-zone formulations and introduce the Shifted Cohesive Zone Method (SCZM), with applications to crystal plasticity on non-interface-fitted meshes. By shifting the enforcement of traction-separation laws from the true interface to a nearby surrogate interface, SCZM enables the use of standard finite element spaces while avoiding the meshing burden associated with interface-conformal discretizations. We present a simplified SCZM weak form defined on the surrogate interface, leading to a straightforward implementation of the nonlinear residual and consistent tangent matrix. The method is implemented in the open-source \moose{} framework and coupled with constitutive models from \neml{}, enabling simulations with linear elasticity, multiple traction-separation laws, and history-dependent crystal plasticity. We further develop a geometry-aware, PCA-enhanced point classification algorithm to accelerate surrogate-domain construction. Verification and benchmark studies in two and three dimensions demonstrate that SCZM achieves first-order convergence for non-interface-fitted interface problems and closely matches interface-fitted reference solutions in terms of reaction forces, surface energy release, deformation, stress fields, and damage evolution. These results indicate that SCZM provides an accurate and efficient framework for modeling interface mechanics in complex microstructures without requiring interface-fitted meshes.
\end{abstract}


\end{frontmatter}

\section{Introduction}
\label{Sec:Intro}

Interface conditions play a central role in many partial differential equations arising in solid and fluid mechanics~\cite{wells2001new,lee2003immersed,bansch2004finite,xiao2020high,zhao2022enriched}. In particular, problems involving material interfaces, cracks, and grain boundaries require the accurate imposition of discontinuities in the solution or its derivatives across lower-dimensional manifolds. In many applications of crystal plasticity, particularly in representative volume element (RVE) simulations of polycrystalline or composite microstructures, interface behavior often controls the macroscopic response. Grain boundaries, phase interfaces, and particle-matrix interfaces can undergo debonding, sliding, or damage, contributing significantly to inelastic deformation, surface energy release, and eventual failure. These interface-mediated mechanisms are especially important in advanced structural materials, including metal matrix composites and high-temperature alloys, where interfacial degradation influences creep, stress relaxation, and long-term aging behavior~\cite{chen2023icme,hu2023mechanistic,baweja2025predicting,baweja2026development,aduloju2025modeling,bhatt2024microstructural}. Moreover, recent efforts to couple crystal plasticity with grain boundary models further highlight the importance of explicitly accounting for interfacial mechanisms in predictive simulations of creep and damage evolution~\cite{nassif2019combined}. While crystal plasticity models accurately capture intragranular deformation through slip-based mechanisms, they typically assume perfect bonding between neighboring material regions unless augmented with an explicit interface model. Cohesive-zone formulations provide a natural and physically consistent framework to incorporate interface softening and damage within crystal plasticity simulations. However, their application in RVE-scale computations is often hindered by the requirement of interface-fitted meshes, particularly for complex three-dimensional microstructures. This limitation motivates the development of methods that enable cohesive interface modeling within crystal plasticity frameworks while retaining the flexibility of non-interface-fitted discretizations.

In recent years, significant effort has been devoted to the development of non-interface-fitted finite element methods that enable the use of meshes independent of interface geometry. A common strategy is to enrich the approximation space locally so that discontinuities can be represented within elements intersected by the interface. Representative approaches include the Immersed Finite Element Method (IFEM)~\cite{li2003new,adjerid2023enriched} and the extended Finite Element Method (XFEM)~\cite{moes1999finite,xiao2020high}. Another class of methods, such as the Enriched Immersed Boundary Method (EIBM)~\cite{zhao2022enriched}, introduces additional degrees of freedom on interface elements and enforces interface conditions weakly, for example via Nitsche's method. While these approaches are effective, they rely on modifying the approximation space locally, leading to non-uniform function spaces across the computational domain.

An alternative paradigm is to enforce interface conditions without altering the underlying finite element space. The Shifted Fracture Method (SFM)~\cite{li2021shifted,li2023simple,li2025crack} follows this philosophy by shifting the enforcement of interface conditions from the true interface to nearby mesh entities that define a surrogate interface. This approach preserves a uniform approximation space while retaining accuracy on non-interface-fitted meshes. SFM can be viewed as a specialization of the more general Shifted Boundary Method (SBM)~\cite{main2018shifted,atallah2024nonlinear,zeng2025shifted}, which has been successfully applied to a wide range of problems, including fluid flow~\cite{main2018shifted,Main2018TheSB,YANG2026114334}, solid and fracture mechanics~\cite{atallah2021shifted,li2021shifted,li2023simple,li2023blended,li2025shifted}, moving-interface problems~\cite{li2020shifted,colomes2021weighted}, thermal flows~\cite{YANG2025114333}, and moving-boundary problems~\cite{xu2024weighted,xu2025weighted}.

Building on this foundation, several SBM variants have been proposed to enhance accuracy, robustness, and flexibility, including O-SBM (optimal surrogate selection)~\cite{yang2024optimal}, W-SBM~\cite{colomes2021weighted}, Gap-SBM~\cite{collins2025gap,antonelli2026isogeometric}, and GSBM~\cite{colomes2026generalized}. In addition, extensions to octree-based discretizations (Octree-SBM) enable adaptive mesh refinement and demonstrate excellent parallel scalability~\cite{yang2024optimal,YANG2026114334,YANG2025114333,yang2025octree,KARKI2025118248}. Despite these advances, existing SFM formulations have primarily been restricted to simplicial meshes (triangles and tetrahedra), limiting their applicability in many engineering settings.

\begin{figure}[htb]
\centering
\begin{subfigure}[b]{0.48\textwidth}
    \centering
    \includegraphics[width=\linewidth, trim={50 10 150 120}, clip]{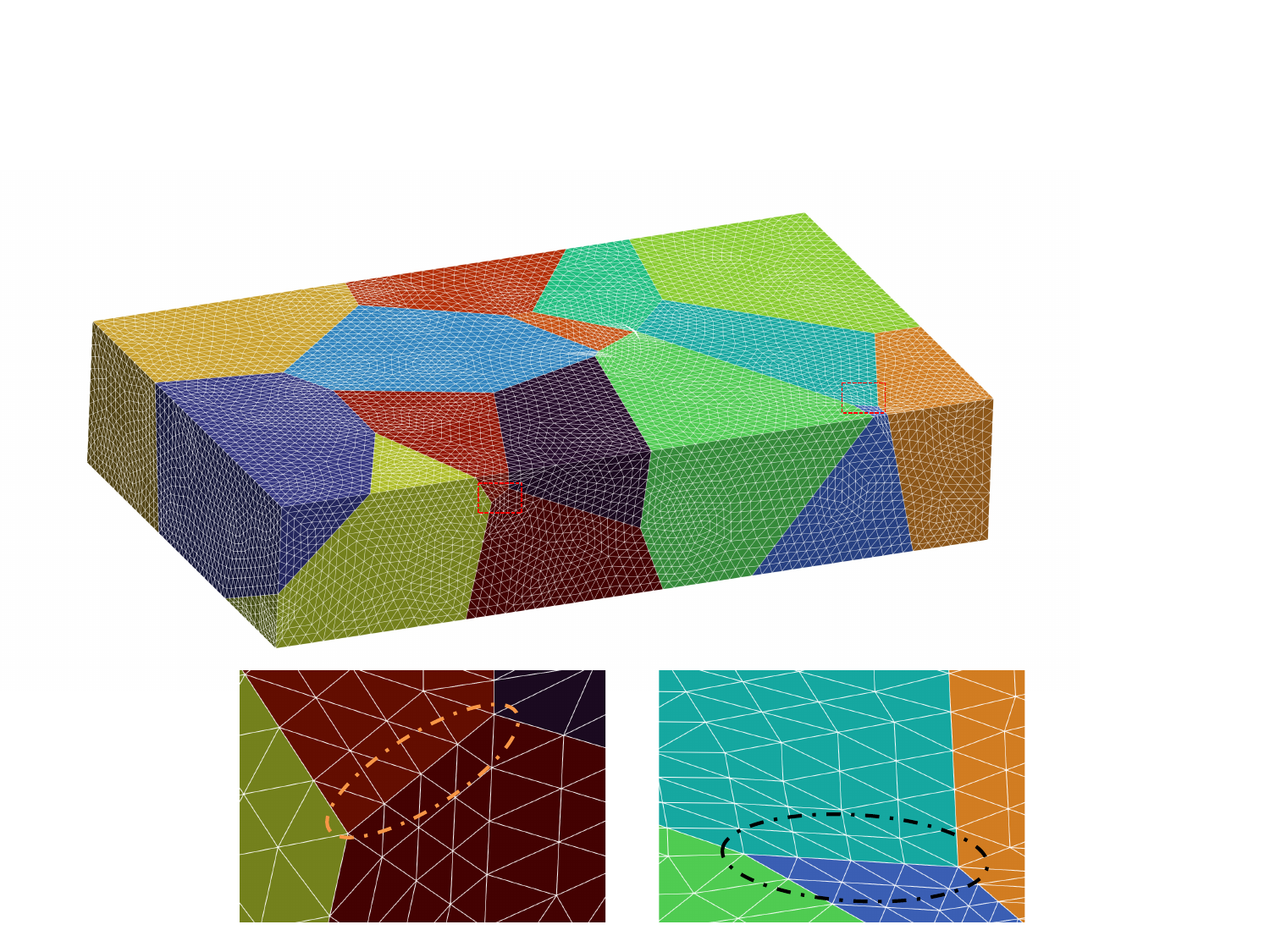}
    \caption{Interface-conformal tetrahedral mesh}
    \label{fig:mesh_comparison_tet}
\end{subfigure}
\hfill
\begin{subfigure}[b]{0.48\textwidth}
    \centering
    \includegraphics[width=\linewidth, trim={50 10 150 120}, clip]{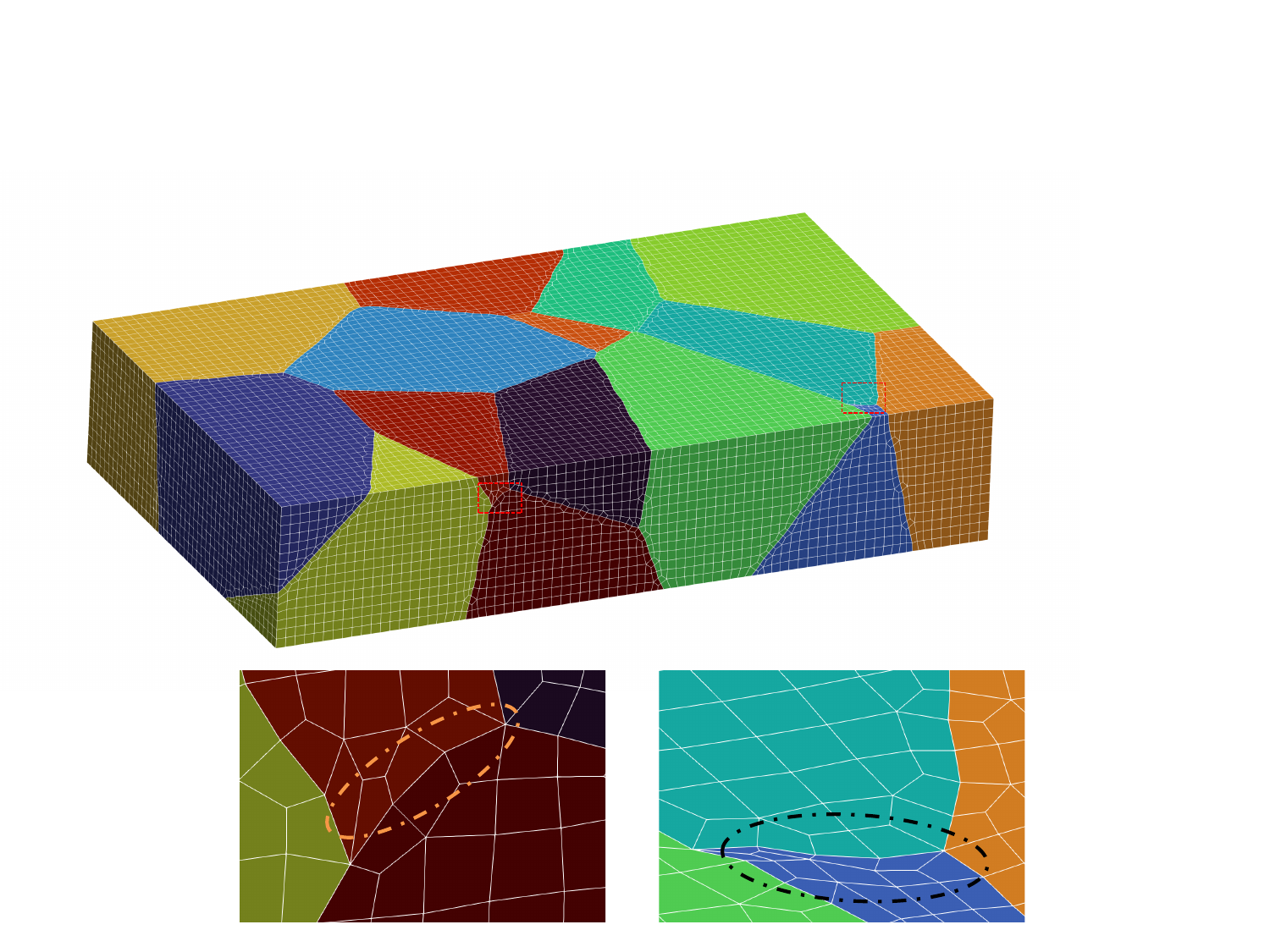}
    \caption{Pseudo-hexahedral mesh with distorted interfaces}
    \label{fig:mesh_comparison_pseudo}
\end{subfigure}
\caption{
Comparison of interface representation. 
(a) Tetrahedral meshes conform well to interfaces. 
(b) Hexahedral meshing may result in pseudo-conformal interfaces with geometric inconsistencies.
}
\label{fig:mesh_comparison_top}
\end{figure}

\begin{figure}[htb]
\centering
\includegraphics[width=0.7\linewidth, trim={50 150 150 120}, clip]{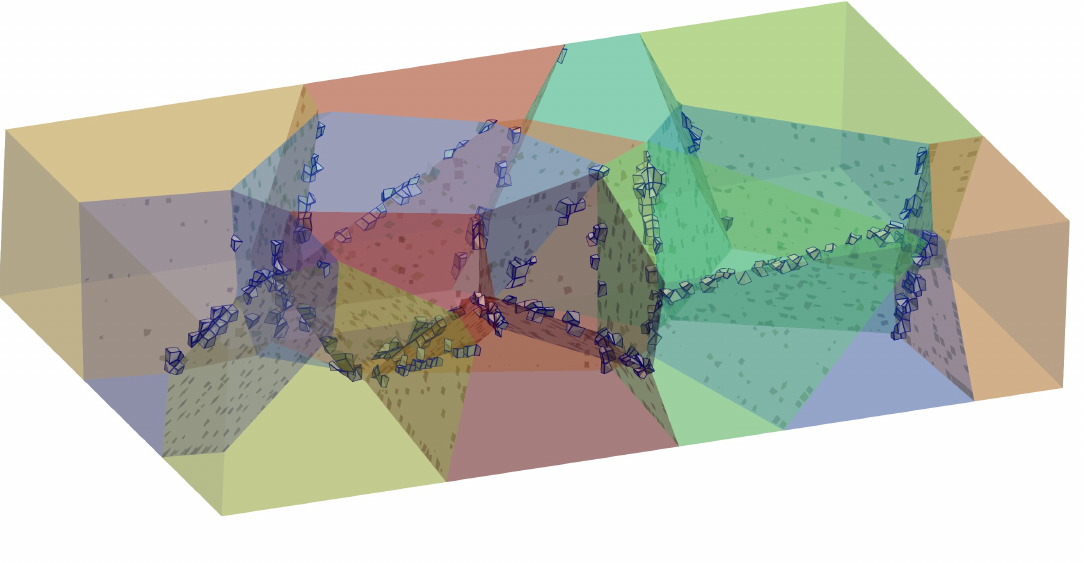}
\caption{
Hexahedral mesh with negative scaled Jacobian elements resulting from enforcing strict interface conformity.
}
\label{fig:mesh_comparison_negative}
\end{figure}

In solid mechanics applications, quadrilateral and hexahedral elements are often preferred over simplicial elements due to their superior accuracy in representing derivatives, even with low-order shape functions~\cite{matouvs2004finite,dohrmann2000node,neto2005f,de2008computational,cheng2016stabilized}. This advantage is particularly important in nonlinear problems, such as plasticity, where these elements help mitigate volumetric locking through appropriate stabilization techniques~\cite{bower2009applied,hughes1987finite,nakshatrala2008fem}. However, generating interface-fitted quadrilateral or hexahedral meshes remains a significant challenge. Standard meshing tools such as \gmsh{} often struggle to produce high-quality interface-conforming meshes for complex geometries (see \figref{fig:mesh_comparison_tet}), and more advanced tools such as \sculpt{} may produce only approximately conforming (pseudo-hexahedral) meshes when sharp features or intricate geometries are present (see \figref{fig:mesh_comparison_pseudo}). Moreover, enforcing strict interface conformity in hexahedral meshes can lead to severely distorted elements, including those with negative scaled Jacobians (see \figref{fig:mesh_comparison_negative}), which compromise numerical stability. These challenges motivate the development of methods that relax the requirement of interface-fitted discretizations while maintaining accuracy. In this context, shifted methods offer a compelling alternative by enabling flexible meshing strategies while enforcing interface conditions in a variationally consistent manner.

In this work, we extend the shifted framework to cohesive-zone formulations and introduce the Shifted Cohesive Zone Method (SCZM). First introduced in~\cite{li2020shifted,li2021shifted,li2023simple,li2023blended,li2025crack}, SCZM has not previously been developed or demonstrated in detail. The proposed method enables the simulation of interface-dominated problems, including crystal plasticity with interface softening mechanisms, on non-interface-fitted meshes. In particular, we demonstrate that SCZM can be applied to quadrilateral and hexahedral discretizations without sacrificing accuracy, thereby addressing a key limitation of existing shifted fracture formulations.

Our key contributions are:
\begin{itemize}
    \item \textit{First demonstration of SCZM for crystal plasticity}. We extend the shifted cohesive-zone framework to history-dependent constitutive models and demonstrate its applicability to crystal plasticity problems.

    \item \textit{Simplified SCZM weak formulation}. We derive a reduced weak form on the surrogate interface that enables a straightforward implementation of the nonlinear residual and its consistent tangent matrix.

    \item \textit{Open-source implementation}. We provide the first implementation of SCZM within the \moose{} framework, enabling simulations on complex microstructures while integrating seamlessly with existing physics modules and constitutive models, including \neml{}.

    \item \textit{Geometry-aware point classification algorithm}. We develop a PCA-enhanced, geometry-aware point classification approach that significantly accelerates surrogate-domain construction within the SCZM/SBM framework.

    \item \textit{Comprehensive verification and validation}. We demonstrate the accuracy and robustness of the proposed method across multiple element types (triangular, quadrilateral, and hexahedral), traction-separation laws, and bulk material models, including linear elasticity and crystal plasticity.
\end{itemize}

\begin{figure}[htbp]
\centering
\includegraphics[width=0.85\linewidth]{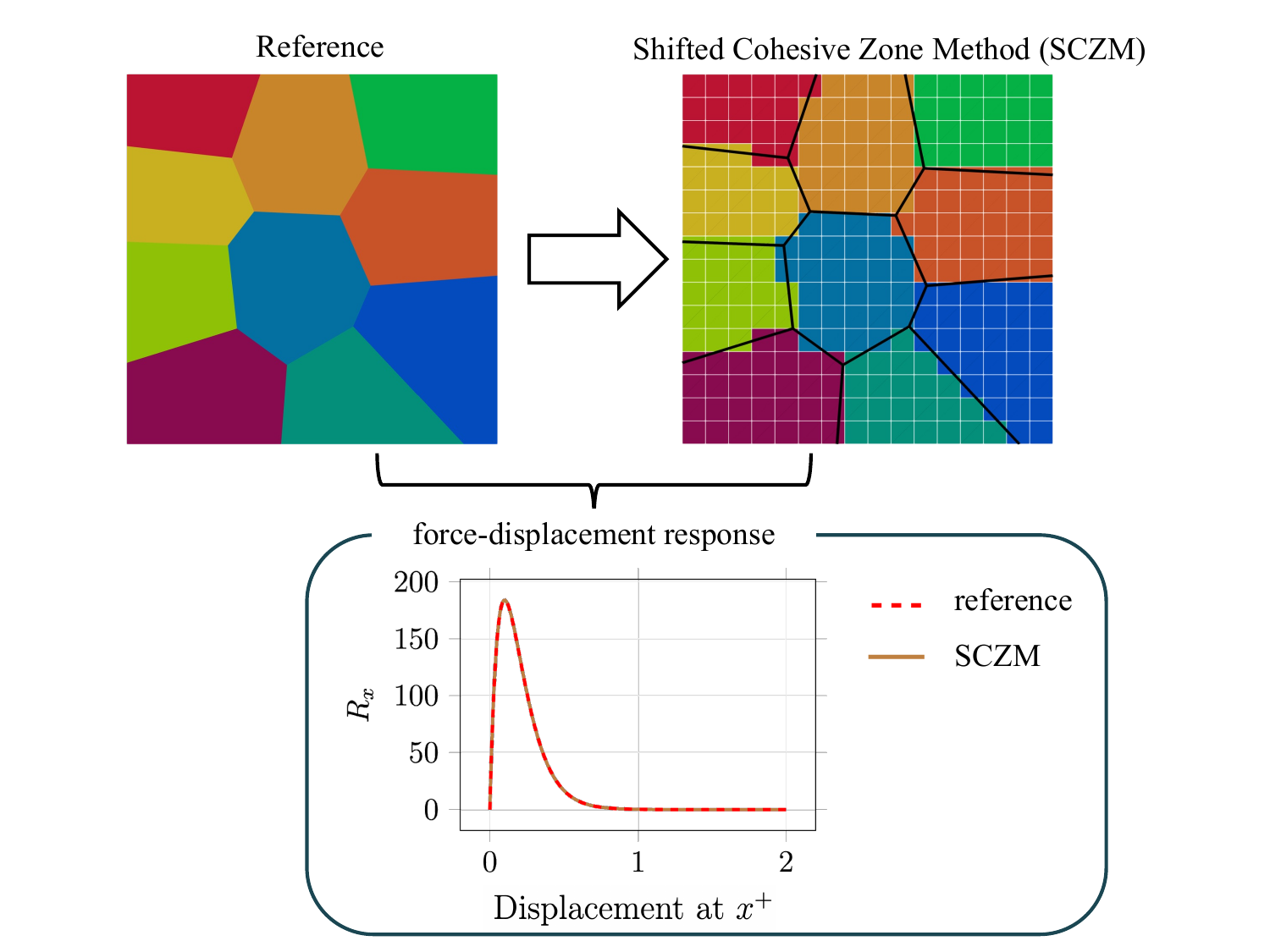}
\caption{
Schematic overview of the proposed SCZM framework. 
Interface conditions are enforced on a surrogate interface defined on a non-interface-fitted mesh. 
The method enables cohesive-zone modeling without requiring interface-fitted discretization, while maintaining accuracy comparable to that of interface-fitted approaches.
}
\label{fig:sczm_overview}
\end{figure}

To provide a high-level understanding of the proposed framework, a schematic overview of the SCZM methodology is presented in \figref{fig:sczm_overview}.

\section{Theory}
\label{sec:Math}

\subsection{Governing equations}

Let $\Omega$ denote a bounded Lipschitz domain in $\mathbb{R}^{d}$ representing a solid, where the spatial dimension, $d$, is either 2 or 3. Let $\partial\Omega$ denote the external boundary of the solid. Let $\partial\Omega_D$ and $\partial\Omega_N$ denote the Dirichlet and Neumann boundaries of the solid, respectively, which jointly partition the external boundary as $\partial\Omega = \partial\Omega_D \cup \partial\Omega_N$. In addition, let $\G$ denote the internal interface of interest, on which the traction continuity is considered to be a constitutive choice. The internal interface $\G$ is assumed to be two-sided and orientable so it admits a unit normal field $\bs{n}$ almost everywhere. Moreover, the internal interface can be embedded, i.e. it does not necessarily intersect with the external boundaries, i.e. $\G \cap \partial\Omega = \emptyset$. Both the external boundary and the interface are assumed to be Lipschitz continuous. Static linear momentum balance of the solid, neglecting body force, can be expressed as:
\begin{subequations}\label{eq:momentum}
\begin{align}
\nabla \cdot \boldsymbol{\sigma} = \mathbf{0} &\quad \text{in } \Omega, \\
\bs{\sigma} \bs{n} = \bs{t}_N  &\quad \text{on } \partial\Omega_N, \\
\bs{u} = \bs{g}  &\quad \text{on } \partial\Omega_D, \\
\bs{\sigma}\bs{n} = \pm\bs{t}_{\text{coh}} &\quad  \text{on } \G^\pm, \label{eq:InterfaceTraction}
\end{align}
\end{subequations}
where $\boldsymbol{\sigma}$ is the Cauchy stress, $\bs{n}$ is the outward normal, $\bs{t}_N$ is the prescribed traction, $\bs{g}$ is the prescribed displacement, and we denote the interface traction by $\bs{t}_{\G} = \bs{t}_{\mathrm{coh}}$. $\nabla$ denotes derivatives with respect to the current configuration. The superscript in $\G^\pm$ signifies the two sides of the interfaces, with no particular choice of sign convention.

\subsection{Weak form}

We introduce the discrete trial and test function spaces as follows:
\begin{subequations}\label{eq:FunctionSpace}
\begin{align}
\mathcal{S} &= \left\{ \bs{v} \in \left[L^2(\Omega)\right]^d : \bs{v} \in \left[H^1(\Omega \setminus \G)\right]^d, \bs{v}|_{\partial\Omega_D} = \bs{g} \right\}, \\
\mathcal{V} &= \left\{ \bs{v} \in \left[L^2(\Omega)\right]^d: \bs{v} \in \left[H^1(\Omega \setminus \G)\right]^d, \bs{v}|_{\partial\Omega_D} = \bs{0} \right\}.
\end{align}
\end{subequations}
In other words, the space $\mathcal{S}$ consists of vector-valued functions that are $H^1$ away from $\G$, i.e., on each connected component of $\Omega\setminus\G$, with no continuity enforced across $\G$.

The weak form of \eqnref{eq:momentum} follows as: Find $\bs{u} \in \mathcal{S}$ such that $a(\bs{u},\bs{\phi}) = \ell(\bs{\phi})$, $\forall \bs{\phi} \in \mathcal{V}$, where
\begin{subequations}
\begin{align}\label{eq:Momentum_full}
a(\bs{u},\bs{\phi}) &= \int\limits_{\Om \setminus \G} \nabla \bs{\phi} : \bs{\sigma} \ \mathrm{d}V - \int\limits_{\G^+} \jumpscalar{\bs{\phi}} \cdot \bs{t}_{\text{coh}} \ \mathrm{d}S, \\
\ell(\bs{\phi}) &= \int\limits_{\partial\Omega_N} \bs{\phi} \cdot \bs{t}_N \ \mathrm{d}S,
\end{align}
\end{subequations}
where $\jumpscalar{z}_\G = \left.z\right\rvert_{\G^+} - \left.z\right\rvert_{\G^-}$ denotes the jump of a function $z$ across a two-sided interface $\G$.

\subsection{The shifted cohesive zone method}

The shifted cohesive zone method (SCZM) approximates the interface contribution in the weak form on a nearby surrogate interface. Starting from the exact interface term in \eqref{eq:Momentum_full}, we seek to express the integral over the true interface $\G$ in terms of an integral over a perturbed surface $\G_\epsilon$.

Let $\G \subset \mathbb{R}^d$ be a compact $C^2$ hypersurface with unit normal $\bs{n}$. Let $\{\G_\epsilon\}_{|\epsilon|<\epsilon_0}$ denote a $C^1$ perturbation of $\G$ generated by a smooth vector field $\bs{D} \in C^1(\mathcal{U};\mathbb{R}^d)$ defined on a tubular neighborhood $\mathcal{U}$ of $\G$, such that
\begin{equation}
    \G_\epsilon = \Phi_\epsilon(\G), \quad \frac{d}{d\epsilon}\Phi_\epsilon(\bs{x}) = \bs{D}(\Phi_\epsilon(\bs{x})), \quad \Phi_0(\bs{x}) = \bs{x}.
\end{equation}
Assuming a well-defined closest-point projection $\Pi: \G_\epsilon \rightarrow \G$, the interface integral in \eqref{eq:Momentum_full} can be expressed on $\G_\epsilon$ via change of variables:
\begin{equation}
\begin{split}
    \int_{\G^+} \jumpscalar{\bs{\phi}}(\bs{x}) \cdot \bs{t}_{\text{coh}}(\bs{x}) \, \mathrm{d}S
    =
    \int_{\G_\epsilon^+}
    \jumpscalar{\bs{\phi}}(\Pi(\bs{y})) \cdot \bs{t}_{\text{coh}}(\Pi(\bs{y}))
    J^\Pi(\bs{y}) \, \mathrm{d}S_\epsilon,
\end{split}
\label{eq:cov_rewrite}
\end{equation}
where $J^\Pi$ is the surface Jacobian of the projection.

To obtain a tractable approximation, we employ a first-order surface-transport expansion. For a sufficiently smooth function $g$ defined in $\mathcal{U}$, the pull-back satisfies
\begin{equation}\label{eq:surftrans}
    g(\Pi(\bs{x})) = g(\bs{x}) - \epsilon(\bs{x}) \nabla_{\bs{n}} g(\bs{x}) + \mathcal{O}(\epsilon^2).
\end{equation}
Applying this expansion to the integrand in \eqref{eq:cov_rewrite} and neglecting higher-order terms yields
\begin{equation}
\begin{split}
    \int_{\G^+} \jumpscalar{\bs{\phi}} \cdot \bs{t}_{\text{coh}} \, \mathrm{d}S
    \approx
    \int_{\G_\epsilon^+}
    \left[
    \jumpscalar{\bs{\phi}} \cdot \bs{t}_{\text{coh}}
    -
    \epsilon \nabla_{\bs{n}} \left( \jumpscalar{\bs{\phi}} \cdot \bs{t}_{\text{coh}} \right)
    \right]
    J^\Pi \, \mathrm{d}S_\epsilon.
\end{split}
\label{eq:SurfTrans_rewrite}
\end{equation}

The surface Jacobian admits the geometric decomposition
\begin{equation}
    J^\Pi = |\bs{n} \cdot \bs{n}_\epsilon| \det(\bs{I} - \epsilon \bs{S})
    =
    |\bs{n} \cdot \bs{n}_\epsilon| \left( 1 - \epsilon \kappa + \mathcal{O}(\epsilon^2) \right),
\end{equation}
where $\bs{n}_\epsilon$ is the normal on $\G_\epsilon$, $\bs{S}$ is the Weingarten map, and $\kappa$ is the mean curvature. Substituting into \eqref{eq:SurfTrans_rewrite} and retaining first-order terms gives
\begin{equation}
\begin{split}
    \int_{\G^+} \jumpscalar{\bs{\phi}} \cdot \bs{t}_{\text{coh}} \, \mathrm{d}S
    \approx
    &\int_{\G_\epsilon^+}
    (1 - \epsilon \kappa)
    \jumpscalar{\bs{\phi}} \cdot \bs{t}_{\text{coh}}
    |\bs{n} \cdot \bs{n}_\epsilon| \, \mathrm{d}S_\epsilon \\
    &-
    \int_{\G_\epsilon^+}
    \epsilon \nabla_{\bs{n}}
    \left( \jumpscalar{\bs{\phi}} \cdot \bs{t}_{\text{coh}} \right)
    |\bs{n} \cdot \bs{n}_\epsilon| \, \mathrm{d}S_\epsilon.
\end{split}
\label{eq:shifted_interface_term}
\end{equation}

Substituting \eqref{eq:shifted_interface_term} into the weak form \eqref{eq:Momentum_full} yields a first-order shifted formulation in which the interface contribution is evaluated on $\G_\epsilon$. The first term represents a geometric \textit{area correction}, accounting for normal misalignment and curvature, while the second term is a \textit{field correction} arising from the Taylor expansion of the integrand.

\paragraph{Discrete formulation and surrogate interface}

In practice, the surrogate interface is constructed from mesh facets. Let $\mathcal{T}_h$ be a tessellation of $\Omega$ that is not required to conform to $\G$, and let $\mathcal{F}_h^\mathrm{int}$ denote the set of internal facets. The discrete surrogate interface is defined as
\begin{equation}
    \G_h = \bigcup_{F \in \mathcal{F}_h^\G} F,
    \quad
    \mathcal{F}_h^\G = \{ F \in \mathcal{F}_h^\mathrm{int} : F \cap \G \ne \emptyset \}.
\end{equation}
Motivated by \eqref{eq:shifted_interface_term}, and noting that the overall error is dominated by the geometric approximation $\G_h \to \G$, we retain the leading-order geometric contribution and neglect higher-order curvature and field-correction terms. This yields the discrete shifted bilinear form
\begin{equation}
\begin{split}
    a_h(\bs{u},\bs{\phi}) =
    \int_{\Omega \setminus \G_h} \nabla \bs{\phi} : \bs{\sigma} \, \mathrm{d}V
    -
    \int_{\G_h^+}
    \jumpscalar{\bs{\phi}} \cdot \bs{t}_{\text{coh}}
    |\bs{n} \cdot \bs{n}_h| \, \mathrm{d}S_h,
\end{split}
\label{eq:discrete_shifted}
\end{equation}
where $\bs{n}_h$ denotes the unit normal to the facet.

\paragraph{Normal mismatch and directional correction}

Because $\G_h$ is composed of mesh facets, its normal $\bs{n}_h$ generally differs from the true interface normal $\bs{n}$. This mismatch admits the decomposition
\begin{equation}
    \bs{n}_h = (\bs{n}_h \cdot \bs{n}) \bs{n} + \bs{\tau}_h,
    \qquad
    \bs{\tau}_h = \bs{n}_h - (\bs{n}_h \cdot \bs{n}) \bs{n},
\end{equation}
where $\bs{\tau}_h$ is tangential to the true interface. Using the interface traction relation $\bs{\sigma}^\pm \bs{n}^\pm = \pm \bs{t}_{\text{coh}}$, the traction associated with the surrogate normal can be expressed as
\begin{equation}
    \bs{\sigma}^\pm \bs{n}_h^\pm
    =
    (\bs{n}_h^\pm \cdot \bs{n}^\pm) \bs{t}_{\text{coh}}
    +
    \bs{\sigma}^\pm \bs{\tau}_h^\pm.
\end{equation}
This decomposition reveals that, in addition to the normal component captured by the geometric factor $|\bs{n}\cdot\bs{n}_h|$, a tangential contribution arises due to the misalignment between $\bs{n}_h$ and $\bs{n}$. Incorporating this effect leads to the expanded discrete formulation
\begin{equation}
\begin{split}
    a_h(\bs{u},\bs{\phi}) =
    &\int_{\Omega \setminus \G_h} \nabla \bs{\phi} : \bs{\sigma} \, \mathrm{d}V
    -
    \int_{\G_h^+}
    \jumpscalar{\bs{\phi}} \cdot \bs{t}_{\text{coh}}
    |\bs{n} \cdot \bs{n}_h| \, \mathrm{d}S_h \\
    &-
    \int_{\G_h^+} \bs{\phi}^+ \cdot \bs{\sigma}^+ \bs{\tau}_h \, \mathrm{d}S_h
    +
    \int_{\G_h^-} \bs{\phi}^- \cdot \bs{\sigma}^- \bs{\tau}_h \, \mathrm{d}S_h.
\end{split}
\label{eq:directional_correction}
\end{equation}
The last two terms constitute the \textit{directional correction}, which accounts for the tangential traction induced by normal mismatch on the surrogate interface. This correction \footnote{This interpretation does not imply exact equivalence to the original interface formulation, but provides a consistent discrete approximation that captures leading-order geometric and traction effects.} arises naturally from the traction vector $\bs{\sigma}\bs{n}_h$ associated with the discrete facet geometry and improves the consistency of the surrogate-interface formulation.

\subsection{Constitutive traction-separation laws}
\label{sec:TSL}

The cohesive response along the interface is described through a traction-separation law (TSL), which relates the interface traction to the displacement jump and, in general, a set of internal variables. In its most general form, the cohesive traction can be written as
\begin{equation}
    \bs{t}_{\text{coh}} = \mathcal{F}\left(\jumpscalar{\bs{u}}, \boldsymbol{\alpha}\right),
\end{equation}
where $\jumpscalar{\bs{u}}$ is the displacement jump across the interface and $\boldsymbol{\alpha}$ denotes a set of internal variables governing history-dependent behavior (e.g., damage, plasticity, or softening).

Within the shifted cohesive zone framework, the interface contribution is evaluated on a surrogate interface rather than the true interface. As a result, the displacement jump must be consistently approximated at the surrogate location. To this end, we employ a first-order surface-transport approximation, yielding the shifted displacement jump approximation
\begin{equation}
    \left.\jumpscalar{\bs{u}}\right\rvert_{\Pi(\bs{x})}
    =
    \left.\left(\jumpscalar{\bs{u}} + \jumpscalar{\nabla \bs{u}} \cdot \bs{d}\right)\right\rvert_{\bs{x}}, \quad \forall \bs{x} \in \Gamma_h,
\end{equation}
where $\bs{d}$ denotes the vector from a point on the surrogate interface to its closest-point projection on the true interface. This expression is obtained by applying the same surface-transport expansion used in the derivation of the shifted weak form, and represents a first-order approximation of the displacement field evaluated at the true interface.

The cohesive traction in SCZM is then defined as
\begin{equation}
    \bs{t}_{\text{coh}} = \mathcal{F}\left(\left.\jumpscalar{\bs{u}}\right\rvert_{\Pi(\bs{x})}, \boldsymbol{\alpha}\right),
\end{equation}
i.e., the functional form of the traction-separation law remains unchanged, and only the kinematic argument is modified through the shifted displacement jump. This construction allows existing cohesive models to be incorporated into the SCZM framework without modification.

In the present study, several representative traction-separation laws are used in the numerical examples, including exponential softening laws~\cite{needleman1987continuum,needleman1990analysis,xu1994numerical}, bilinear mixed-mode formulations~\cite{Camanho2002}, and coupled three-dimensional models by \citet{salehani2018coupled}. These models are employed solely for demonstration purposes and do not restrict the generality of the proposed framework.

\section{Implementation details}
\label{sec:implementation}

We consider polycrystalline representative volume elements (RVEs) in which the microstructure is described by a set of grain domains separated by interfaces. In the absence of interface-fitted meshes, the construction of a surrogate interface requires (i) classification of mesh entities with respect to the underlying geometry, (ii) assignment of material regions (grains), and (iii) evaluation of geometric quantities required by the SCZM formulation. Together, these steps define a consistent pipeline for constructing surrogate interfaces and associated geometric quantities in polycrystalline RVEs, enabling the application of SCZM on non-interface-fitted meshes.

\subsection{Point classification}

\begin{figure}[htb]
  \centering
  \begin{subfigure}[b]{0.45\textwidth}
    \centering
    \includegraphics[width=\textwidth,trim={0 0 0 0},clip]{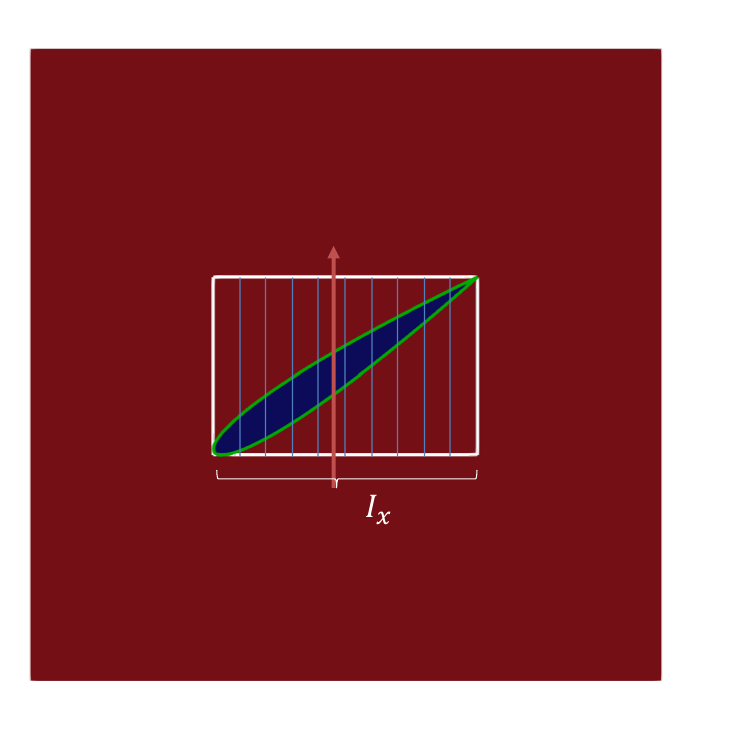}
    \caption{Control-cell method.
             A fixed ray direction is employed; the domain is partitioned into
             user-defined control cells, each linked to the \SBMElem{}s it contains.
             When a ray is cast, its control cell is located and only the mapped
             \SBMElem{}s are checked for intersection. Here, \SBMElem{} refers to elements belonging
             to the true boundary representation.}
    \label{fig:cell_approach}
  \end{subfigure}
  \hfill
  \begin{subfigure}[b]{0.45\textwidth}
    \centering
    \includegraphics[width=\textwidth,trim={0 0 0 0},clip]{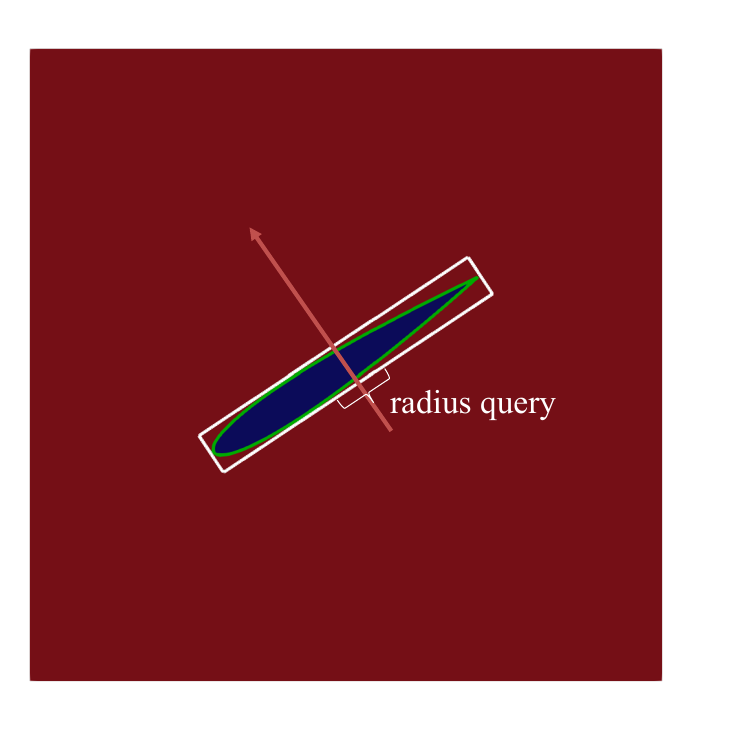}
    \caption{Proposed PCA-enhanced method.
             The ray direction is aligned with the least-variance principal
             component, and no spatial binning is required.
             Candidate \SBMElem{}s are obtained from a KD-tree radius query and
             then subjected to ray-element intersection tests.}
    \label{fig:pca_approach}
  \end{subfigure}

  \caption{Comparison of element-selection strategies for point-in-polygon
           classification.
           (a) Control-cell method~\cite{borazjani2008curvilinear}.
           (b) Present PCA-enhanced, ray-adaptive method.
           The proposed approach eliminates user-defined grid parameters and
           reduces the number of ray-element intersection checks.}
  \label{fig:inout_methods}
\end{figure}

Point classification is used to determine the relation between mesh entities and the underlying microstructure. Given a point $\bs{x} \in \mathbb{R}^d$, we determine whether it lies inside a given grain domain or in the exterior by means of ray casting. Specifically, a semi-infinite ray $\bs{r}(t) = \bs{x} + t \bs{\ell}$ is cast along a direction $\bs{\ell}$, and the number of intersections with the grain boundary representation $\mathcal{T}$ is evaluated. The point is classified based on the parity of the intersection count. 
To improve robustness and efficiency, the ray direction is chosen adaptively using principal component analysis (PCA) of the boundary geometry, rather than using a fixed direction as in the control cell method~\cite{borazjani2008curvilinear}. Let $\{\lambda_i, \bs{e}_i\}_{i=1}^d$ denote the eigenpairs of the covariance matrix of the boundary point cloud. The ray direction is selected as $\bs{\ell} = \bs{e}_d$, corresponding to the direction of minimum variance. This choice reduces the expected number of ray--facet intersections. \figref{fig:inout_methods} illustrates the difference between this approach and the control cell method procedure employed by, e.g.,~\citet{borazjani2008curvilinear}. \ref{app:inout_algorithms} provides implementation details.

For each query point, intersection tests are restricted to a dynamically constructed set of candidate facets $\mathcal{T}_{\text{cand}}(\bs{x}) \subset \mathcal{T}$, obtained via a geometric filtering criterion (e.g., bounding radius or spatial search structure). This reduces the computational complexity from $\mathcal{O}(N_{\text{elem}} N_{\mathcal{T}})$ to $\mathcal{O}(N_{\text{elem}} N_{\text{cand}})$, where $N_{\mathcal{T}}$ denotes the total number of facets in the true boundary representation, and $N_{\text{cand}} \ll N_{\mathcal{T}}$ in typical cases.

\subsection{Grain assignment and surrogate RVE construction}

Once point classification is available, each non-interface-fitted element $K \in \mathcal{T}_h$ is assigned a grain label. Rather than constructing conformal subcells, we adopt a dominant-volume criterion: each element is assigned the ID of the grain occupying the largest fraction of its volume. 

Formally, let $\chi_g(\bs{x})$ denote the indicator function of grain $g$. The assigned grain for element $K$ is
\begin{equation}
    g_K = \arg\max_g \int_K \chi_g(\bs{x}) \, dV.
\end{equation}
This construction yields a surrogate RVE in which each element is associated with a single grain, while the grain boundaries are implicitly represented by the collection of inter-element facets separating elements with different grain IDs.

A detailed description of the algorithm is provided in~\ref{app:rve_generation}.

\subsection{Distance function and normal evaluation.}

The SCZM formulation requires the distance vector and normal direction associated with the true interface. For each surrogate quadrature point $\bs{q}_h$, we compute its closest-point projection onto the true interface. 

To this end, we construct a k-d tree using the centroids of the interface facets in $\mathcal{T}$. For a given $\bs{q}_h$, the nearest interface element is identified via a nearest-neighbor query. The distance function is then defined as
\begin{equation}
    d(\bs{q}_h) = \min_{\bs{y} \in \mathcal{T}} \|\bs{q}_h - \bs{y}\|,
\end{equation}
and the associated normal vector is taken as the normal of the nearest interface element. Note the distance vector $\bs{d}$ obtained here is used in the definition of the shifted displacement jump in \secref{sec:TSL}. This procedure follows the approach in~\cite{yang2024optimal}.

\subsection{Generation of an interface-fitted mesh and solution projection}
\label{sec:SCZM2IFM}

To transfer the solution from the non-interface-fitted SCZM mesh to an interface-fitted mesh, we employ two main algorithmic components. The target IFM is constructed directly from the SCZM mesh, without assuming an \emph{a priori} interface-fitted discretization.

First, we construct the interface-fitted mesh by conformalizing the SCZM mesh along the interfaces. In particular, we perform only partial mesh surgery: most SCZM elements are copied directly to the IFM, except for elements intersected by the interface. Each intersected element is repartitioned and triangulated so that the resulting mesh conforms to the material interface. These operations are widely supported in utility packages such as \textsc{\href{https://github.com/shapely/shapely}{Shapely}}, \textsc{\href{https://pypi.org/project/manifold3d/}{manifold3d}}, and \textsc{\href{https://pypi.org/project/tetgen/}{TetGen}}.

Second, we project the SCZM solution onto the IFM. The SCZM mesh is preprocessed by constructing a centroid-based k-d tree and, when applicable, a direct lookup structure for structured uniform meshes. If the target point lies inside the source element belonging to the same sub-domain, standard finite-element interpolation is used. Otherwise, a recovery procedure based on the closest element is employed using the same surface--transport expansion as in \eqnref{eq:surftrans}.

\begin{figure}[htb]
  \centering
  \includegraphics[width=0.85\textwidth]{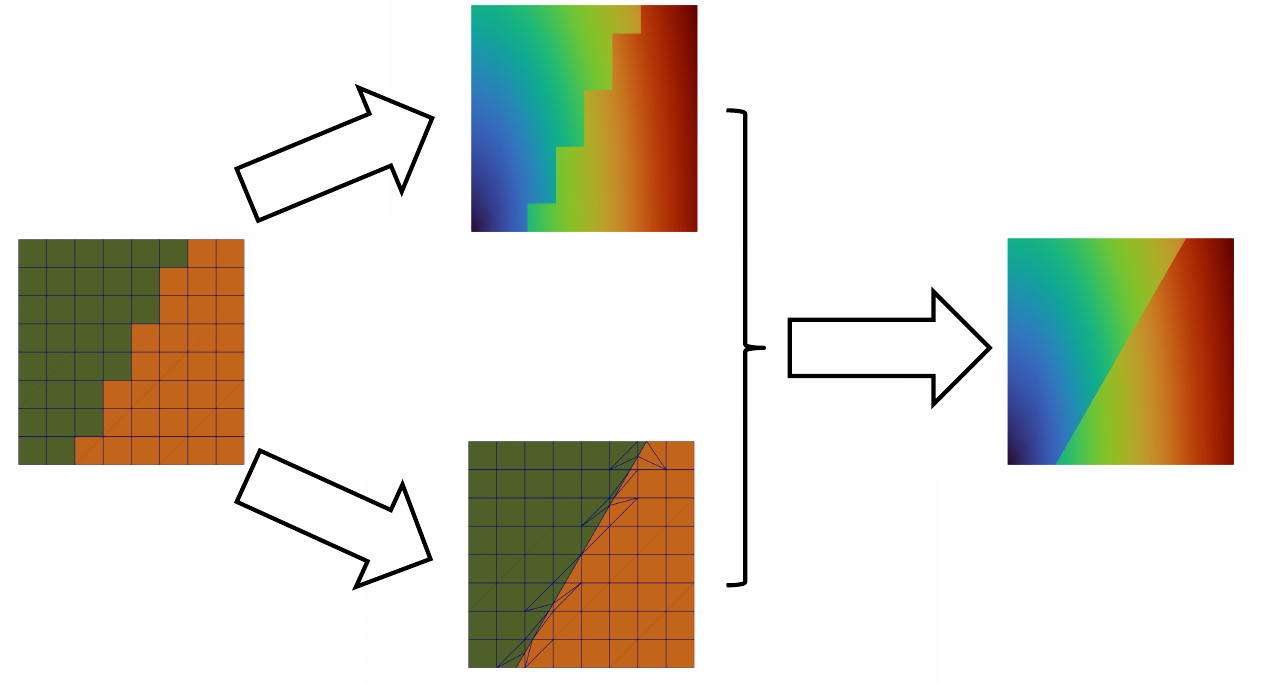}
  \caption{
    Schematic overview of the SCZM-to-IFM transfer procedure. The SCZM mesh is first locally conformalized near the interface to generate an interface-fitted mesh. Since the solution is already available on the original non-interface-fitted SCZM mesh, it is then projected onto the generated IFM mesh for interface-aligned visualization.
    }
  \label{fig:sczm_to_ifm_overview}
\end{figure}

\figref{fig:sczm_to_ifm_overview} provides a schematic overview of this procedure. Detailed algorithmic descriptions and supporting visualizations are deferred to \ref{app:sczm_to_ifm_algorithms}, where they collectively validate both the geometric reconstruction and the solution transfer procedure.

\section{Numerical results}
\label{sec: Results}

In this section, we assess the accuracy, robustness, and generality of the proposed SCZM framework through a series of numerical studies. The validation is organized in a progressive manner, beginning with a manufactured solution (MMS) to verify the correctness and convergence properties of the formulation, followed by single-interface benchmark problems to evaluate the enforcement of traction-separation laws on non-interface-fitted meshes. We then consider polycrystalline representative volume elements (RVEs) with increasing complexity, including both linear elasticity and history-dependent crystal plasticity models. Note that the unit system used in the numerical examples is $\left( \mathrm{N}, \mathrm{mm}, \mathrm{s} \right)$, unless otherwise stated.

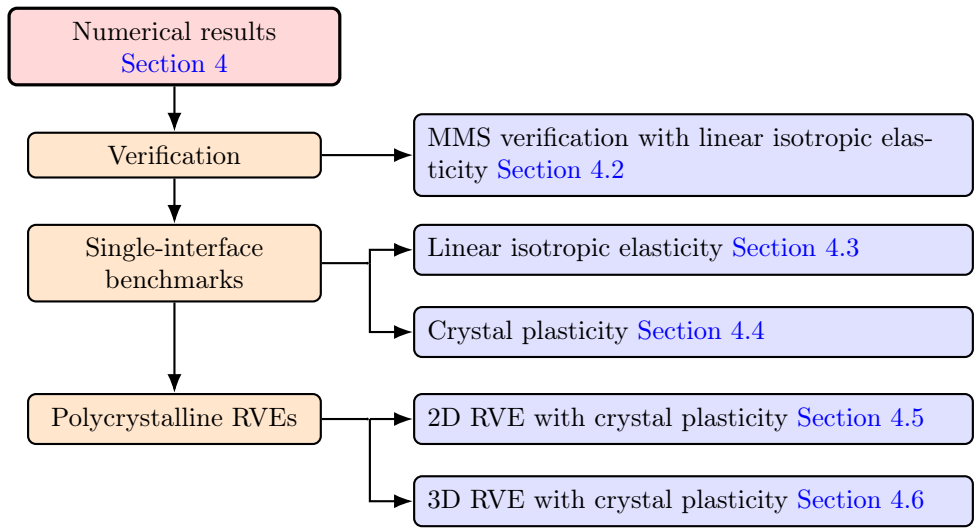
\begin{figure}[htb]
    \centering
    \begin{tikzpicture}[
        >=Latex,
        every node/.style={font=\footnotesize},
        box/.style={
            rectangle,
            rounded corners=3pt,
            draw=black,
            thick,
            align=center,
            inner sep=5pt
        },
        root/.style={
            box,
            fill=red!15,
            very thick,
            text width=4cm
        },
        category/.style={
            box,
            fill=orange!20,
            text width=3.5cm 
        },
        subcat/.style={
            box,
            fill=blue!12,
            text width=7cm,   
            align=left        
        },
        line/.style={draw, thick, -{Latex[length=2.5mm]}},
        dot/.style={circle, fill=black, inner sep=0pt, minimum size=1pt}
    ]


    \node[root] (root) {Numerical results\\\secref{sec: Results}};

    \node[category, below=6mm of root] (verif) {Verification};

    \node[category, below=6mm of verif] (single) {Single-interface\\benchmarks};

    \node[category, below=12mm of single] (poly) {Polycrystalline RVEs};


    \node[subcat, right=12mm of verif] (mms)
    {MMS verification with linear isotropic elasticity \secref{sec:sczm_mms}};

    \node[subcat, right=12mm of single.north east, anchor=north west] (cp)
    {Linear isotropic elasticity \secref{sec:sczm_isotropic}};

    \node[subcat, below=4mm of cp] (neml_single)
    {Crystal plasticity \secref{sec:single_neml2_2D}};

    \node[subcat, right=12mm of poly.north east, anchor=north west] (poly2d)
    {2D RVE with crystal plasticity \secref{sec:multi_neml2_2D}};

    \node[subcat, below=4mm of poly2d] (poly3d)
    {3D RVE with crystal plasticity \secref{sec:multi_neml2_3D}};


    \draw[line] (root.south) -- (verif.north);
    \draw[line] (verif.south) -- (single.north);
    \draw[line] (single.south) -- (poly.north);

    \draw[line] (verif.east) -- (mms.west);

    \draw[thick] (single.east) -- ++(0.6cm, 0) node[dot] (d1) {};
    \draw[line] (d1.east) |- (cp.west);
    \draw[line] (d1.east) |- (neml_single.west);

    \draw[thick] (poly.east) -- ++(0.6cm, 0) node[dot] (d2) {};
    \draw[line] (d2.east) |- (poly2d.west);
    \draw[line] (d2.east) |- (poly3d.west);

    \end{tikzpicture}
    \caption{Tree diagram summarizing the numerical studies presented in \secref{sec: Results}.}
    \label{fig:sczm_tree_optimized}
\end{figure}

The test cases span multiple spatial dimensions (2D and 3D), element types (triangular, quadrilateral, and hexahedral), and traction-separation laws, allowing for a comprehensive evaluation of the method across a broad range of settings. A schematic overview of the numerical studies is provided in \figref{fig:sczm_tree_optimized}.

\subsection{Algorithmic performance of point classification}
\label{sec:wall_time}

We assess the efficiency of the proposed PCA-enhanced point classification algorithm by comparing it with a brute-force full-intersection-scan ray-casting approach~\cite{haines1994point}, in which each query point is tested against all boundary elements.

The test is conducted on a two-dimensional domain containing a NACA0012 airfoil geometry, with the chord length non-dimensionalized to unity. The 2D domain is discretized using uniform quadrilateral elements with mesh size $h = 1/256$. To evaluate the effect of geometric complexity, the number of boundary segments is varied from $N_{\mathcal{T}} = 100$ to $1600$. Reported runtimes include both surrogate interface construction and point classification, averaged over five runs.

\begin{figure}[htb]
  \centering
  \begin{subfigure}[b]{0.48\textwidth}
    \centering
    \begin{tikzpicture}
      \begin{axis}[
        width=\textwidth,
        height=0.76\textwidth,
        ymode=log,
        xmin=-0.2, xmax=4.2,
        ymin=0.1, ymax=2.3,
        xtick={0,1,2,3,4},
        xticklabels={100,200,400,800,1600},
        ytick={0.1,0.21,0.46,0.97,2.1},
        yticklabels={0.1,0.21,0.46,0.97,2.1},
        xlabel={$N_{\mathcal{T}}$},
        ylabel={Runtime ($s$)},
        grid=both,
        grid style={dashed, opacity=0.55},
        tick align=outside,
        tick label style={font=\normalsize},
        label style={font=\normalsize},
        legend style={
          font=\fontsize{12}{14}\selectfont,
          at={(-0.3,1.38)},
          anchor=north west,
          draw=black,
          fill=white,
          cells={align=left},
          inner xsep=4pt,
          inner ysep=2pt
        },
      ]

      \addplot[
        color=mplblue,
        mark=square*,
        line width=1.6pt,
        mark size=2.4pt,
        error bars/.cd,
          y dir=both,
          y explicit,
      ]
      table[
        x expr=\coordindex,
        y=pca_method,
        y error=pca_method_err,
        col sep=space
      ]{runtime_data_2D.txt};
      \addlegendentry{PCA-enhanced method}

      \addplot[
        color=mplorange,
        mark=square*,
        line width=1.4pt,
        mark size=2.2pt,
        error bars/.cd,
          y dir=both,
          y explicit,
      ]
      table[
        x expr=\coordindex,
        y=brute_force,
        y error=brute_force_err,
        col sep=space
      ]{runtime_data_2D.txt};
      \addlegendentry{Full-intersection-scan approach~\cite{haines1994point}}

      \end{axis}
    \end{tikzpicture}
    \caption{Runtime vs. number of boundary elements.}
    \label{fig:runtime}
  \end{subfigure}
  \hfill
  \begin{subfigure}[b]{0.48\textwidth}
    \centering
    \begin{tikzpicture}
      \begin{axis}[
        width=\textwidth,
        height=0.76\textwidth,
        xmin=-0.2, xmax=4.2,
        ymin=1.6, ymax=10.5,
        xtick={0,1,2,3,4},
        xticklabels={100,200,400,800,1600},
        xlabel={$N_{\mathcal{T}}$},
        ylabel={Speedup},
        grid=major,
        grid style={dashed, opacity=0.55},
        tick align=outside,
        tick label style={font=\normalsize},
        label style={font=\normalsize},        
        ]

      \addplot[
        color=mplblue,
        mark=diamond*,
        line width=1.8pt,
        mark size=2.6pt
      ]
      table[
        x expr=\coordindex,
        y=speedup,
        col sep=space
      ]{speedup_data_2D.txt};

      \end{axis}
    \end{tikzpicture}
    \caption{Speedup achieved by avoiding brute force.}
    \label{fig:speedup}
  \end{subfigure}
  \caption{Performance comparison between brute force and PCA-enhanced point-in-polygon testing (\secref{sec:wall_time}).}
  \label{fig:runtime_vs_speedup}
\end{figure}
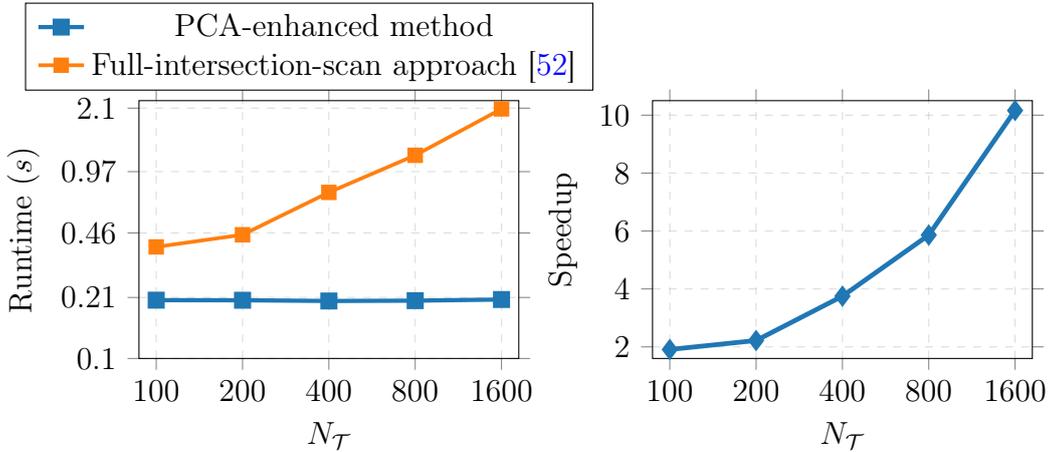

As shown in \figref{fig:runtime_vs_speedup}, the proposed method achieves significantly improved performance compared to the brute-force approach. While the brute-force runtime increases approximately linearly with $N_{\mathcal{T}}$, the PCA-enhanced method exhibits near-constant runtime, indicating that the number of candidate intersection checks remains effectively bounded. This behavior is consistent with the reduced complexity $\mathcal{O}(N_{\text{elem}} N_{\text{cand}})$, where $N_{\text{cand}} \ll N_{\mathcal{T}}$.

\begin{figure}[htb]
  \centering

\begin{subfigure}[b]{0.45\textwidth}
  \centering
  \begin{tikzpicture}
    \begin{axis}[
      width=\textwidth,
      height=0.76\textwidth,
      xlabel={Angle},
      ylabel={Runtime (s)},
      grid=both,
      grid style={dashed, opacity=0.5},
      tick align=outside,
      tick label style={font=\normalsize},
      label style={font=\normalsize},
      legend style={
        at={(-0.5,1.39)},
        anchor=north west,
        draw=black,
        fill=white
      },
    ]

    \addplot[
      color=mplblue,
      mark=diamond*,
      line width=1.6pt,
      mark size=2.3pt
    ]
    table[
      x index=0,
      y index=1,
    ]{runtime_auto_true.txt};
    \addlegendentry{PCA-enhanced (adaptive-direction)}

    \addplot[
      color=mplorange,
      mark=square*,
      line width=1.6pt,
      mark size=2.2pt
    ]
    table[
      x index=0,
      y index=1,
    ]{runtime_auto_false.txt};
    \addlegendentry{Fixed $y$-direction}

    \end{axis}
  \end{tikzpicture}
  \caption{Runtime comparison between PCA-enhanced (adaptive-direction) and fixed-direction methods.}
\end{subfigure}
  \hfill
  \begin{subfigure}[b]{0.45\textwidth}
    \centering
    \begin{tikzpicture}
      \begin{axis}[
        width=\textwidth,
        height=0.76\textwidth,
        xlabel={Angle},
        ylabel={Speedup (fixed / adaptive)},
        grid=both,
        grid style={dashed, opacity=0.5},
        tick align=outside,
        tick label style={font=\normalsize},
        label style={font=\normalsize},
      ]

      \addplot[
        color=mplblue,
        mark=diamond*,
        line width=1.8pt,
        mark size=2.4pt
      ]
      coordinates {
        (0,  1.047450)
        (5,  1.025334)
        (10, 1.014124)
        (15, 1.009085)
        (20, 1.047820)
        (25, 1.081119)
        (30, 1.087137)
        (35, 1.128966)
        (40, 1.171809)
        (45, 1.185444)
      };

      \end{axis}
    \end{tikzpicture}
    \caption{Speedup of PCA-enhanced (adaptive-direction) method over fixed $y$-direction method.}
  \end{subfigure}

  \caption{Performance evaluation under varying airfoil rotation angles (\secref{sec:wall_time}).}
  \label{fig:timing_rotation_study}
\end{figure}
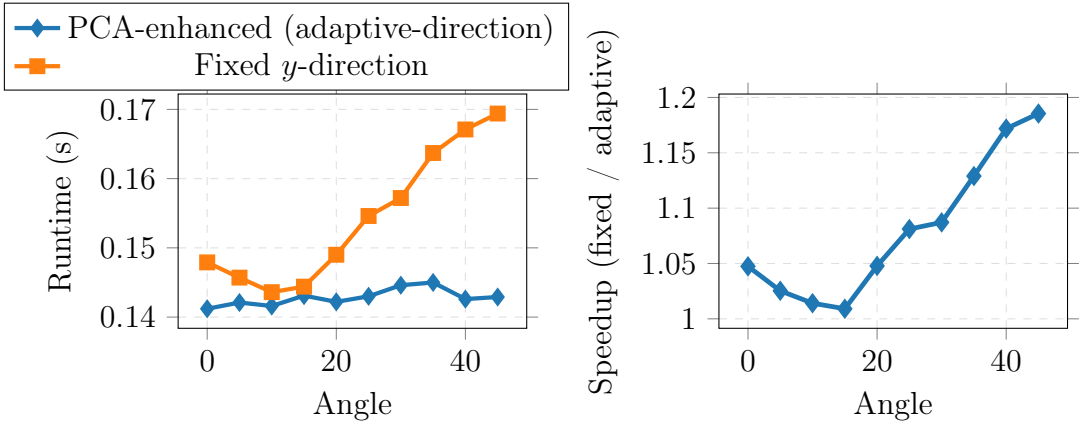

To further evaluate the effect of ray direction, we compare the PCA-based adaptive approach with a fixed-direction method in which rays are cast along the $y$-axis. The airfoil is rotated from $0^\circ$ to $45^\circ$, and the corresponding runtimes are shown in \figref{fig:timing_rotation_study}. The PCA-based method consistently outperforms the fixed-direction approach across all orientations. The minimum runtime for the fixed-direction method occurs near $10^\circ$, where the geometry is nearly aligned with the $y$-axis, confirming that performance is sensitive to the choice of ray direction.

These results demonstrate that selecting the ray direction based on the principal components of the geometry reduces unnecessary ray--facet intersection tests and leads to substantial computational savings in point classification.

\subsection{MMS-based verification with linear isotropic elasticity}
\label{sec:sczm_mms}

We consider a two-dimensional domain $\Omega = [-0.5,0.5]^2$ with a planar interface located at $x = x_0 = 0.25$, separating subdomains $\Omega^+$ (left subdomain) and $\Omega^-$ (right subdomain). 
All quantities in this section are nondimensionalized for simplicity in constructing the manufactured solution. 
We define a unit normal vector $\boldsymbol{n} = (1,0)^\mathsf{T}$ pointing from $\Omega^+$ to $\Omega^-$. Accordingly, we set $\boldsymbol{n}_+ = \boldsymbol{n}$ and $\boldsymbol{n}_- = -\boldsymbol{n}$.

A manufactured displacement field is prescribed in $\Omega^+$ as $\bs{u}^+ = (-\sin(\pi x),\,0)^\mathsf{T}$. The material response is linear isotropic elasticity under plane strain with Poisson's ratio $\nu = 0$, and Young's moduli $E^+ = 0.1$ and $E^- = 1$ in $\Omega^+$ and $\Omega^-$, respectively.

The displacement field in $\Omega^-$ is constructed using a prescribed jump $\bs{g}(x) = (ax^2 + b,\,0)^\mathsf{T}$ such that $\bs{u}^- = \bs{u}^+ - \bs{g}$. The corresponding stress fields follow from the constitutive law, yielding $\sigma_{xx}^+ = -\pi E^+ \cos(\pi x)$ in $\Omega^+$ and $\sigma_{xx}^- = E^-(-\pi \cos(\pi x) - 2ax)$ in $\Omega^-$. To ensure consistency, traction continuity is enforced at the interface, i.e., $\bs{\sigma}^+ \bs{n}_+ + \bs{\sigma}^- \bs{n}_- = \bs{0}$ at $x = x_0$. In addition, we prescribe a linear cohesive relation under the present nondimensional setting as $\bs{t}_{\text{coh}} = -\bs{g}$, and enforce $\bs{t}_{\text{coh}} = \bs{\sigma}^+ \bs{n}_+ = -\bs{\sigma}^- \bs{n}_-$ at the interface. These conditions uniquely determine the coefficients $a$ and $b$. The body force is obtained by substituting the manufactured solution into the governing equations. This construction yields an exact solution satisfying both the bulk equations and the interface law, providing a suitable benchmark for verification.

To assess the convergence behavior of SCZM, we perform simulations on Delaunay triangular meshes with varying resolution. The characteristic mesh size is defined as $h_\Omega = (|\Omega|/N_e)^{1/d}$ with $d=2$.

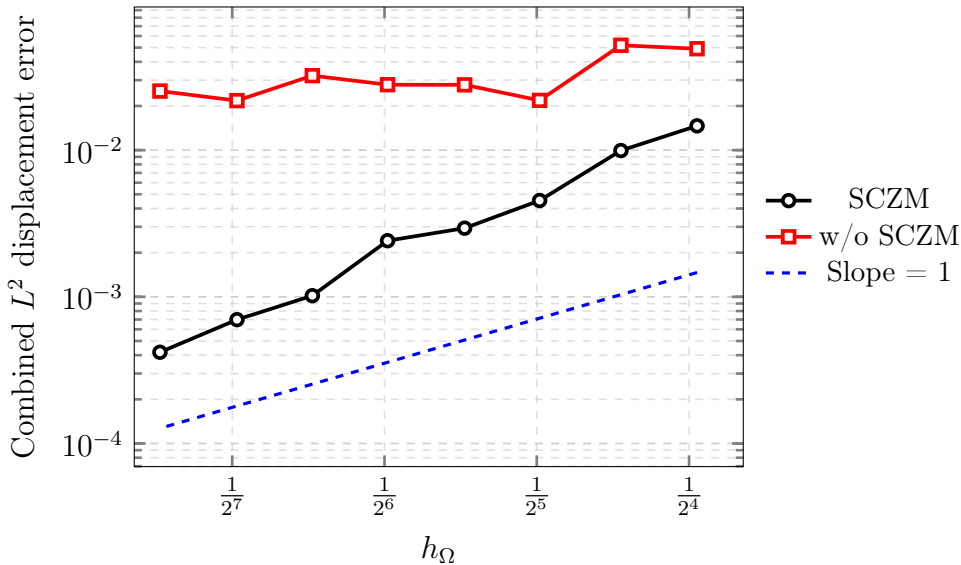
\begin{figure}[htb]
  \centering
  \begin{tikzpicture}
    \begin{axis}[
      width=0.7\textwidth,
      height=0.56\textwidth,
      xmode=log,
      ymode=log,
      xlabel={$h_\Omega$},
      ylabel={Combined $L^2$ displacement error},
      xmin=5e-3, xmax=8e-2,
      grid=both,
      grid style={dashed, line width=0.6pt, opacity=0.5},
      tick align=inside,
      tick style={line width=1pt},
      tick label style={font=\normalsize},
      label style={font=\normalsize},
      legend style={
        font=\small,
        at={(1.02,0.5)},
        anchor=west,
        draw=none,
        fill=none,
      },
      xtick={1/2^3,1/2^4,1/2^5,1/2^6,1/2^7,1/2^8},
      xticklabels={$ \frac{1}{2^{3}} $,$ \frac{1}{2^{4}} $,$ \frac{1}{2^{5}} $,$ \frac{1}{2^{6}} $,$ \frac{1}{2^{7}} $,$ \frac{1}{2^{8}} $},
    ]

      \addplot[
        color=black,
        mark=*,
        mark options={fill=white},
        line width=1.4pt,
        mark size=2.2pt
      ]
      table[x=h,y=errorsSCZM,col sep=space]{SCZM_solid_mechanics_data.txt};
      \addlegendentry{SCZM}

      \addplot[
        color=red,
        mark=square*,
        mark options={fill=white},
        line width=1.4pt,
        mark size=2.2pt
      ]
      table[x=h,y=errorsNo_SCZM,col sep=space]{SCZM_solid_mechanics_data.txt};
      \addlegendentry{w/o SCZM}

      \addplot[
        color=blue,
        dashed,
        line width=1.2pt
      ]
      table[x=h,y=slope1,col sep=space]{SCZM_solid_mechanics_data.txt};
      \addlegendentry{Slope = 1}

    \end{axis}
  \end{tikzpicture}
    \caption{Mesh convergence study (\secref{sec:sczm_mms}): the combined $L^2$ displacement error, computed from the $x$- and $y$-direction displacement errors, decreases linearly with mesh refinement for the SBM-SCZM formulation, demonstrating first-order accuracy. In contrast, without SCZM, the error remains $\mathcal{O}(1)$.}  
    \label{fig:mesh_convergence_sczm}
\end{figure}

The convergence results are shown in \figref{fig:mesh_convergence_sczm}, where the combined $L^2$ displacement error (including both $x$- and $y$-components) is plotted as a function of the characteristic mesh size $h_\Omega$. The SCZM formulation exhibits clear first-order convergence with mesh refinement. This behavior is consistent with the underlying approximation, in which the dominant error arises from the geometric mismatch between the true interface and its surrogate representation. In particular, the use of a facet-based surrogate interface yields a first-order accurate approximation of the interface geometry, which in turn limits the overall convergence rate, even though the bulk discretization employs standard finite elements. The results confirm that SCZM correctly captures both bulk and interface contributions and achieves the expected convergence behavior for non-interface-fitted formulations. For comparison, simulations without SCZM fail to exhibit consistent convergence, highlighting the importance of properly accounting for the shifted interface contribution.

Next, we consider a simplified test case by prescribing $\bs{g}(x)=(ax+b,\,0)^\mathsf{T}$ and setting $E^+=E^-=0.1$. The cohesive traction is then $\bs{t}_{\text{coh}}=-\bs{g}=-(ax+b,\,0)^\mathsf{T}$. For this parameter choice, $E^+=E^-$ implies $a=0$, so $\bs{t}_{\text{coh}}=-(b,\,0)^\mathsf{T}$ is spatially constant and $\nabla \bs{t}_{\text{coh}}=\bs{0}$.

As shown in \figref{fig:mesh_convergence_sczm_quad}, this simplified test case exhibits a convergence rate close to second order. This behavior can be explained by examining the role of the \textit{field correction} term. In this simplified case, given that $\nabla_{\bs{n}}\bs{t}_{\text{coh}}=\bs{0}$, the \textit{field correction} term is close to zero. Consequently, omitting the \textit{field correction} term does not result in a loss of accuracy; the system maintains the second-order convergence expected when using first-order finite element basis functions.



\pgfplotstableread[col sep=space]{errorsNo_SCZM.txt}\sczmdata

\pgfplotstablegetelem{0}{h}\of{\sczmdata}
\pgfmathsetmacro{\hzero}{\pgfplotsretval}

\pgfplotstablegetelem{0}{errorsSCZM}\of{\sczmdata}
\pgfmathsetmacro{\ezero}{\pgfplotsretval}

\pgfmathsetmacro{\Ctwo}{\ezero/(\hzero*\hzero)}

\begin{figure}[htb]
  \centering
  \begin{tikzpicture}
    \begin{axis}[
      width=0.7\textwidth,
      height=0.56\textwidth,
      xmode=log,
      ymode=log,
      xlabel={$h_\Omega$},
      ylabel={$L^2$ error},
      xmin=5e-3, xmax=8e-2,
      grid=both,
      grid style={dashed, line width=0.6pt, opacity=0.5},
      tick align=inside,
      tick style={line width=1pt},
      tick label style={font=\normalsize},
      label style={font=\normalsize},
      legend style={
        font=\small,
        at={(1.02,0.5)},
        anchor=west,
        draw=none,
        fill=none,
      },
      xtick={1/2^3,1/2^4,1/2^5,1/2^6,1/2^7,1/2^8},
      xticklabels={$ \frac{1}{2^{3}} $,$ \frac{1}{2^{4}} $,$ \frac{1}{2^{5}} $,$ \frac{1}{2^{6}} $,$ \frac{1}{2^{7}} $,$ \frac{1}{2^{8}} $},
    ]

      \addplot[
        color=black,
        mark=*,
        mark options={fill=white},
        line width=1.4pt,
        mark size=2.2pt
      ]
      table[x=h,y=errorsSCZM,col sep=space]{errorsNo_SCZM.txt};
      \addlegendentry{SCZM}

      \addplot[
        color=red,
        mark=square*,
        mark options={fill=white},
        line width=1.4pt,
        mark size=2.2pt
      ]
      table[x=h,y=errorsNo_SCZM,col sep=space]{errorsNo_SCZM.txt};
      \addlegendentry{w/o SCZM}

        \addplot[
          color=blue,
          dashed,
          line width=1.2pt
        ]
        table[
          x=h,
          y expr=\Ctwo*(\thisrow{h})^2,
          col sep=space
        ]{\sczmdata};
        \addlegendentry{Slope = 2}    
        \end{axis}
  \end{tikzpicture}
    \caption{Mesh convergence study (\secref{sec:sczm_mms}): the combined $L^2$ displacement error -- calculated from the $x$- and $y$-components -- decreases quadratically with mesh refinement, demonstrating nearly second-order accuracy. In contrast, without SCZM, the error remains approximately $\mathcal{O}(1)$.}
\label{fig:mesh_convergence_sczm_quad}
\end{figure}
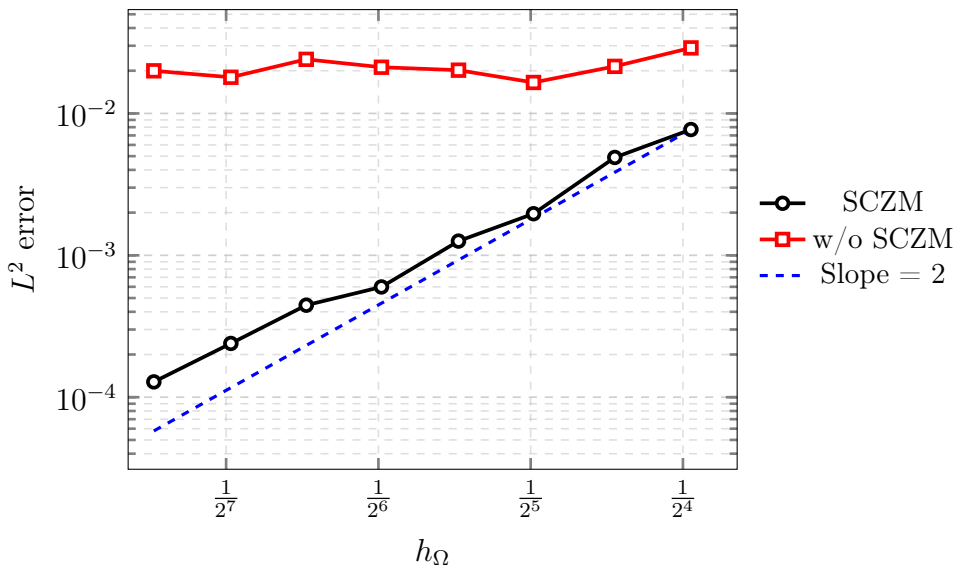

\subsection{Single-interface benchmark with isotropic elasticity}
\label{sec:sczm_isotropic}

We consider a two-dimensional domain $[0,1]^2$ with a single interface located at $x=0.5$. The material is homogeneous on both sides, with Young's modulus $E = 10^{3}$ and Poisson's ratio $\nu = 0.3$. An exponential traction-separation law is prescribed along the interface, with parameters $\mathcal{G}_c = 50$, $\delta_0 = 0.1$, and $\beta = 0$. Boundary conditions enforce $u_x=0$ on the left boundary, $u_y=0$ on the bottom boundary, and a displacement-controlled loading $u_x = 10^{-2} t$ on the right boundary. Simulations are performed up to $t=200$.

As a reference solution, we employ an interface-fitted discretization using a uniform quadrilateral mesh with resolution $1/256$. SCZM results are obtained on non-interface-fitted Delaunay triangular meshes with mesh size $h_{\Omega} \approx 4.58 \times 10^{-2}$, for which the surrogate interface is only an approximation of the true interface.

\begin{figure}[htb]
  \centering
  \begin{tikzpicture}
    \begin{axis}[
      width=0.62\textwidth,
      height=0.52\textwidth,
      xmode=log,
      ymode=log,
      xlabel={Load increment size ($\Delta t$)},
      ylabel={$|\text{energy release} - 50|$},
    grid=major,
    major grid style={
      line width=0.25pt,
      draw=gray!25
    },
    minor grid style={draw=none},
    tick align=outside,
      tick label style={font=\normalsize},
      label style={font=\normalsize},
      legend style={
        font=\small,
        at={(1.02,0.7)},
        anchor=south west,
        draw=none,
        fill=none
      },
    ]

      \addplot[
        color=black,
        mark=*,
        mark options={fill=white},
        line width=1.4pt,
        mark size=2.2pt
      ]
      table[
        x=dt,
        y=err_sczm,
        col sep=space
      ]{error_vs_dt_compare_min_a_data.txt};
      \addlegendentry{SCZM}

      \addplot[
        color=red,
        mark=square*,
        mark options={fill=white},
        line width=1.4pt,
        mark size=2.2pt
      ]
      table[
        x=dt,
        y=err_no_sczm,
        col sep=space
      ]{error_vs_dt_compare_min_a_data.txt};
      \addlegendentry{w/o SCZM}

      \addplot[
        color=blue,
        dashed,
        line width=1.2pt
      ]
      table[
        x=dt,
        y=slope2,
        col sep=space
      ]{error_vs_dt_compare_min_a_data.txt};
      \addlegendentry{slope = 2}

    \end{axis}
  \end{tikzpicture}
    \caption{Convergence of the energy release error with respect to the load increment size (\secref{sec:sczm_isotropic}). The SCZM case exhibits second-order convergence, whereas the error without SCZM remains essentially constant (i.e., $\mathcal{O}(1)$).}    
    \label{fig:time_convergence_SCZM}
\end{figure}
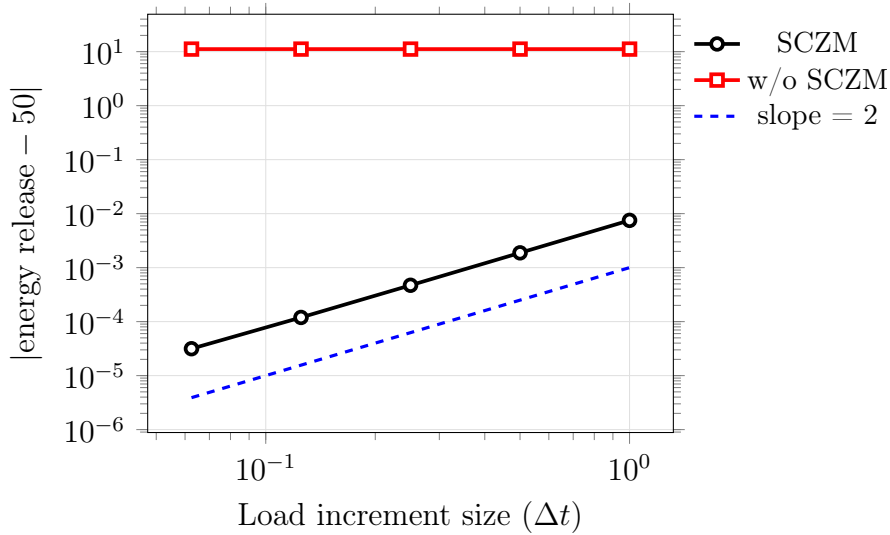

\figref{fig:reactive_effective} compares the reactive force--displacement response. SCZM exhibits excellent agreement with the interface-fitted solution, whereas simulations without SCZM show significant deviation. This demonstrates that the shifted formulation correctly enforces the cohesive response on non-interface-fitted meshes. To further assess accuracy, we examine the convergence of energy release. The energy release is computed as $\int_0^\infty \bs{t}_{\text{coh}}(\jumpscalar{\bs{u}})\cdot\jumpscalar{\bs{u}}$, with the SCZM formulation using the shifted jump. The exact energy release is $\mathcal{G}_c L = 50$, and we report the error $|\text{energy release} - 50|$. \figref{fig:time_convergence_SCZM} shows that SCZM achieves second-order convergence with respect to load increment size, whereas the error without SCZM remains essentially constant. This highlights the importance of consistently accounting for the shifted interface kinematics.

\begin{figure}[htb]
\centering

\begin{subfigure}[b]{0.47\textwidth}
\centering
\begin{tikzpicture}
\begin{axis}[
    width=\textwidth,
    height=0.78\textwidth,
    xlabel={$\delta_{\mathrm{eff}}$},
    ylabel={$R_x$},
    grid=major,
    major grid style={line width=0.25pt, draw=gray!25},
    minor grid style={draw=none},
    tick align=outside,
    tick label style={font=\normalsize},
    label style={font=\normalsize},
    legend to name=sharedlegend_reactive,
    legend columns=2,
    legend style={
        draw=none,
        fill=none,
        font=\footnotesize,
        column sep=1em
    }
]

\addplot[
    color=black,
    line width=1.2pt
]
table[
    x=delta_eff,
    y=react_x,
    col sep=space,
    restrict expr to domain={\thisrow{curve_id}}{1:1}
]{all_area_curves_no_sczm_only_data.txt};
\addlegendentry{$h_{\Omega}=4.583492\times10^{-2}$}

\addplot[
    color=blue,
    line width=1.2pt
]
table[
    x=delta_eff,
    y=react_x,
    col sep=space,
    restrict expr to domain={\thisrow{curve_id}}{2:2}
]{all_area_curves_no_sczm_only_data.txt};
\addlegendentry{$h_{\Omega}=6.482037\times10^{-2}$}

\addplot[
    color=teal,
    line width=1.2pt
]
table[
    x=delta_eff,
    y=react_x,
    col sep=space,
    restrict expr to domain={\thisrow{curve_id}}{3:3}
]{all_area_curves_no_sczm_only_data.txt};
\addlegendentry{$h_{\Omega}=8.770580\times10^{-2}$}

\addplot[
    color=orange,
    line width=1.2pt
]
table[
    x=delta_eff,
    y=react_x,
    col sep=space,
    restrict expr to domain={\thisrow{curve_id}}{4:4}
]{all_area_curves_no_sczm_only_data.txt};
\addlegendentry{$h_{\Omega}=1.270001\times10^{-1}$}

\addplot[
    color=violet,
    line width=1.2pt
]
table[
    x=delta_eff,
    y=react_x,
    col sep=space,
    restrict expr to domain={\thisrow{curve_id}}{5:5}
]{all_area_curves_no_sczm_only_data.txt};
\addlegendentry{$h_{\Omega}=1.825742\times10^{-1}$}

\addplot[
    color=brown,
    line width=1.2pt
]
table[
    x=delta_eff,
    y=react_x,
    col sep=space,
    restrict expr to domain={\thisrow{curve_id}}{6:6}
]{all_area_curves_no_sczm_only_data.txt};
\addlegendentry{$h_{\Omega}=2.500000\times10^{-1}$}

\addplot[
    color=red,
    dashed,
    line width=1.4pt
]
table[
    x=delta_eff,
    y=react_x,
    col sep=space
]{all_area_curves_ifm_data.txt};
\addlegendentry{IFM}

\end{axis}
\end{tikzpicture}
\caption{Without SCZM}
\end{subfigure}
\hfill
\begin{subfigure}[b]{0.47\textwidth}
\centering
\begin{tikzpicture}
\begin{axis}[
    width=\textwidth,
    height=0.78\textwidth,
    xlabel={$\delta_{\mathrm{eff}}$},
    ylabel={$R_x$},
    grid=major,
    major grid style={line width=0.25pt, draw=gray!25},
    minor grid style={draw=none},
    tick align=outside,
    tick label style={font=\normalsize},
    label style={font=\normalsize},
    legend style={draw=none}
    ]

\addplot[
    color=black,
    line width=1.2pt
]
table[
    x=delta_eff,
    y=react_x,
    col sep=space,
    restrict expr to domain={\thisrow{curve_id}}{1:1}
]{all_area_curves_sczm_only_data.txt};

\addplot[
    color=blue,
    line width=1.2pt
]
table[
    x=delta_eff,
    y=react_x,
    col sep=space,
    restrict expr to domain={\thisrow{curve_id}}{2:2}
]{all_area_curves_sczm_only_data.txt};

\addplot[
    color=teal,
    line width=1.2pt
]
table[
    x=delta_eff,
    y=react_x,
    col sep=space,
    restrict expr to domain={\thisrow{curve_id}}{3:3}
]{all_area_curves_sczm_only_data.txt};

\addplot[
    color=orange,
    line width=1.2pt
]
table[
    x=delta_eff,
    y=react_x,
    col sep=space,
    restrict expr to domain={\thisrow{curve_id}}{4:4}
]{all_area_curves_sczm_only_data.txt};

\addplot[
    color=violet,
    line width=1.2pt
]
table[
    x=delta_eff,
    y=react_x,
    col sep=space,
    restrict expr to domain={\thisrow{curve_id}}{5:5}
]{all_area_curves_sczm_only_data.txt};

\addplot[
    color=brown,
    line width=1.2pt
]
table[
    x=delta_eff,
    y=react_x,
    col sep=space,
    restrict expr to domain={\thisrow{curve_id}}{6:6}
]{all_area_curves_sczm_only_data.txt};

\addplot[
    color=red,
    dashed,
    line width=1.4pt
]
table[
    x=delta_eff,
    y=react_x,
    col sep=space
]{all_area_curves_ifm_data.txt};

\end{axis}
\end{tikzpicture}
\caption{With SCZM}
\end{subfigure}

\vspace{0.6em}

\pgfplotslegendfromname{sharedlegend_reactive}

\caption{Evolution of the reactive force with respect to the effective displacement, comparing simulations with and without SCZM (\secref{sec:sczm_isotropic}).}
\label{fig:reactive_effective}
\end{figure}
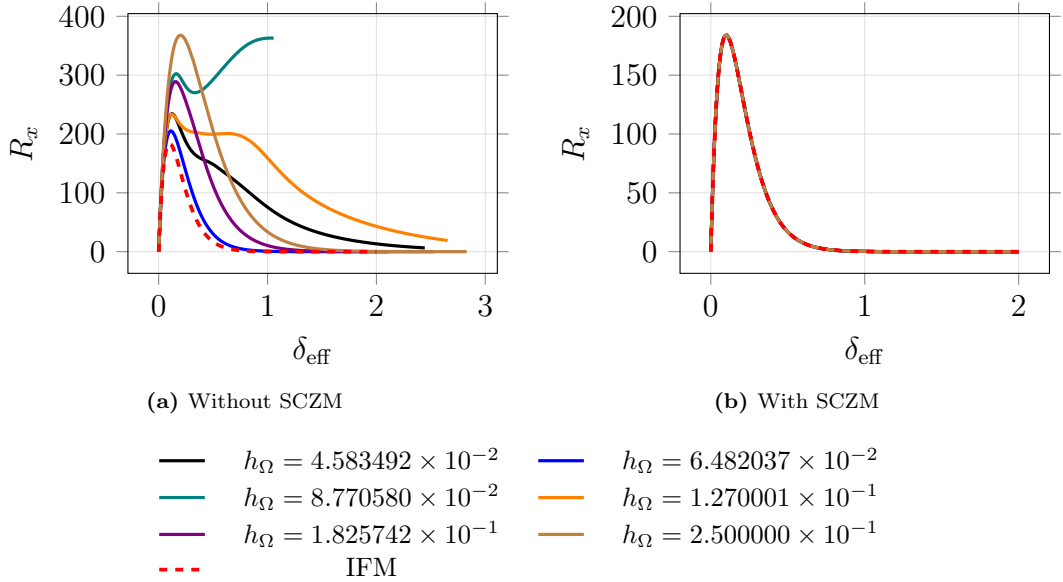

\begin{figure}[htb]
\centering

\begin{subfigure}[b]{0.45\textwidth}
\centering
\begin{tikzpicture}
\begin{axis}[
    width=\textwidth,
    height=0.75\textwidth,
    xlabel={$l$},
    ylabel={von Mises stress},
    ymode=log,
    ymin=200, ymax=270,
    grid=major,
    tick align=outside,
    label style={font=\normalsize},
    tick label style={font=\normalsize},
    legend to name=sharedlegend_stress,
    legend columns=2,
    legend style={draw=none, fill=none, font=\footnotesize},
    every axis plot/.append style={line width=1.5pt},
    line width=1.2pt
]

\addplot[color=black] table[x=l,y=vonmises_stress,col sep=space,
restrict expr to domain={\thisrow{curve_id}}{1:1}]
{stress_vs_l_theta_t30_all_data.txt};
\addlegendentry{IFM}

\addplot[color=orange] table[x=l,y=vonmises_stress,col sep=space,
restrict expr to domain={\thisrow{curve_id}}{2:2}]
{stress_vs_l_theta_t30_all_data.txt};
\addlegendentry{w/o SCZM}

\addplot[color=blue] table[x=l,y=vonmises_stress,col sep=space,
restrict expr to domain={\thisrow{curve_id}}{3:3}]
{stress_vs_l_theta_t30_all_data.txt};
\addlegendentry{SCZM (w/o directional correction term)}

\addplot[color=red,dashed] table[x=l,y=vonmises_stress,col sep=space,
restrict expr to domain={\thisrow{curve_id}}{4:4}]
{stress_vs_l_theta_t30_all_data.txt};
\addlegendentry{SCZM}

\end{axis}
\end{tikzpicture}
\caption{$t = 30$ (all models)}
\end{subfigure}
\hfill
\begin{subfigure}[b]{0.45\textwidth}
\centering
\begin{tikzpicture}
\begin{axis}[
    width=\textwidth,
    height=0.75\textwidth,
    xlabel={$l$},
    ylabel={von Mises stress},
    ymode=log,
    ymin=40, ymax=1000,
    grid=major,
    tick align=outside,
    label style={font=\normalsize},
    tick label style={font=\normalsize},
    line width=1.2pt
]

\addplot[color=black] table[x=l,y=vonmises_stress,col sep=space,
restrict expr to domain={\thisrow{curve_id}}{1:1}]
{stress_vs_l_theta_t50_all_data.txt};

\addplot[color=orange] table[x=l,y=vonmises_stress,col sep=space,
restrict expr to domain={\thisrow{curve_id}}{2:2}]
{stress_vs_l_theta_t50_all_data.txt};

\addplot[color=blue] table[x=l,y=vonmises_stress,col sep=space,
restrict expr to domain={\thisrow{curve_id}}{3:3}]
{stress_vs_l_theta_t50_all_data.txt};

\addplot[color=red,dashed] table[x=l,y=vonmises_stress,col sep=space,
restrict expr to domain={\thisrow{curve_id}}{4:4}]
{stress_vs_l_theta_t50_all_data.txt};

\end{axis}
\end{tikzpicture}
\caption{$t = 50$ (all models)}
\end{subfigure}

\vspace{0.8em}

\begin{subfigure}[b]{0.45\textwidth}
\centering
\begin{tikzpicture}
\begin{axis}[
    width=\textwidth,
    height=0.75\textwidth,
    xlabel={$l$},
    ylabel={von Mises stress},
    ymode=log,
    ymin=200, ymax=225,
    grid=major,
    tick align=outside,
    line width=1.2pt
]

\addplot[color=black] table[x=l,y=vonmises_stress,col sep=space,
restrict expr to domain={\thisrow{curve_id}}{1:1}]
{stress_vs_l_theta_t30_compare_data.txt};

\addplot[color=blue] table[x=l,y=vonmises_stress,col sep=space,
restrict expr to domain={\thisrow{curve_id}}{2:2}]
{stress_vs_l_theta_t30_compare_data.txt};

\addplot[color=red,dashed] table[x=l,y=vonmises_stress,col sep=space,
restrict expr to domain={\thisrow{curve_id}}{3:3}]
{stress_vs_l_theta_t30_compare_data.txt};

\end{axis}
\end{tikzpicture}
\caption{$t = 30$ (IFM vs SCZM variants)}
\end{subfigure}
\hfill
\begin{subfigure}[b]{0.45\textwidth}
\centering
\begin{tikzpicture}
\begin{axis}[
    width=\textwidth,
    height=0.75\textwidth,
    xlabel={$l$},
    ylabel={von Mises stress},
    ymode=log,
    ymin=43, ymax=49,
    grid=major,
    tick align=outside,
    line width=1.2pt
]

\addplot[color=black] table[x=l,y=vonmises_stress,col sep=space,
restrict expr to domain={\thisrow{curve_id}}{1:1}]
{stress_vs_l_theta_t50_compare_data.txt};

\addplot[color=blue] table[x=l,y=vonmises_stress,col sep=space,
restrict expr to domain={\thisrow{curve_id}}{2:2}]
{stress_vs_l_theta_t50_compare_data.txt};

\addplot[color=red,dashed] table[x=l,y=vonmises_stress,col sep=space,
restrict expr to domain={\thisrow{curve_id}}{3:3}]
{stress_vs_l_theta_t50_compare_data.txt};

\end{axis}
\end{tikzpicture}
\caption{$t = 50$ (IFM vs SCZM variants)}
\end{subfigure}

\vspace{0.6em}

\pgfplotslegendfromname{sharedlegend_stress}

\caption{Von Mises stress along the diagonal length (\secref{sec:sczm_isotropic}). Top row: four numerical results are compared (IFM, w/o SCZM, SCZM w/o directional correction term, and SCZM). Bottom row: a focused comparison among IFM, SCZM w/o directional correction term, and SCZM. Across both times, the SCZM solution shows excellent agreement with IFM, whereas the other approaches exhibit noticeable discrepancies.}
\label{fig:stress_over_l_theta30}
\end{figure}
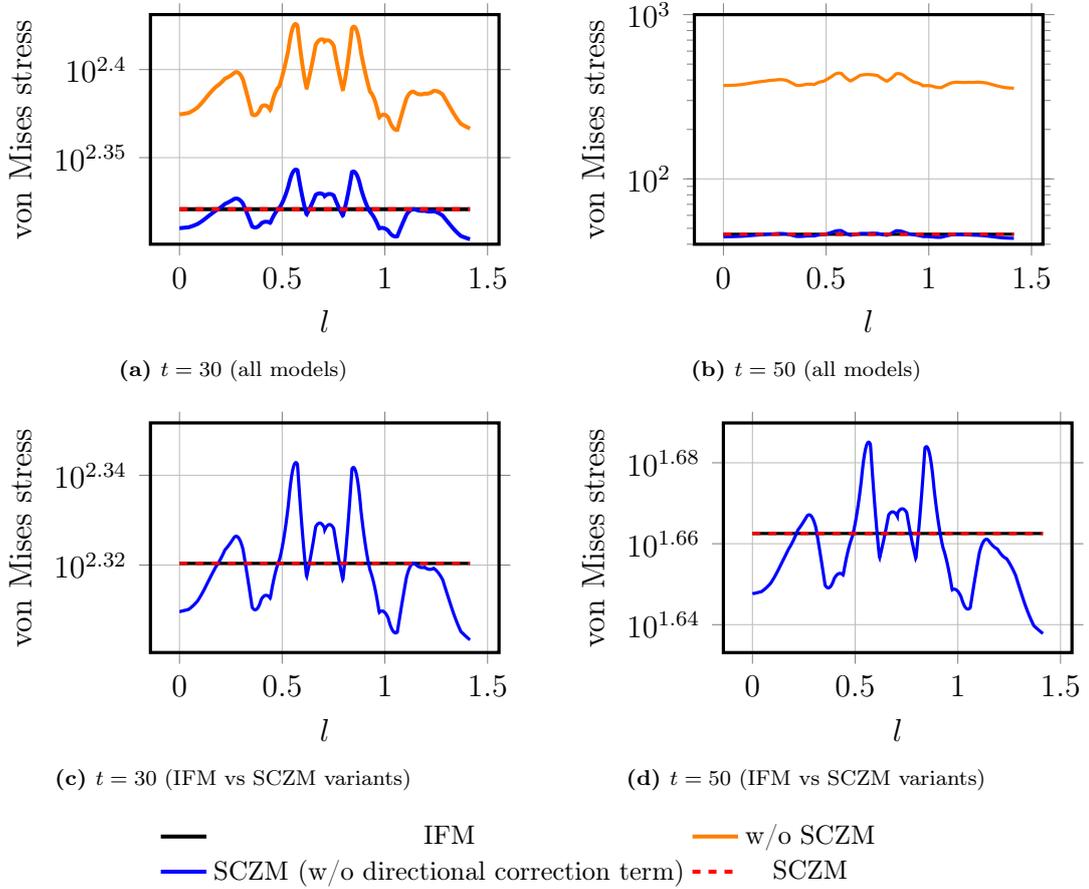

We next consider a rotated interface defined by $y = \sqrt{3}(x - 0.5) + 0.5$, discretized on a non-interface-fitted square mesh with mesh size $h=1/32$. The von Mises stress is evaluated along the diagonal direction and compared against a highly refined interface-fitted reference solution with mesh size $h_\Omega \approx 0.0028$. As shown in \figref{fig:stress_over_l_theta30}, SCZM yields accurate stress predictions along the interface, while the solution without SCZM exhibits systematic offsets.

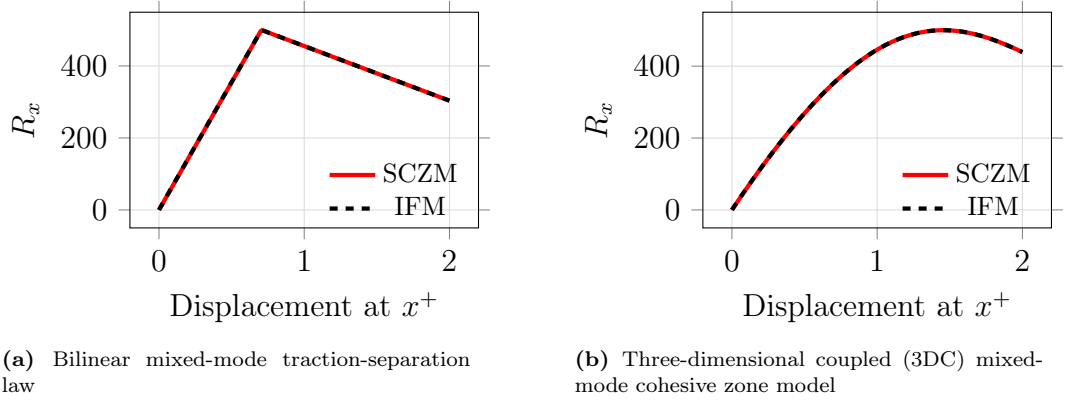
\begin{figure}[htb]
\centering
\begin{subfigure}[b]{0.45\textwidth}
\centering
\begin{tikzpicture}
\begin{axis}[
    width=\textwidth,
    height=0.72\textwidth,
    xlabel={Displacement at $x^+$},
    ylabel={$R_x$},
    grid=major,
    major grid style={line width=0.25pt, draw=gray!30},
    minor grid style={draw=none},
    tick align=outside,
    tick label style={font=\normalsize},
    label style={font=\normalsize},
    legend style={
        draw=none,
        fill=none,
        font=\footnotesize,
        at={(0.98,0.35)},
        anchor=north east
    },
    every axis plot/.append style={line width=1.6pt},
]

\addplot[
    color=red,
    solid
]
table[
    x expr=0.01*\thisrow{time},
    y=react_x,
    col sep=space,
    restrict expr to domain={\thisrow{curve_id}}{1:1}
]{react_x_vs_time_sczm_vs_bfm_data_Bilinear.txt};
\addlegendentry{SCZM}

\addplot[
    color=black,
    dashed
]
table[
    x expr=0.01*\thisrow{time},
    y=react_x,
    col sep=space,
    restrict expr to domain={\thisrow{curve_id}}{2:2}
]{react_x_vs_time_sczm_vs_bfm_data_Bilinear.txt};
\addlegendentry{IFM}

\end{axis}
\end{tikzpicture}
\caption{Bilinear mixed-mode traction-separation law}
\end{subfigure}
\hfill
\begin{subfigure}[b]{0.45\textwidth}
\centering
\begin{tikzpicture}
\begin{axis}[
    width=\textwidth,
    height=0.72\textwidth,
    xlabel={Displacement at $x^+$},
    ylabel={$R_x$},
    grid=major,
    major grid style={line width=0.25pt, draw=gray!30},
    minor grid style={draw=none},
    tick align=outside,
    tick label style={font=\normalsize},
    label style={font=\normalsize},
    legend style={
        draw=none,
        fill=none,
        font=\footnotesize,
        at={(0.98,0.35)},
        anchor=north east
    },
    every axis plot/.append style={line width=1.6pt},
]

\addplot[
    color=red,
    solid
]
table[
    x expr=0.01*\thisrow{time},
    y=react_x,
    col sep=space,
    restrict expr to domain={\thisrow{curve_id}}{1:1}
]{react_x_vs_time_sczm_vs_bfm_data_SI.txt};
\addlegendentry{SCZM}

\addplot[
    color=black,
    dashed
]
table[
    x expr=0.01*\thisrow{time},
    y=react_x,
    col sep=space,
    restrict expr to domain={\thisrow{curve_id}}{2:2}
]{react_x_vs_time_sczm_vs_bfm_data_SI.txt};
\addlegendentry{IFM}

\end{axis}
\end{tikzpicture}
\caption{Three-dimensional coupled (3DC) mixed-mode cohesive zone model}
\end{subfigure}
\caption{$x$-direction reactive forces versus the prescribed displacement at $x^+$ for two traction-separation laws (\secref{sec:sczm_isotropic}): (a) bilinear mixed-mode traction-separation law and (b) three-dimensional coupled (3DC) mixed-mode cohesive zone model. In both cases, the SCZM results are in close agreement with the corresponding IFM results.}
\label{fig:different_traction_separations}
\end{figure}

Finally, we assess the robustness of SCZM with respect to different traction-separation laws. In addition to the baseline model, we consider (i) the bilinear mixed-mode model of \citet{Camanho2002} and (ii) the coupled three-dimensional model of \citet{salehani2018coupled}. \figref{fig:different_traction_separations} shows that SCZM consistently reproduces the interface-fitted results for both models, confirming that the framework is independent of the specific cohesive law.

\subsection{Single-interface benchmark with crystal plasticity}
\label{sec:single_neml2_2D}

We next consider a history-dependent bulk constitutive model to assess the performance of SCZM in a more realistic setting. The solid response assumes a finite-strain single-crystal viscoplasticity response with isotropic elasticity and slip-based plasticity with Voce hardening. Refer to \cite{messner2025fully} for a complete description of the constitutive model. The elastic response is defined by $E = 10^3$ and $\nu = 0.3$, while plastic deformation is governed by rate-dependent crystallographic slip. All internal variables, including slip resistance and lattice orientation, are integrated implicitly using a backward Euler scheme, ensuring stable constitutive updates.

\begin{figure}[htb]
\centering
\begin{subfigure}{0.48\textwidth}
    \includegraphics[width=\linewidth]{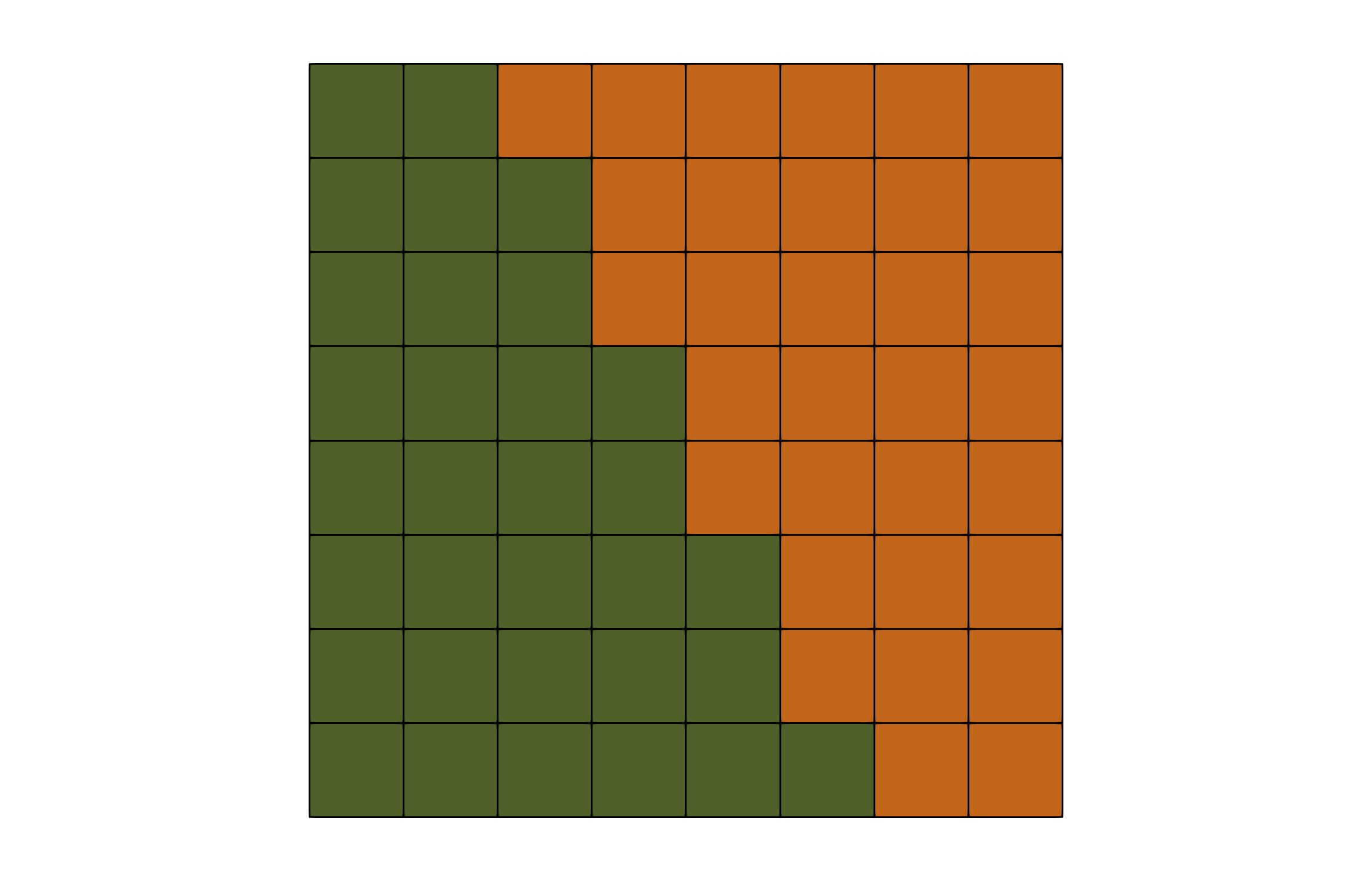}
    \caption{Non-interface-fitted mesh}
\end{subfigure}
\hfill
\begin{subfigure}{0.48\textwidth}
    \includegraphics[width=\linewidth]{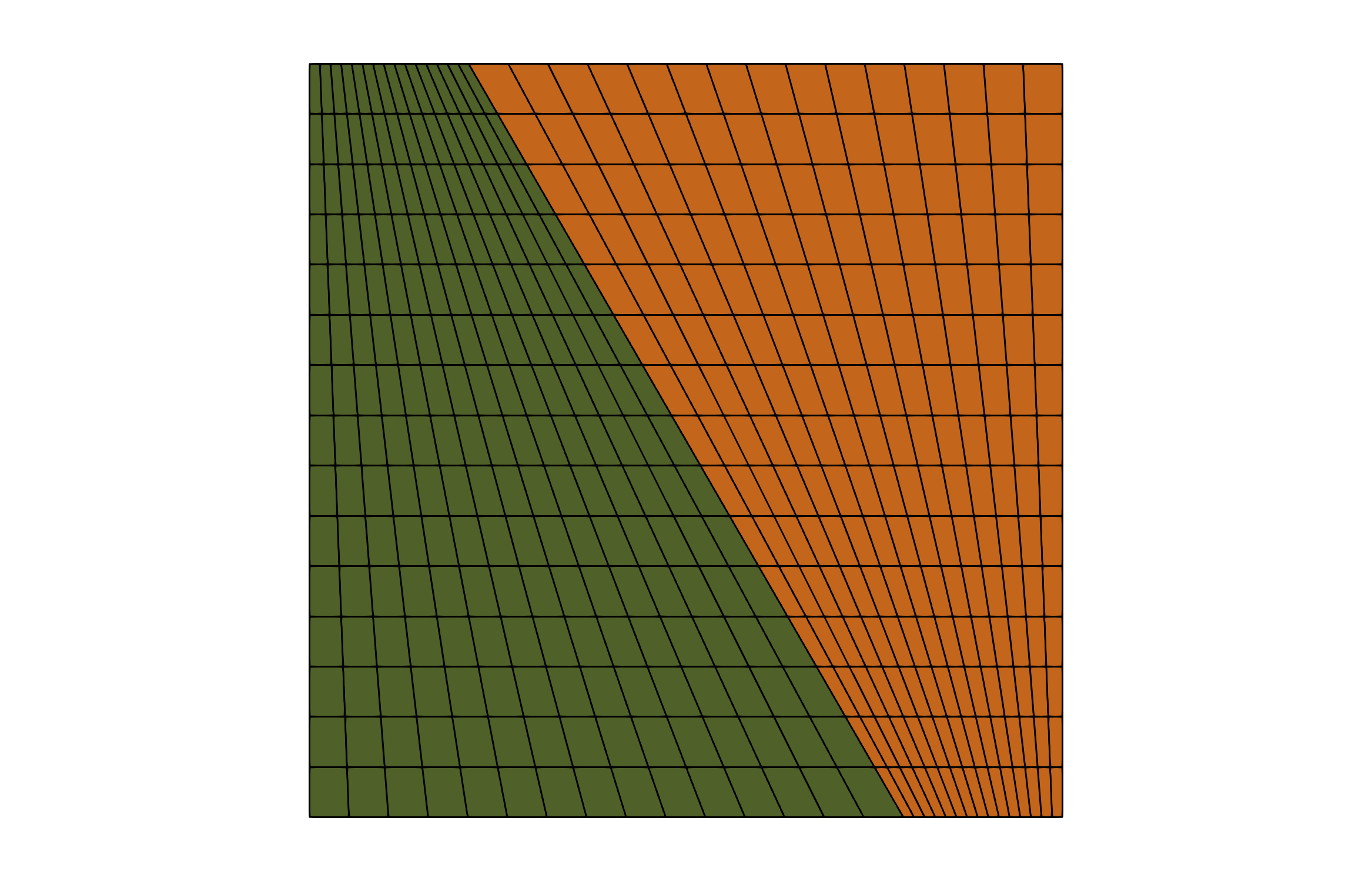}
    \caption{Interface-fitted mesh}
\end{subfigure}
\caption{Comparison of mesh configurations used in this study (\secref{sec:single_neml2_2D}).}
\label{fig:mesh_comparison}
\end{figure}

The computational domain is $[0,1]^2$, with an inclined interface defined by $y = -\sqrt{3}(x - 0.5) + 0.5$. Both interface-fitted and non-interface-fitted quadrilateral meshes are considered (see \figref{fig:mesh_comparison}). The boundary conditions follow the same setup as in \secref{sec:sczm_isotropic}

We compare four simulations:
\begin{itemize}
    \item[(a)] Non-interface-fitted mesh with na\"ive CZM enforcement,
    \item[(b)] SCZM without directional correction,
    \item[(c)] SCZM with directional correction,
    \item[(d)] Interface-fitted mesh (IFM), used as the reference solution.
\end{itemize}

The IFM simulation employs a finer mesh ($h_{\Omega} \approx 0.0471$), while the non-interface-fitted cases use $h = 0.125$. A bilinear mixed-mode traction-separation law~\cite{Camanho2002} is used for all cases with the following parameters: $K = 2\times10^{3}$, $G_{I c} = 30$, $G_{II c} = 10$, $N = 200$, $S = 100$, $\eta = 2.2$, and $\mu = 10^{-3}$. Since cohesive tractions depend on stresses from both sides of the interface, the bulk constitutive model is evaluated independently at interface quadrature points in each adjacent region, ensuring consistent evolution of history-dependent internal variables.

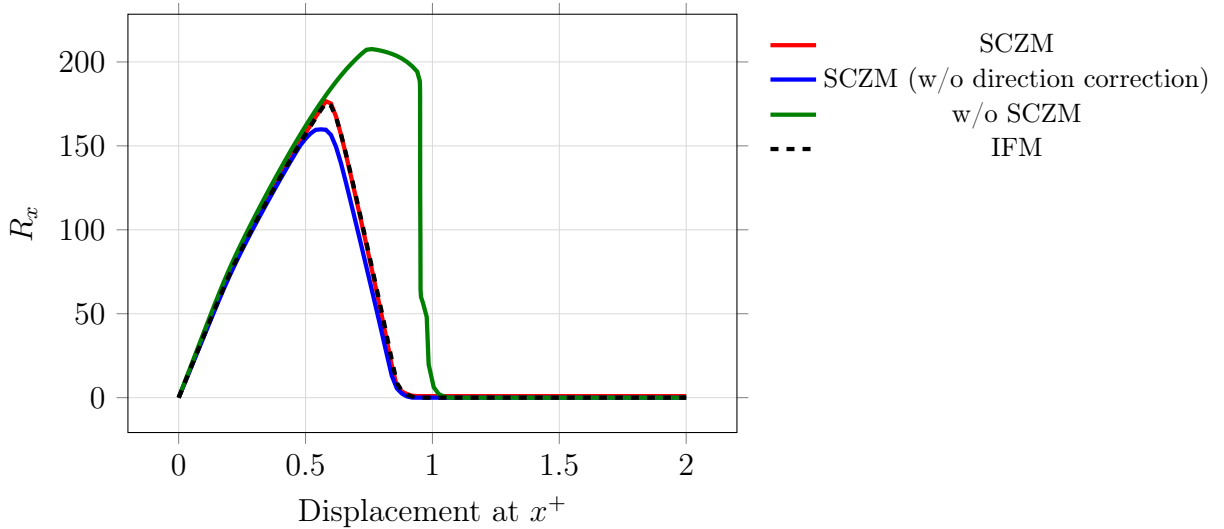
\begin{figure}[htb]
\centering
\begin{tikzpicture}
\begin{axis}[
    width=0.7\linewidth,
    height=0.52\linewidth,
    xlabel={Displacement at $x^+$},
    ylabel={$R_x$},
    grid=major,
    major grid style={line width=0.25pt, draw=gray!30},
    minor grid style={draw=none},
    tick align=outside,
    tick label style={font=\normalsize},
    label style={font=\normalsize},
    legend style={
        draw=none,
        fill=none,
        font=\footnotesize,
        at={(1.8,0.98)},
        anchor=north east
    },
    every axis plot/.append style={line width=1.6pt},
]

\addplot[
    red,
    solid
]
table[
    x=disp_xplus,
    y=react_x,
    col sep=space,
    restrict expr to domain={\thisrow{curve_id}}{1:1}
]{react_x_vs_time_mesh_study_theta30_neml2_2D_few.txt};
\addlegendentry{SCZM}

\addplot[
    blue
]
table[
    x=disp_xplus,
    y=react_x,
    col sep=space,
    restrict expr to domain={\thisrow{curve_id}}{2:2}
]{react_x_vs_time_mesh_study_theta30_neml2_2D_few.txt};
\addlegendentry{SCZM (w/o direction correction)}

\addplot[
    green!50!black
]
table[
    x=disp_xplus,
    y=react_x,
    col sep=space,
    restrict expr to domain={\thisrow{curve_id}}{3:3}
]{react_x_vs_time_mesh_study_theta30_neml2_2D_few.txt};
\addlegendentry{w/o SCZM}

\addplot[
    black,
    dashed
]
table[
    x=disp_xplus,
    y=react_x,
    col sep=space,
    restrict expr to domain={\thisrow{curve_id}}{4:4}
]{react_x_vs_time_mesh_study_theta30_neml2_2D_few.txt};
\addlegendentry{IFM}

\end{axis}
\end{tikzpicture}
\caption{Reaction force in the $x$-direction ($R_x$) as a function of the prescribed displacement at $x^+$ for 2D plasticity simulations (\secref{sec:single_neml2_2D}).}
\label{fig:2D_grain_t_Rx_2D_plasticity}
\end{figure}

\figref{fig:2D_grain_t_Rx_2D_plasticity} shows the reaction force in the $x$-direction as a function of time. The SCZM formulation with directional correction (case (c)) closely matches the IFM reference, while the formulation without directional correction (case (b)) exhibits noticeable deviation. The naive CZM treatment on non-interface-fitted meshes (case (a)) leads to significant errors.

\begin{figure}[htb]
\centering

\begin{subfigure}{0.49\textwidth}
    \includegraphics[width=\linewidth,trim = {500 100 350 100}, clip]{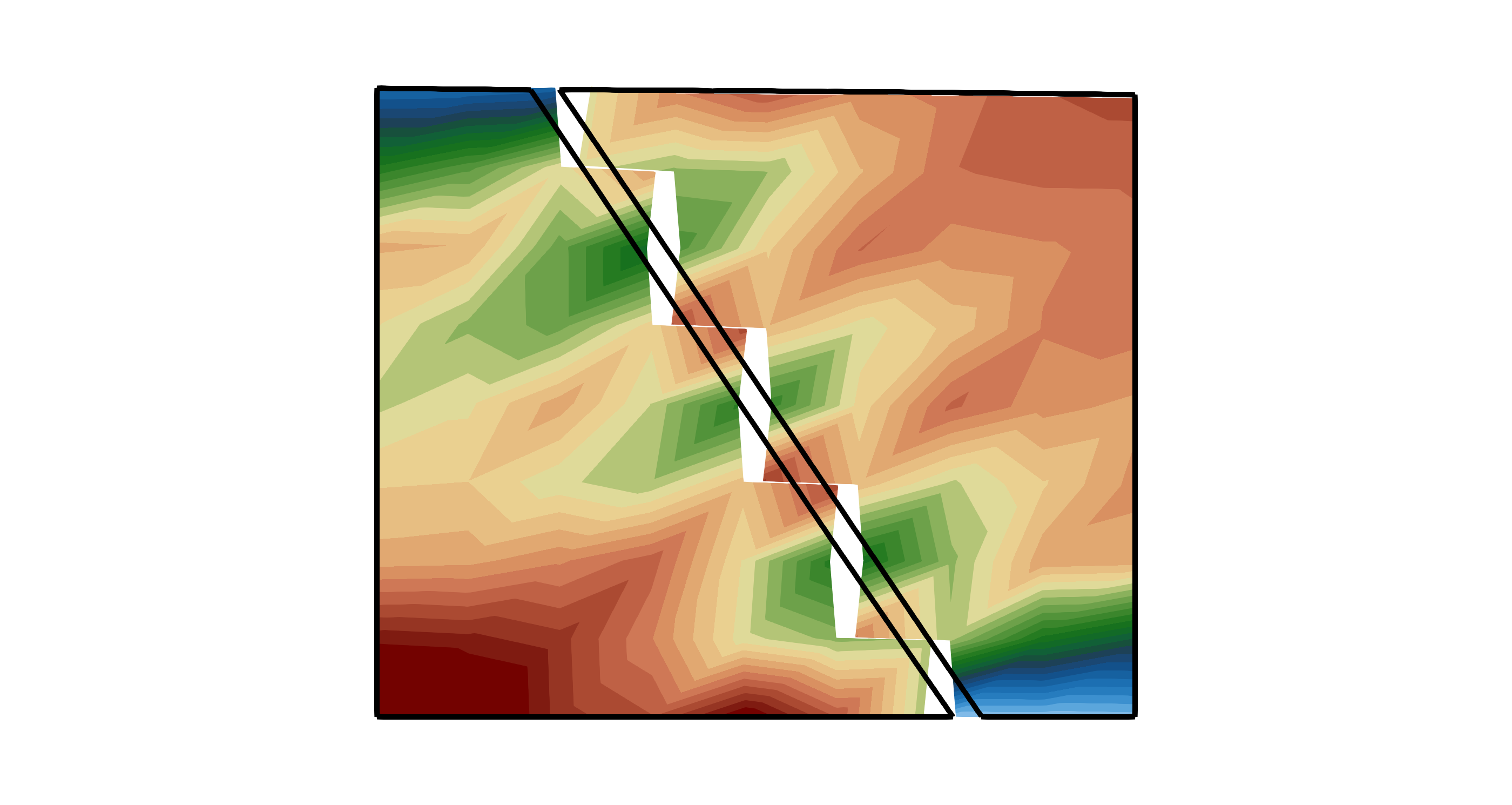}
    \caption{w/o SCZM}
\end{subfigure}
\hfill
\begin{subfigure}{0.49\textwidth}
    \includegraphics[width=\linewidth,trim = {500 100 350 100}, clip]{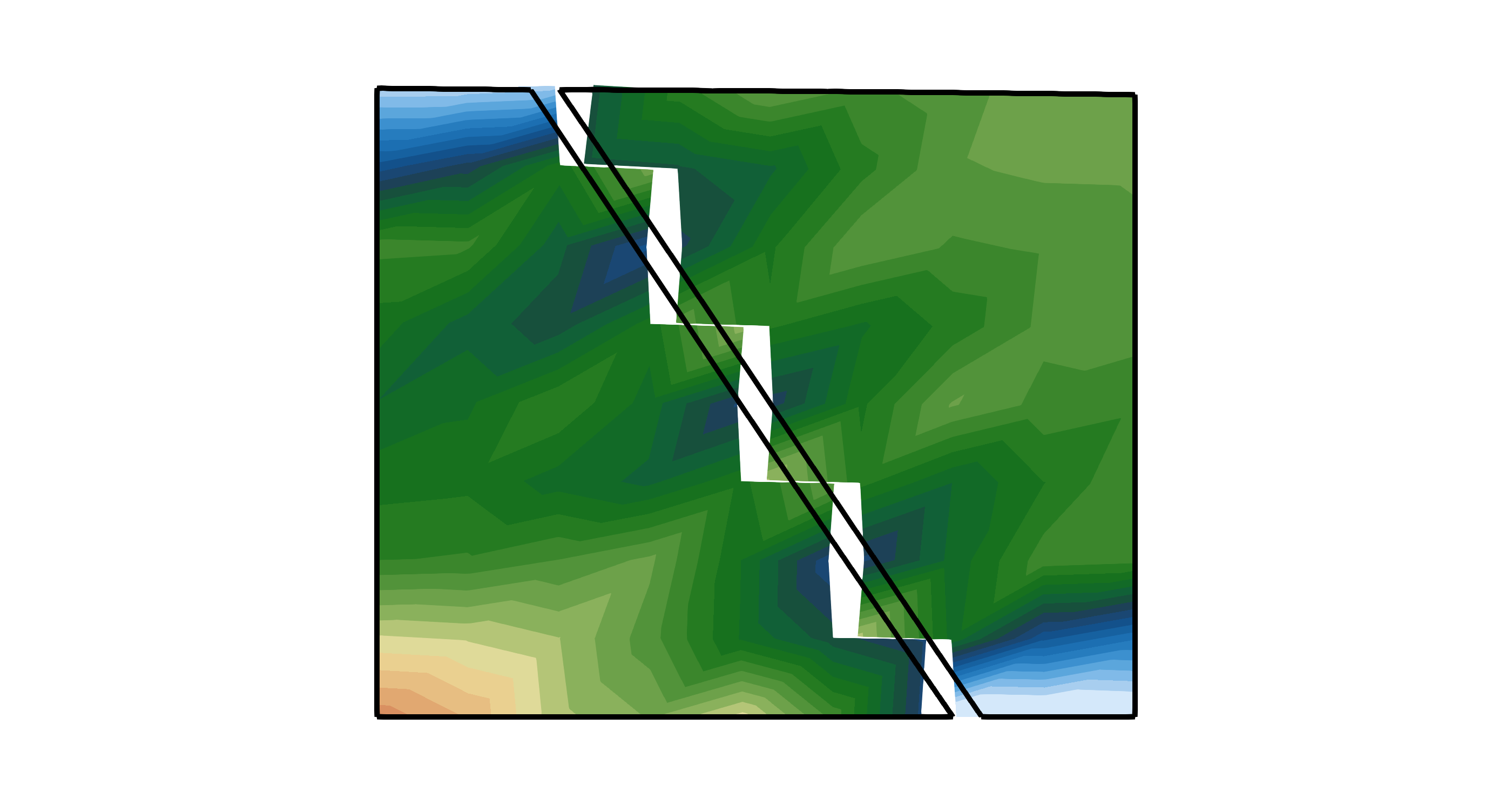}
    \caption{SCZM w/o directional correction}
\end{subfigure}
\\
\begin{subfigure}{0.49\textwidth}
    \includegraphics[width=\linewidth,trim = {500 100 350 100}, clip]{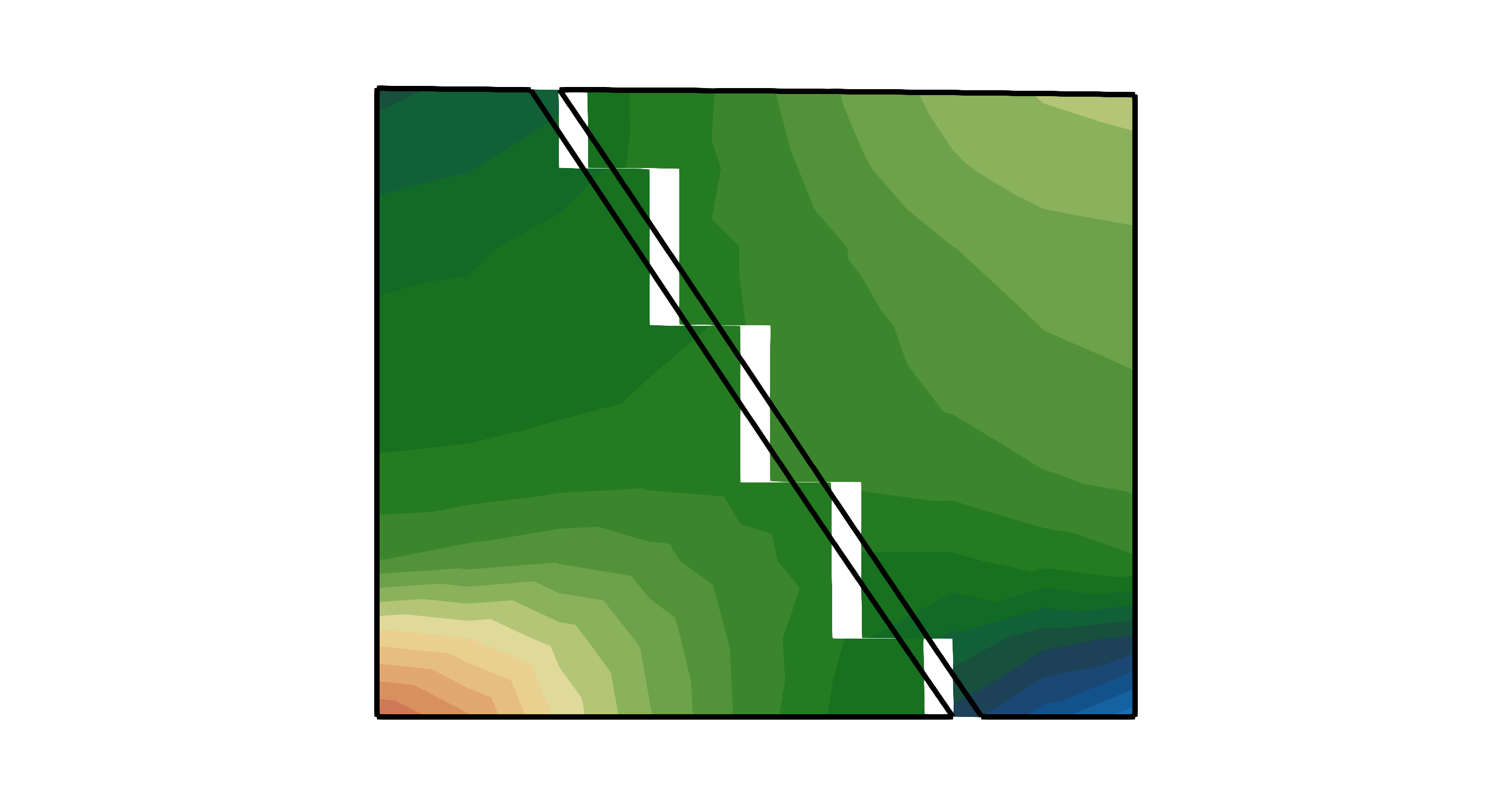}
    \caption{SCZM}
\end{subfigure}
\hfill
\begin{subfigure}{0.49\textwidth}
    \includegraphics[width=\linewidth,trim = {500 100 350 100}, clip]{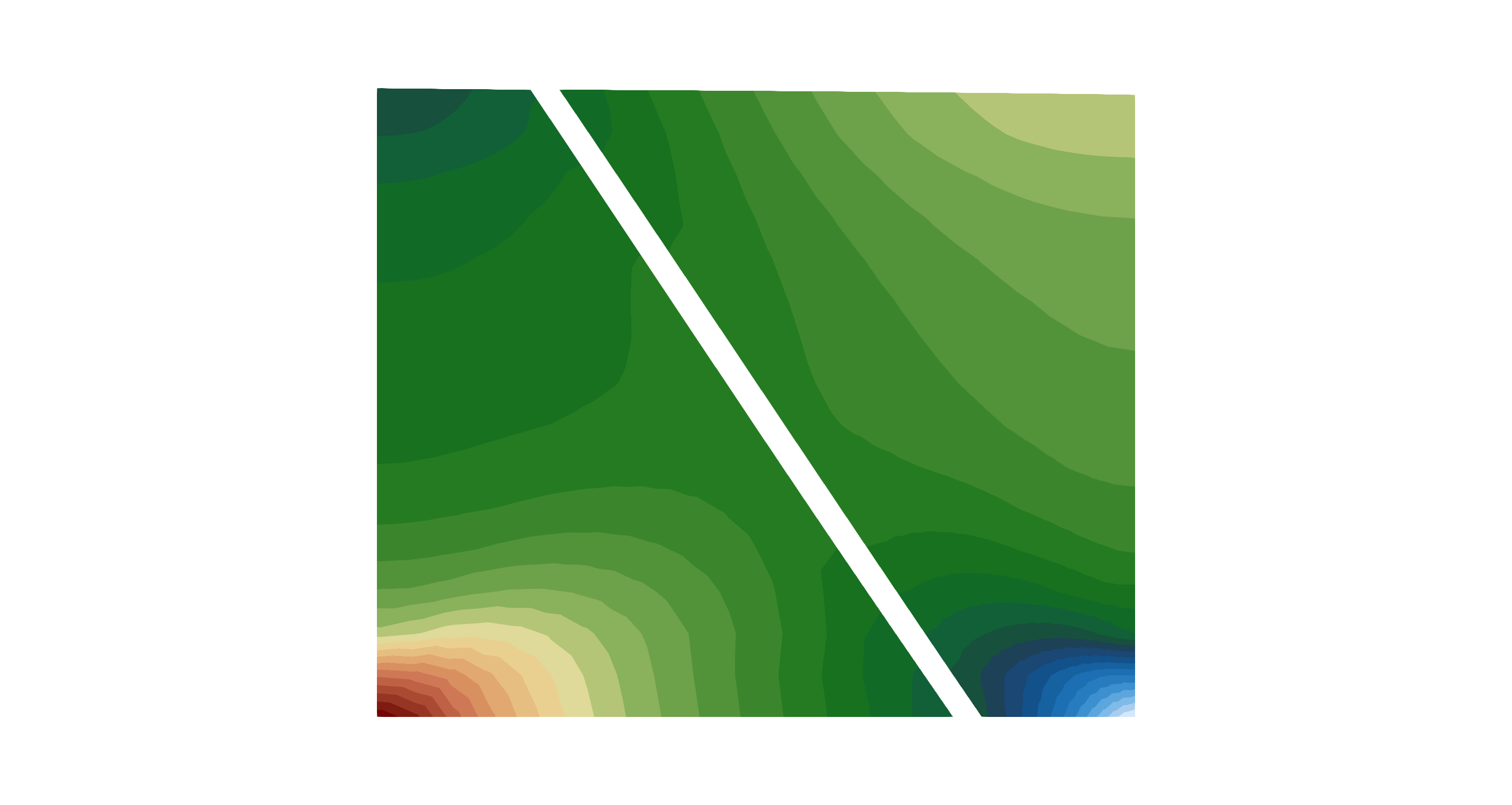}
    \caption{IFM benchmark}
\end{subfigure}

\includegraphics[width=0.15\linewidth]{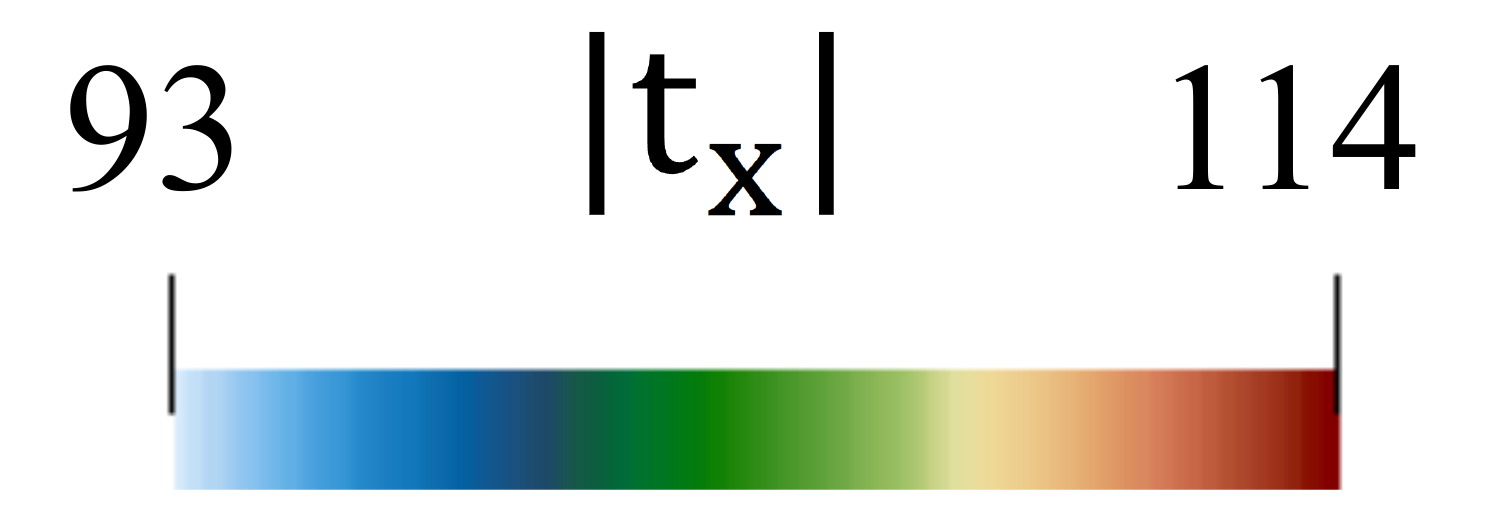}

\caption{
Comparison of $\abs{t_x}$ distributions at $t = 30$ (\secref{sec:single_neml2_2D}): naive CZM without SCZM, SCZM without directional correction, SCZM, and the IFM benchmark. The absence of directional correction leads to visible stress distortion along the surrogate interface. The SCZM formulation restores geometric consistency and closely matches the IFM reference solution.
The black lines indicate the reference interface-fitted configuration under the applied displacement.
}
\label{fig:2D_plasticity_t30}
\end{figure}

\begin{figure}[htb]
\centering

\begin{subfigure}{0.49\textwidth}
    \includegraphics[width=\linewidth,trim = {500 100 350 100}, clip]{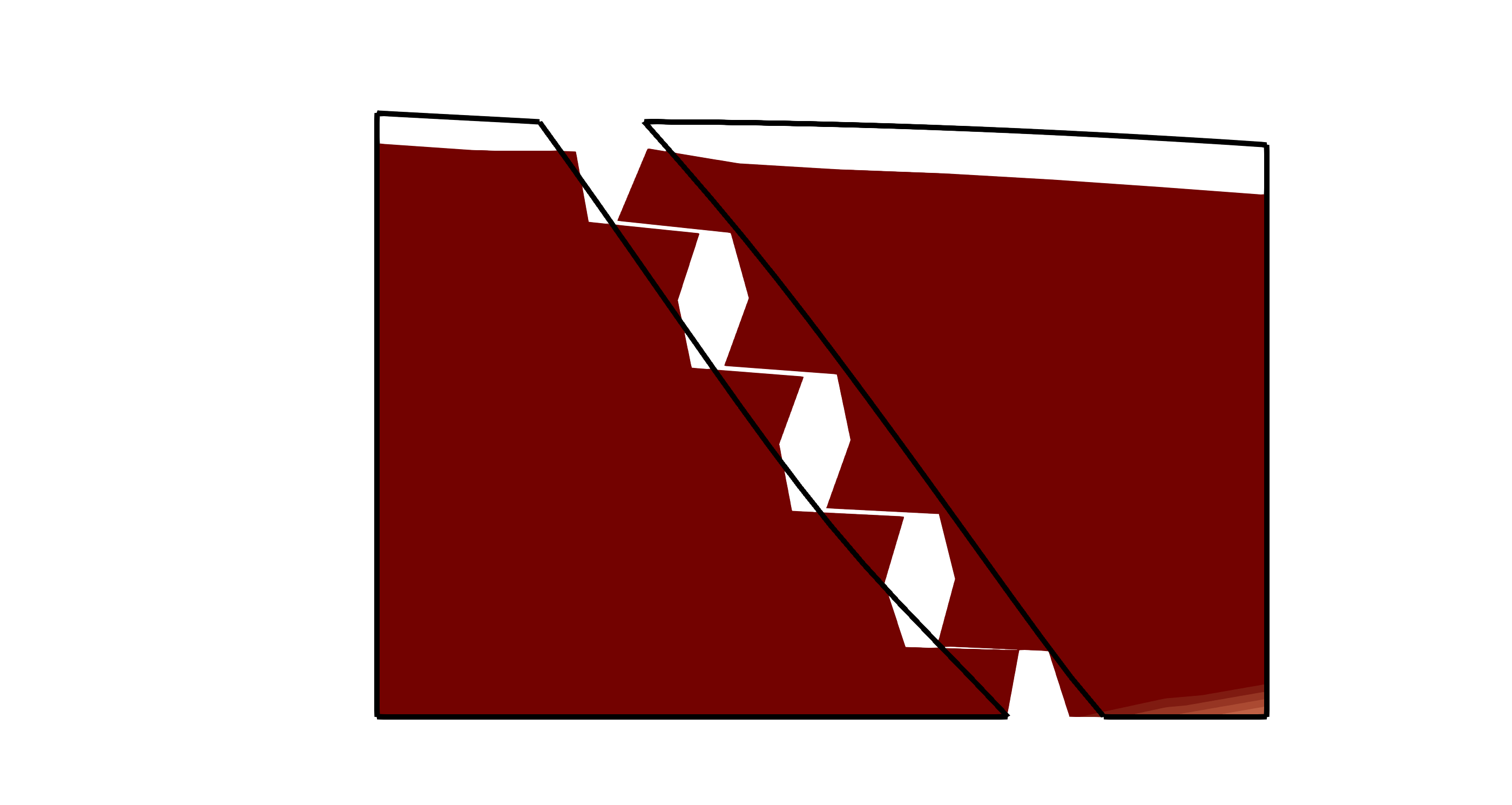}
    \caption{w/o SCZM}
\end{subfigure}
\hfill
\begin{subfigure}{0.49\textwidth}
    \includegraphics[width=\linewidth,trim = {500 100 350 100}, clip]{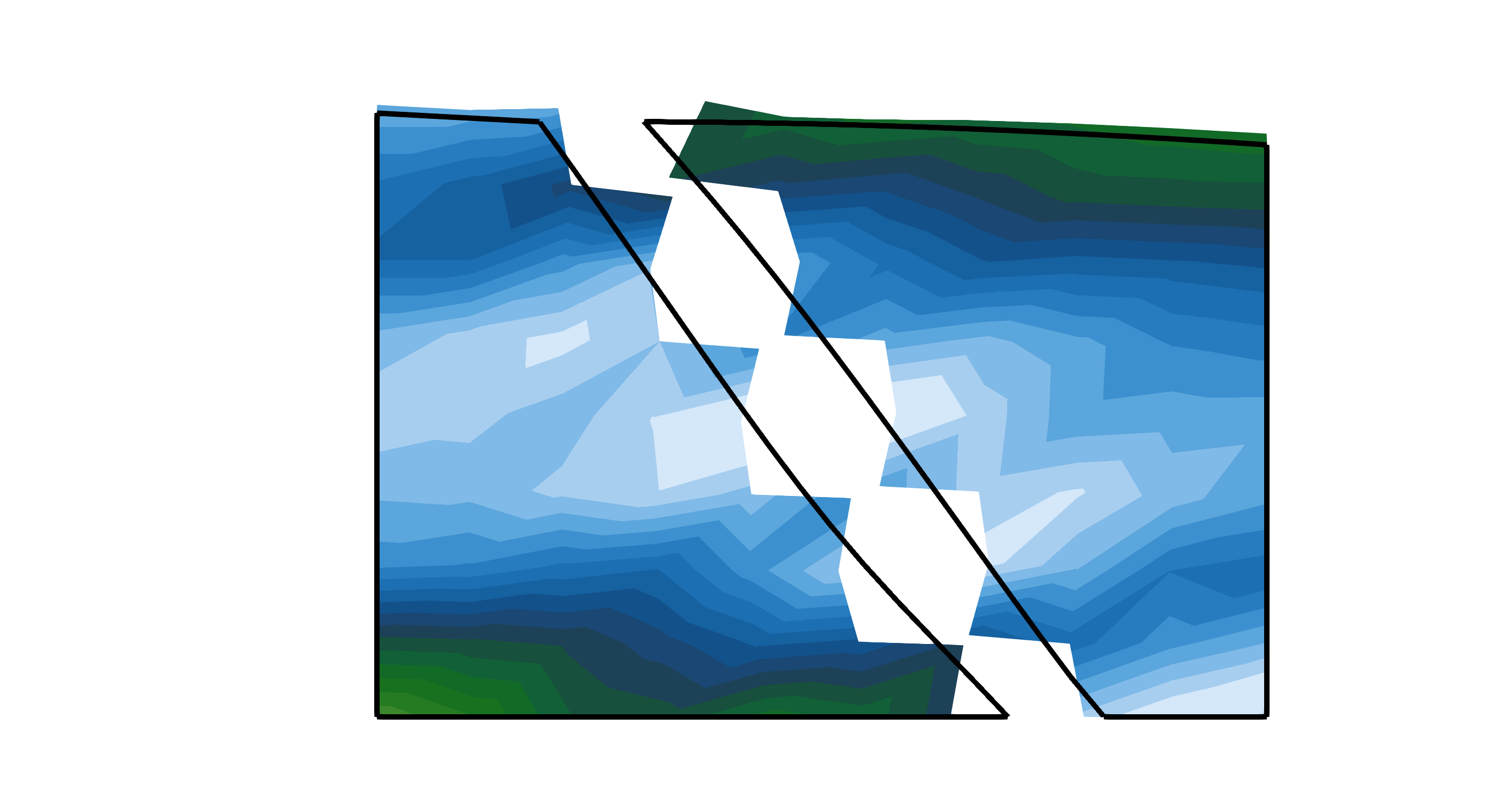}
    \caption{SCZM w/o directional correction}
\end{subfigure}
\\
\begin{subfigure}{0.49\textwidth}
    \includegraphics[width=\linewidth,trim = {500 100 350 100}, clip]{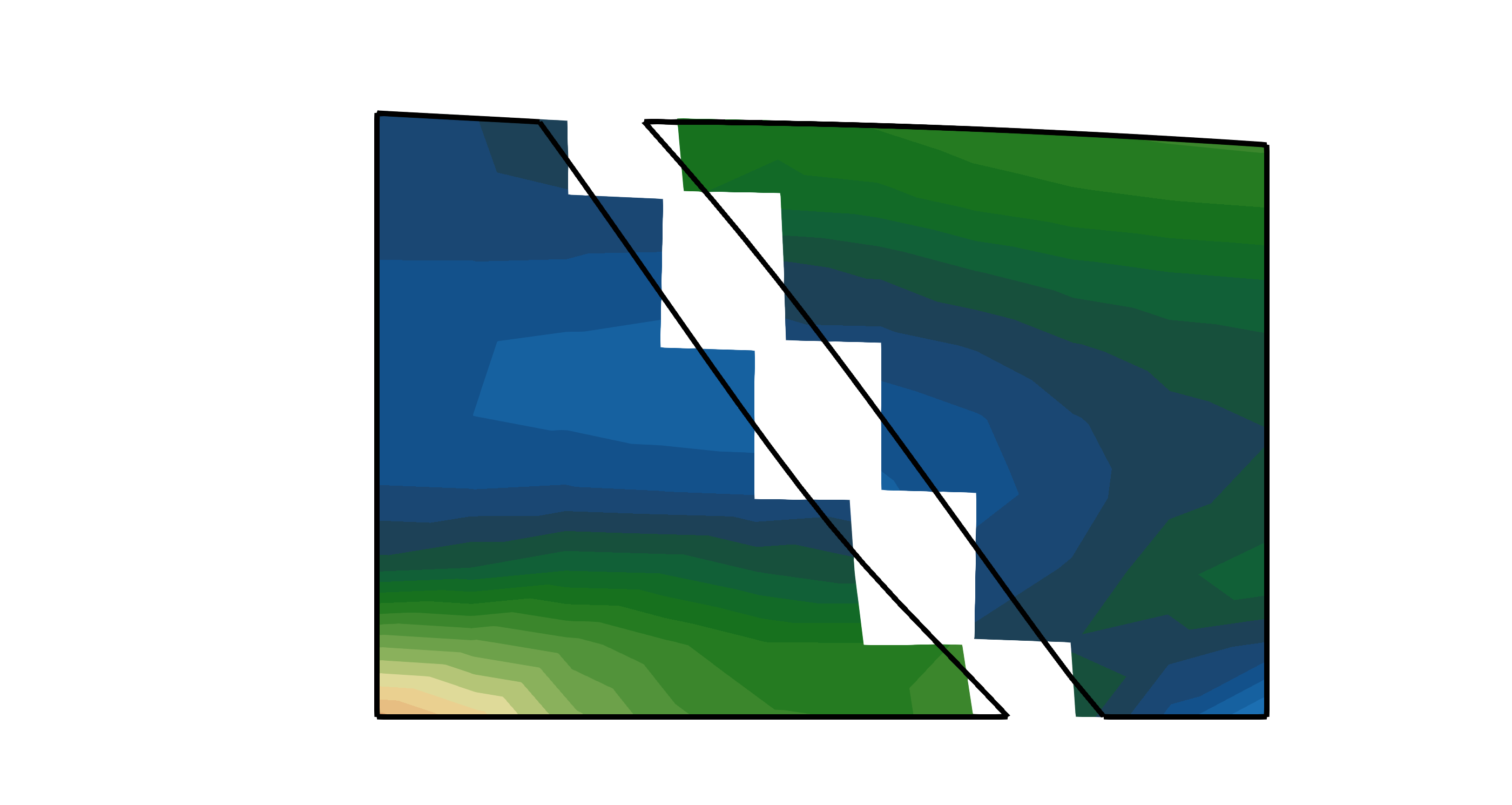}
    \caption{SCZM}
\end{subfigure}
\hfill
\begin{subfigure}{0.49\textwidth}
    \includegraphics[width=\linewidth,trim = {500 100 350 100}, clip]{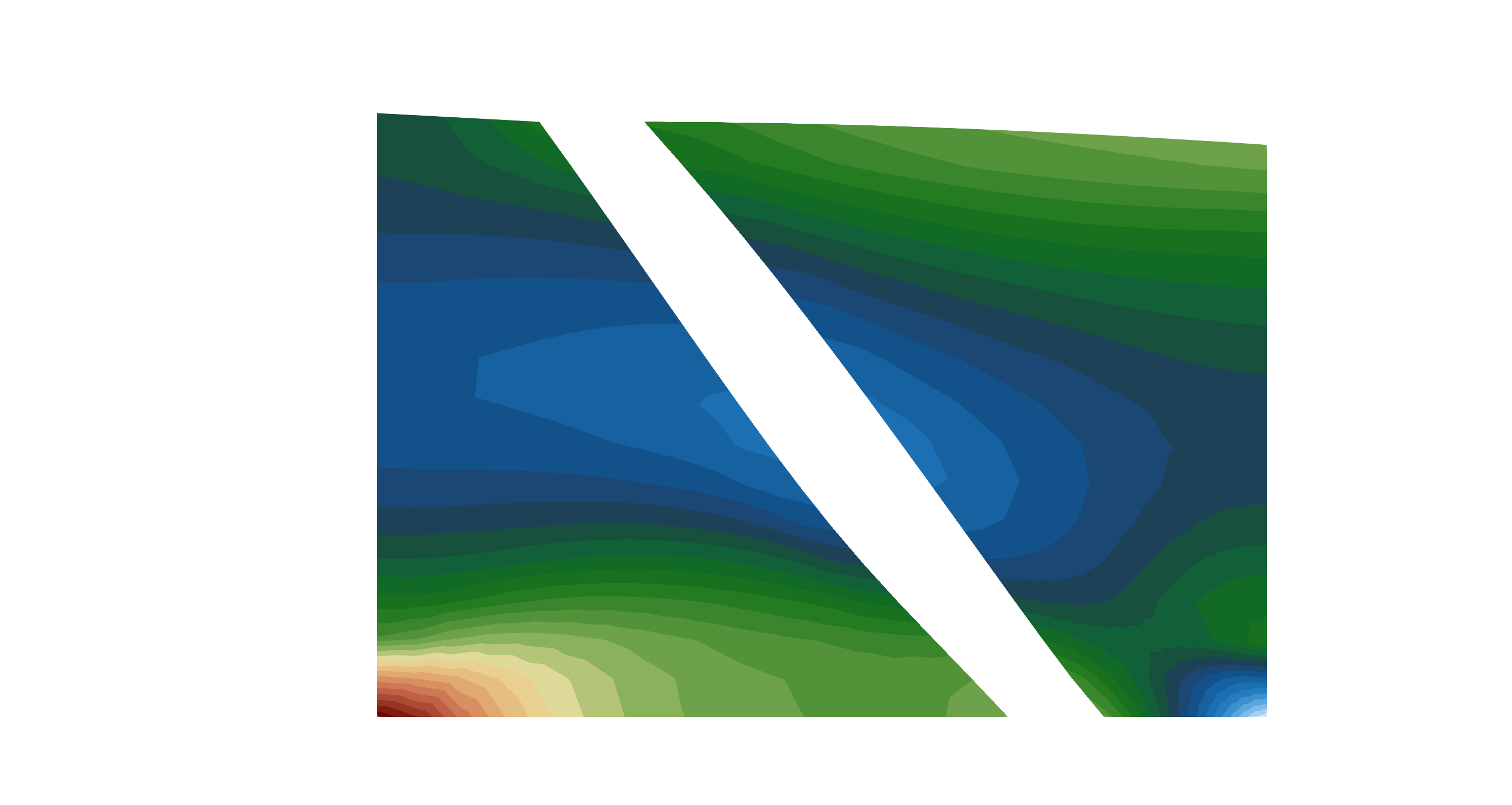}
    \caption{IFM benchmark}
\end{subfigure}

\includegraphics[width=0.15\linewidth]{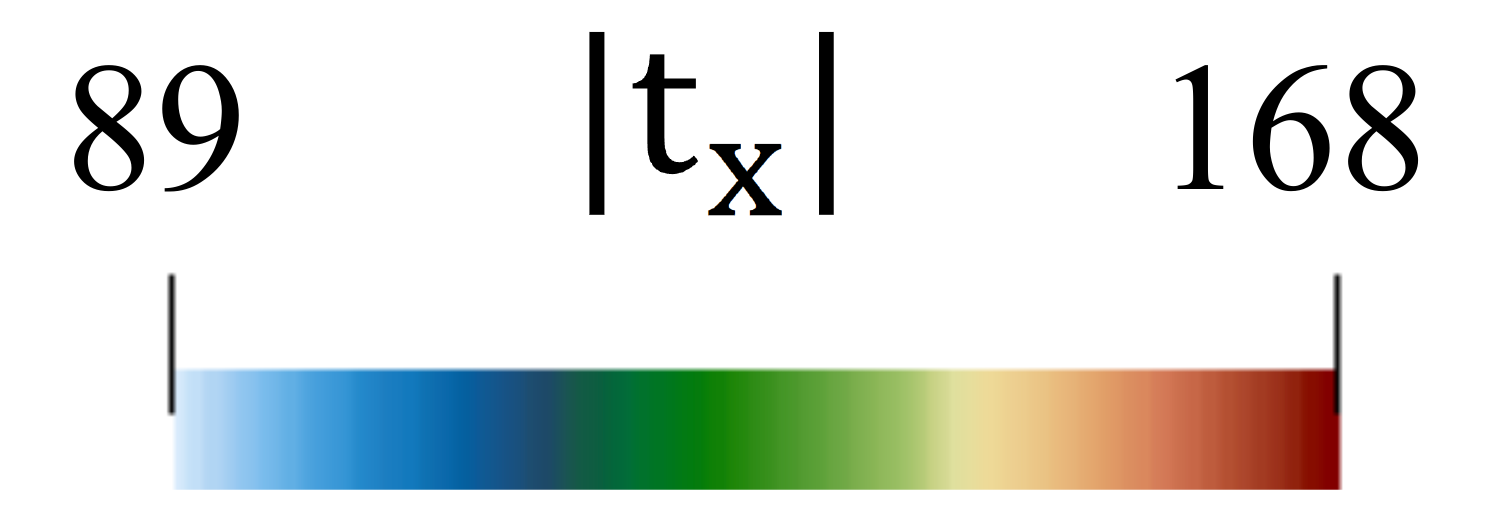}

\caption{
Comparison of $\abs{t_x}$ distributions at $t = 70$ (\secref{sec:single_neml2_2D}): naive CZM without SCZM, SCZM without directional correction, SCZM, and the IFM benchmark. At a later time, the stress distortion caused by the absence of directional correction becomes more pronounced along the surrogate interface. The consistent SCZM formulation maintains geometric consistency and remains in close agreement with the IFM reference solution.
The black lines indicate the reference interface-fitted configuration under the applied displacement.
}
\label{fig:2D_plasticity_t70}
\end{figure}

To further assess the role of directional correction, we examine the traction magnitude $|\bs{t}_x| = |\bs{\sigma} \cdot \bs{e}_x|$ along the interface at selected times (see \figref{fig:2D_plasticity_t30} and \figref{fig:2D_plasticity_t70}). Without directional correction, the solution exhibits pronounced stress distortion along the surrogate interface. In contrast, the full SCZM formulation yields stress fields that are in close agreement with the IFM solution.

These results demonstrate that, in the presence of nonlinear and history-dependent constitutive behavior, accurate enforcement of interface conditions requires both geometric scaling and directional correction. The SCZM framework provides a consistent means to achieve this on non-interface-fitted meshes.

\subsection{Crystal plasticity on 2D polycrystalline RVE}
\label{sec:multi_neml2_2D}

We consider a two-dimensional polycrystalline RVE consisting of five grains, with an interface-fitted discretization used as the reference solution. The bulk constitutive response is modeled using the same model as in \secref{sec:single_neml2_2D}. The computational domain is $[-0.6,0.6]\times[-0.4,0.4]$, with boundary conditions $u_x=0$ on the left, $u_y=0$ on the bottom, and $u_x = 10^{-2}t$ on the right. Simulations are performed for $t \in [0,100]$ with $\Delta t = 1$. A bilinear mixed-mode traction-separation law is used at grain boundaries with the following parameters: $K = 50$, $G_{I c} = 1$, $G_{II c} = 3$, $N = 3$, $S = 5$, $\eta = 2.2$, and $\mu = 0$.

\begin{figure}[htb]
\centering
\includegraphics[width=0.75\linewidth]{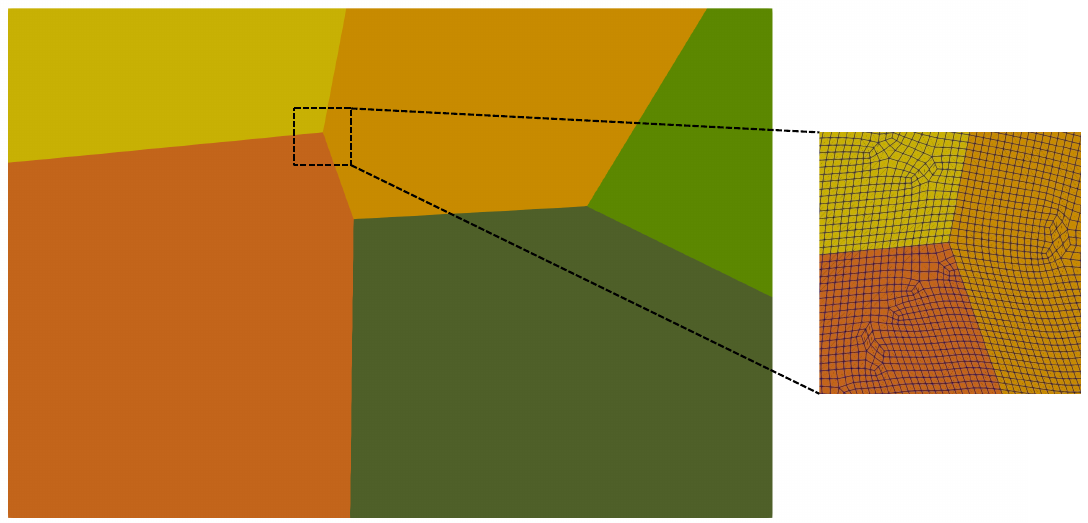}
\caption{2D polycrystalline interface-fitted meshes with a characteristic mesh size of $h_{\Omega} \approx 2.92\times10^{-3}$ (\secref{sec:multi_neml2_2D}).}
\label{fig:2D_grains_bfm}
\end{figure}

The reference solution is obtained using an interface-fitted quadrilateral mesh with characteristic size $h_\Omega \approx 2.92\times10^{-3}$ (see \figref{fig:2D_grains_bfm}). SCZM is applied on non-interface-fitted quadrilateral meshes with resolutions $h = 0.2/2^{\text{lvl}}$, $\text{lvl}=0,\dots,6$, to assess convergence toward the reference solution.

\begin{figure}[htb]
\centering

\begin{subfigure}{0.45\linewidth}
\centering
\begin{tikzpicture}
\begin{axis}[
    width=\linewidth,
    height=0.72\linewidth,
    xlabel={Displacement at $x^+$},
    ylabel={$R_x$},
    xmin=0, xmax=1,
    ymin=-0.5, ymax=3.5,
    grid=major,
    major grid style={line width=0.25pt, draw=gray!30},
    minor grid style={draw=none},
    tick align=outside,
    tick label style={font=\normalsize},
    label style={font=\normalsize},
    every axis plot/.append style={line width=1.4pt}
]

\addplot[color=c1]
table[
    x expr=0.01*\thisrow{time}, y=react_x,
    col sep=space,
    restrict expr to domain={\thisrow{mesh_ny}}{4:4}
]{react_x_vs_time_mesh_study.txt};

\addplot[color=c2]
table[
    x expr=0.01*\thisrow{time}, y=react_x,
    col sep=space,
    restrict expr to domain={\thisrow{mesh_ny}}{8:8}
]{react_x_vs_time_mesh_study.txt};

\addplot[color=c3]
table[
    x expr=0.01*\thisrow{time}, y=react_x,
    col sep=space,
    restrict expr to domain={\thisrow{mesh_ny}}{16:16}
]{react_x_vs_time_mesh_study.txt};

\addplot[color=c4]
table[
    x expr=0.01*\thisrow{time}, y=react_x,
    col sep=space,
    restrict expr to domain={\thisrow{mesh_ny}}{32:32}
]{react_x_vs_time_mesh_study.txt};

\addplot[color=c5]
table[
    x expr=0.01*\thisrow{time}, y=react_x,
    col sep=space,
    restrict expr to domain={\thisrow{mesh_ny}}{64:64}
]{react_x_vs_time_mesh_study.txt};

\addplot[color=c6]
table[
    x expr=0.01*\thisrow{time}, y=react_x,
    col sep=space,
    restrict expr to domain={\thisrow{mesh_ny}}{128:128}
]{react_x_vs_time_mesh_study.txt};

\addplot[color=c7]
table[
    x expr=0.01*\thisrow{time}, y=react_x,
    col sep=space,
    restrict expr to domain={\thisrow{mesh_ny}}{256:256}
]{react_x_vs_time_mesh_study.txt};

\addplot[color=black,dashed]
table[
    x expr=0.01*\thisrow{time}, y=react_x,
    col sep=space,
    restrict expr to domain={\thisrow{mesh_ny}}{-1:-1}
]{react_x_vs_time_mesh_study.txt};

\end{axis}
\end{tikzpicture}
\caption{Evolution of the reaction force in the $x$-direction ($R_x$).}
\label{fig:2D_grain_t_Rx}
\end{subfigure}
\hfill
\begin{subfigure}{0.45\linewidth}
\centering
\begin{tikzpicture}
\begin{axis}[
    width=\linewidth,
    height=0.72\linewidth,
    xlabel={Displacement at $x^+$},
    ylabel={$\max(d)$},
    xmin=0, xmax=1,
    grid=major,
    major grid style={line width=0.25pt, draw=gray!30},
    minor grid style={draw=none},
    tick align=outside,
    tick label style={font=\normalsize},
    label style={font=\normalsize},
    every axis plot/.append style={line width=1.4pt}
]

\addplot[color=c1]
table[
    x expr=0.01*\thisrow{time}, y=max_damage,
    col sep=space,
    restrict expr to domain={\thisrow{mesh_ny}}{4:4}
]{max_damage_vs_time_mesh_study.txt};

\addplot[color=c2]
table[
    x expr=0.01*\thisrow{time}, y=max_damage,
    col sep=space,
    restrict expr to domain={\thisrow{mesh_ny}}{8:8}
]{max_damage_vs_time_mesh_study.txt};

\addplot[color=c3]
table[
    x expr=0.01*\thisrow{time}, y=max_damage,
    col sep=space,
    restrict expr to domain={\thisrow{mesh_ny}}{16:16}
]{max_damage_vs_time_mesh_study.txt};

\addplot[color=c4]
table[
    x expr=0.01*\thisrow{time}, y=max_damage,
    col sep=space,
    restrict expr to domain={\thisrow{mesh_ny}}{32:32}
]{max_damage_vs_time_mesh_study.txt};

\addplot[color=c5]
table[
    x expr=0.01*\thisrow{time}, y=max_damage,
    col sep=space,
    restrict expr to domain={\thisrow{mesh_ny}}{64:64}
]{max_damage_vs_time_mesh_study.txt};

\addplot[color=c6]
table[
    x expr=0.01*\thisrow{time}, y=max_damage,
    col sep=space,
    restrict expr to domain={\thisrow{mesh_ny}}{128:128}
]{max_damage_vs_time_mesh_study.txt};

\addplot[color=c7]
table[
    x expr=0.01*\thisrow{time}, y=max_damage,
    col sep=space,
    restrict expr to domain={\thisrow{mesh_ny}}{256:256}
]{max_damage_vs_time_mesh_study.txt};

\addplot[color=black,dashed]
table[
    x expr=0.01*\thisrow{time}, y=max_damage,
    col sep=space,
    restrict expr to domain={\thisrow{mesh_ny}}{-1:-1}
]{max_damage_vs_time_mesh_study.txt};

\end{axis}
\end{tikzpicture}
\caption{Evolution of the maximum damage.}
\label{fig:2D_grain_t_damage_max}
\end{subfigure}

\vspace{0.6em}

\begin{tikzpicture}
\begin{axis}[
    hide axis,
    xmin=0, xmax=1,
    ymin=0, ymax=1,
    legend columns=3,
    legend style={
        draw=none,
        fill=none,
        font=\footnotesize,
        /tikz/every even column/.append style={column sep=0.35cm}
    }
]
\addlegendimage{color=c1, line width=1.4pt}
\addlegendentry{SCZM ($h=0.8/4$)}

\addlegendimage{color=c2, line width=1.4pt}
\addlegendentry{SCZM ($h=0.8/8$)}

\addlegendimage{color=c3, line width=1.4pt}
\addlegendentry{SCZM ($h=0.8/16$)}

\addlegendimage{color=c4, line width=1.4pt}
\addlegendentry{SCZM ($h=0.8/32$)}

\addlegendimage{color=c5, line width=1.4pt}
\addlegendentry{SCZM ($h=0.8/64$)}

\addlegendimage{color=c6, line width=1.4pt}
\addlegendentry{SCZM ($h=0.8/128$)}

\addlegendimage{color=c7, line width=1.4pt}
\addlegendentry{SCZM ($h=0.8/256$)}

\addlegendimage{color=black,dashed, line width=1.4pt}
\addlegendentry{IFM ($h_{\Omega}\approx 0.00255$)}
\end{axis}
\end{tikzpicture}

\caption{Reaction force in the $x$-direction ($R_x$) and maximum damage as functions of the prescribed displacement at $x^+$ for 2D polycrystalline RVEs (\secref{sec:multi_neml2_2D}) across different mesh resolutions.}
\label{fig:2D_grain_t_Rx_damage}
\end{figure}
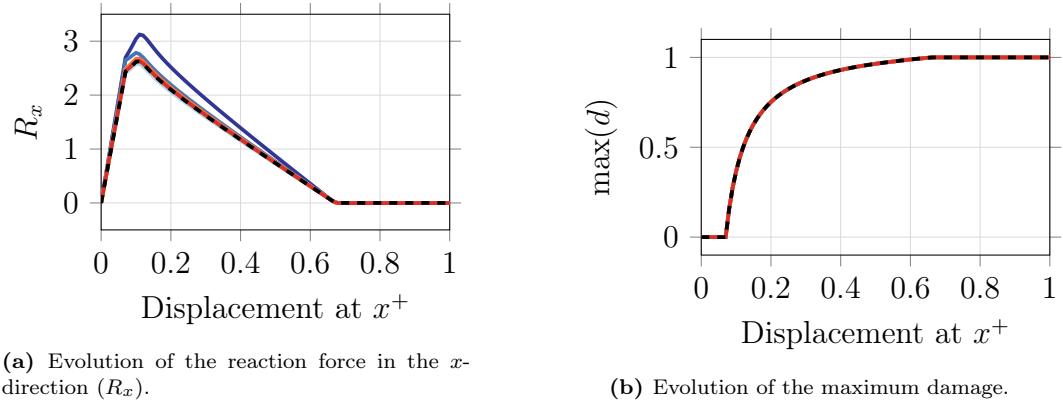
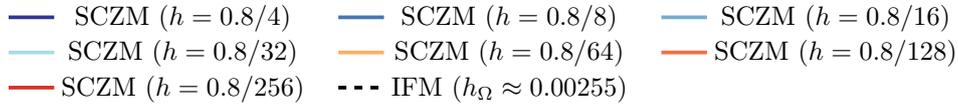

\figref{fig:2D_grain_t_Rx} shows the force-displacement curve in the loading direction. While coarse meshes exhibit noticeable discrepancies, the SCZM solution converges systematically to the interface-fitted result with mesh refinement. A similar level of agreement is observed in the evolution of the maximum damage (\figref{fig:2D_grain_t_damage_max}), where all non-interface-fitted simulations closely match the reference solution.

\begin{figure}[htb]
    \centering
    \begin{subfigure}{0.32\textwidth}
        \includegraphics[width=\linewidth,trim = {50 200 100 200}, clip]{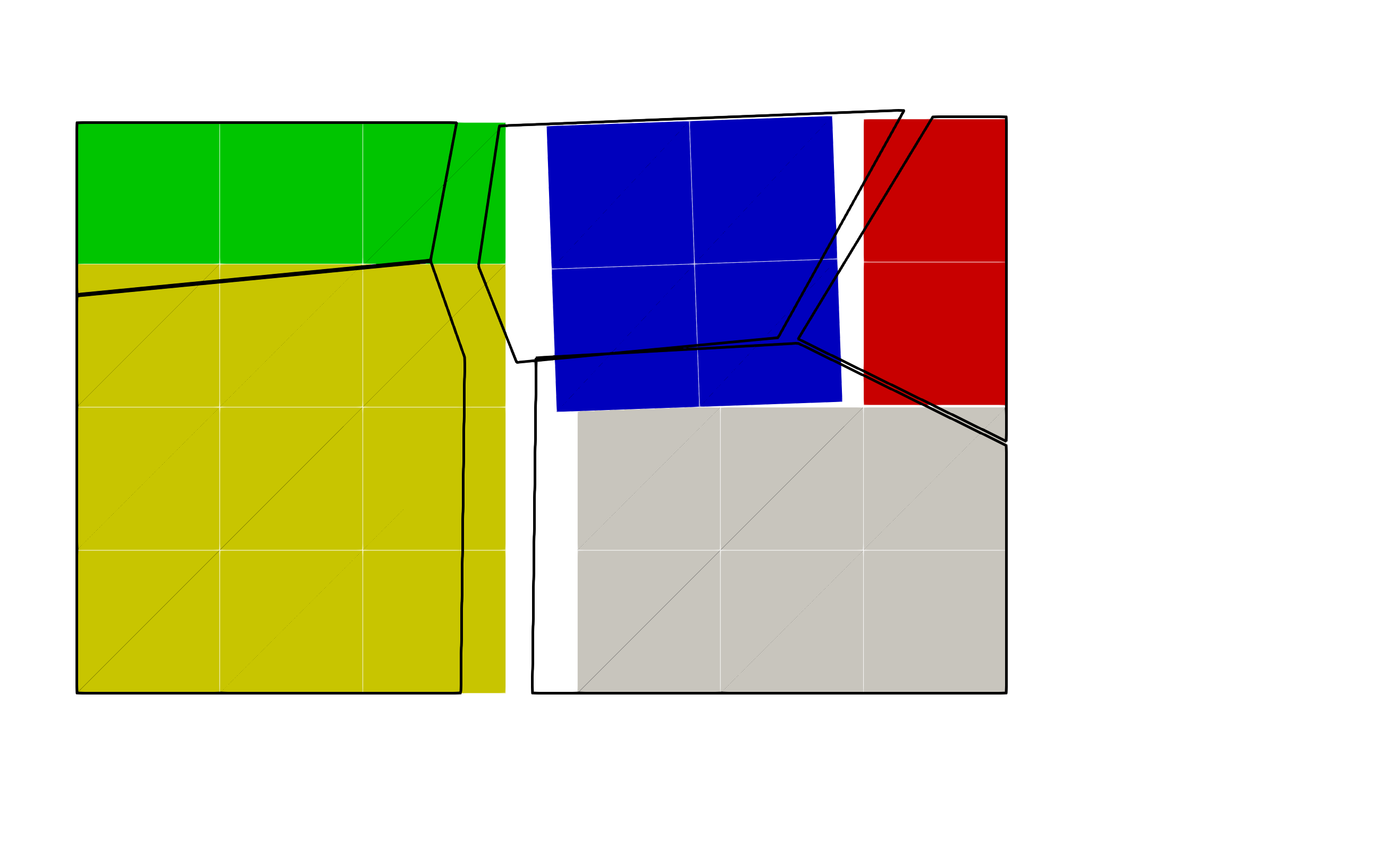}
        \caption{$h=0.2$}
        \label{fig:ny4_t10}
    \end{subfigure}
    \hfill
    \begin{subfigure}{0.32\textwidth}
        \includegraphics[width=\linewidth,trim = {50 200 100 200}, clip]{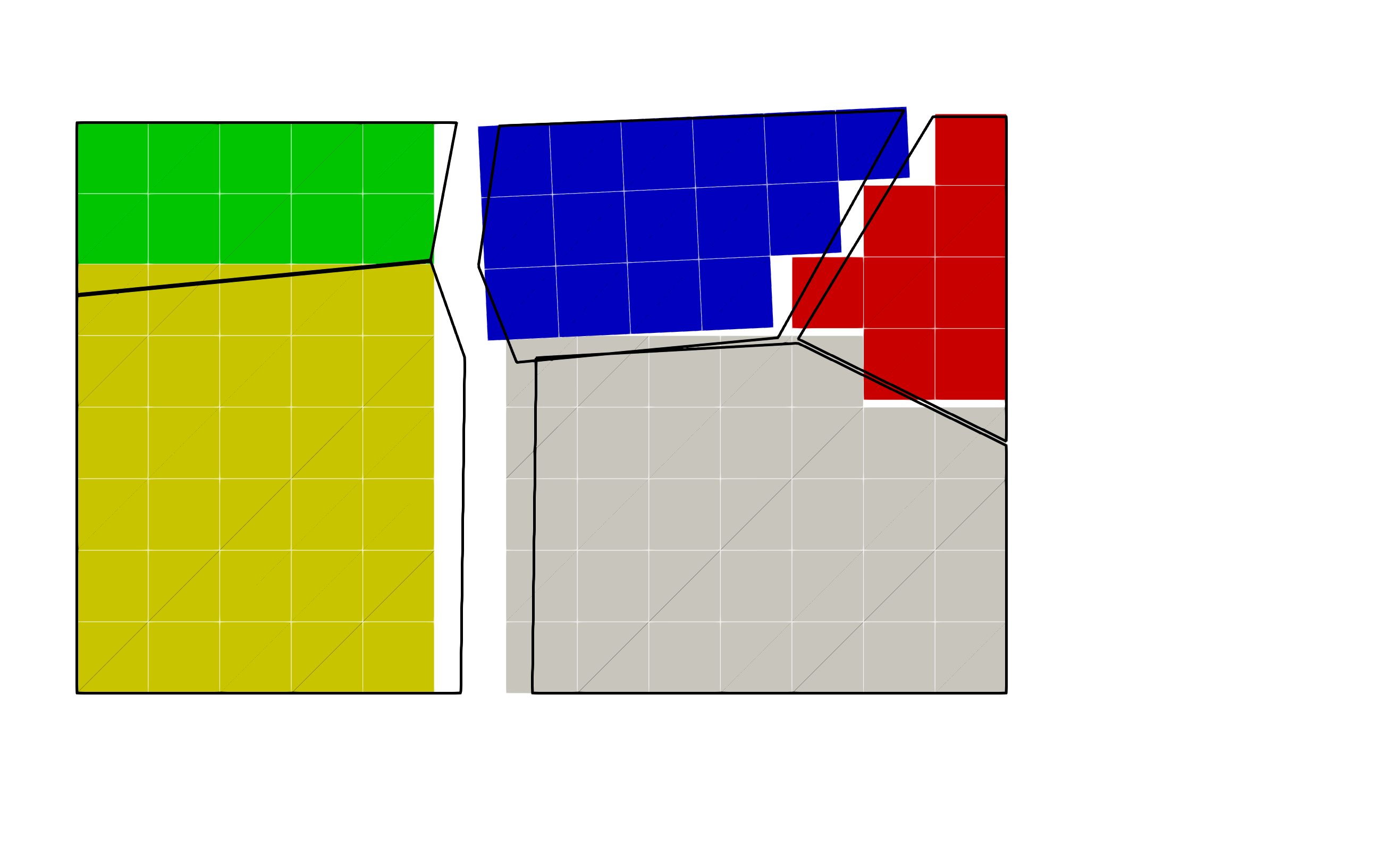}
        \caption{$h=0.1$}
        \label{fig:ny8_t10}
    \end{subfigure}
    \hfill
    \begin{subfigure}{0.32\textwidth}
        \includegraphics[width=\linewidth,trim = {50 200 100 200}, clip]{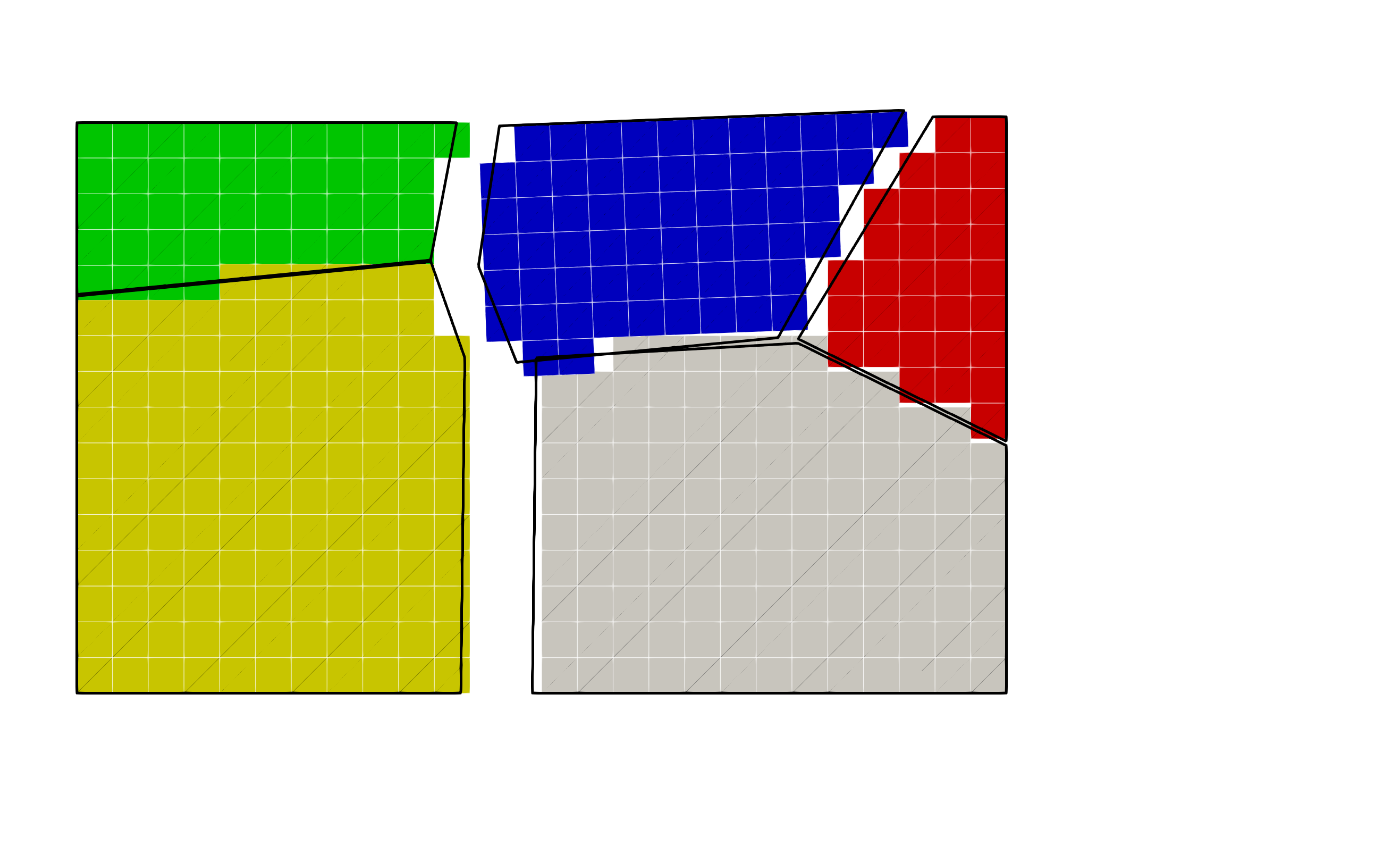}
        \caption{$h=0.05$}
        \label{fig:ny16_t10}
    \end{subfigure}

    \vspace{1em}

    \begin{subfigure}{0.32\textwidth}
        \includegraphics[width=\linewidth,trim = {50 200 100 200}, clip]{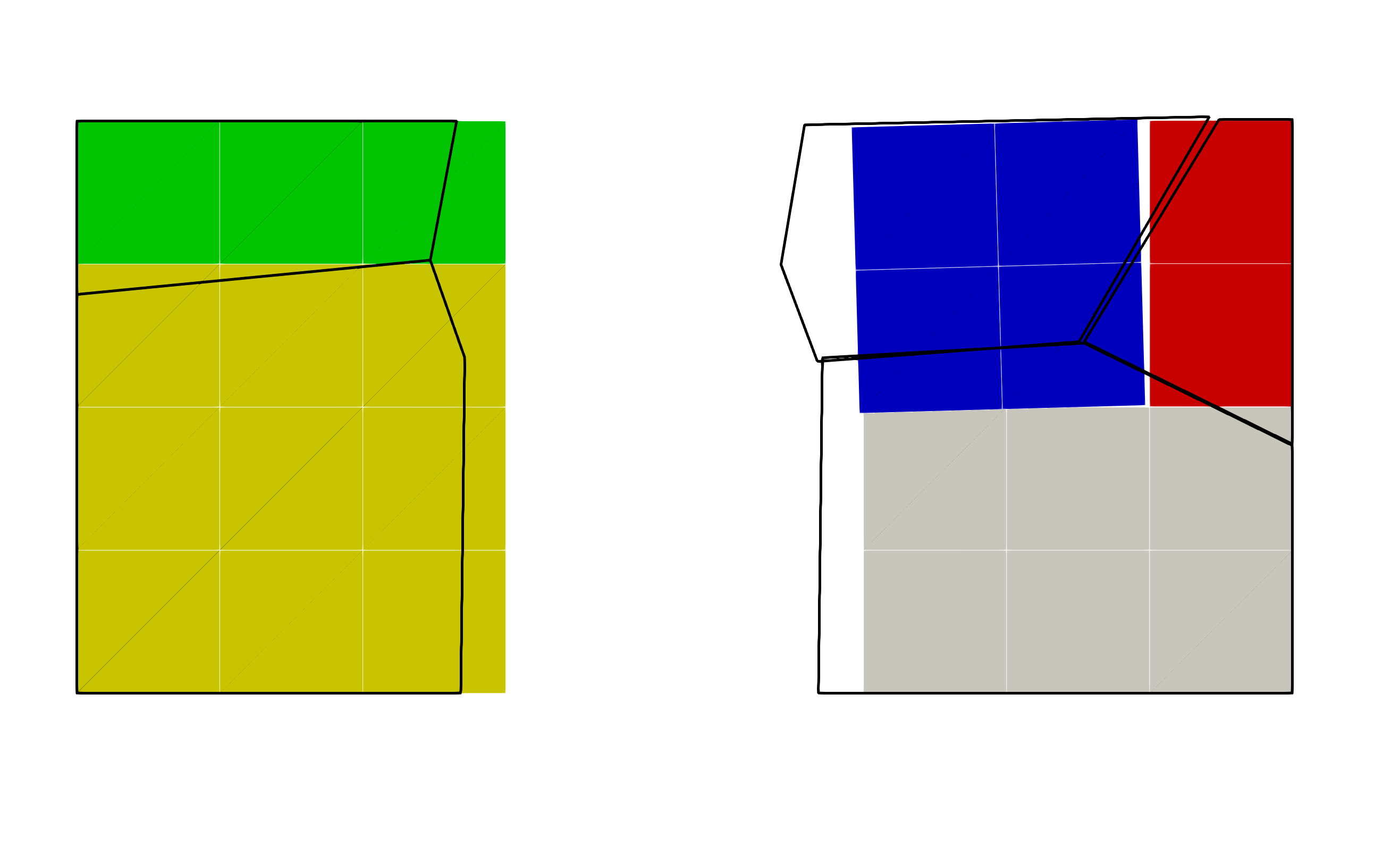}
        \caption{$h=0.2$}
        \label{fig:ny4_t50}
    \end{subfigure}
    \hfill
    \begin{subfigure}{0.32\textwidth}
            \includegraphics[width=\linewidth,trim = {50 200 100 200}, clip]{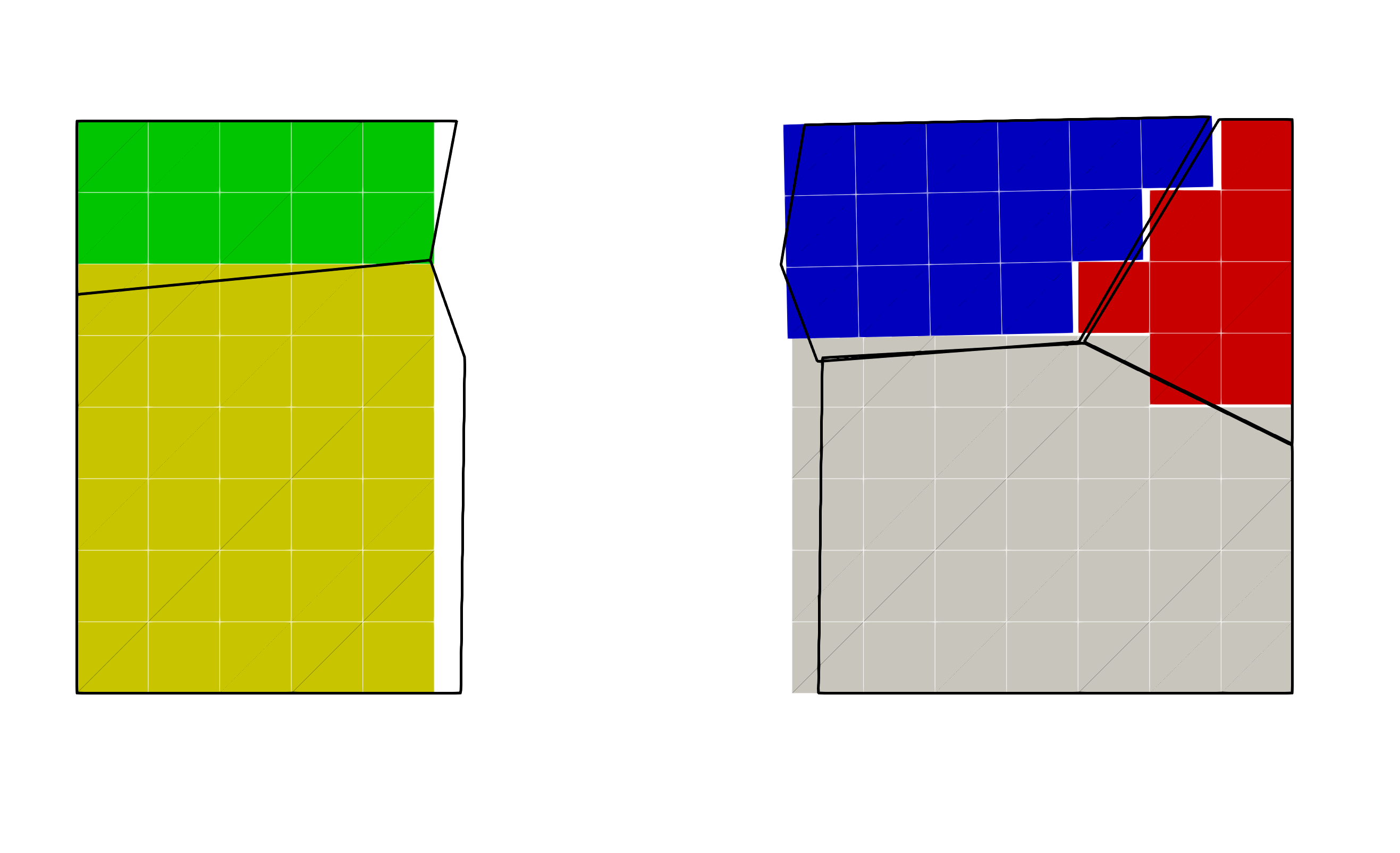}
        \caption{$h=0.1$}
        \label{fig:ny8_t50}
    \end{subfigure}
    \hfill
    \begin{subfigure}{0.32\textwidth}
    \includegraphics[width=\linewidth,trim = {50 200 100 200}, clip]
    {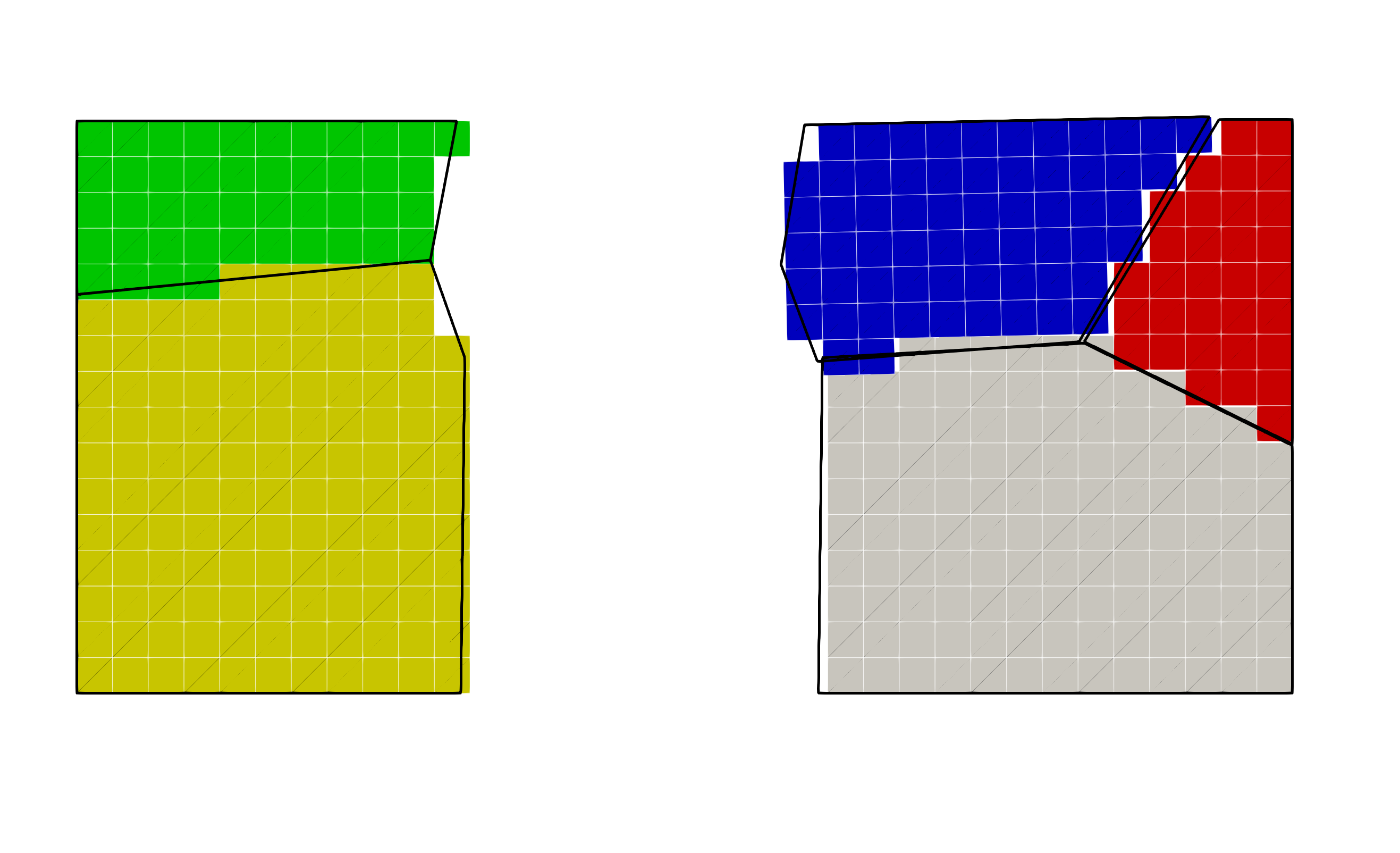}
        \caption{$h=0.05$}
        \label{fig:ny16_t50}
    \end{subfigure}

    \caption{Evolution of the polycrystalline response at two representative time instances, $t=10$ in the first row and $t=50$ in the second row, for different mesh resolutions (\secref{sec:multi_neml2_2D}). From left to right, the results correspond to $h=0.2$, $h=0.1$, and $h=0.05$. The black lines denote the reference interface-fitted solution.}
    \label{fig:sczm_mesh_time_comparison}
\end{figure}

To examine spatial accuracy, \figref{fig:sczm_mesh_time_comparison} compares deformation fields at selected time instances across mesh resolutions. SCZM reproduces the reference response with increasing fidelity as the mesh is refined, despite the lack of interface conformity.

\begin{figure}[htb]
\centering

\begin{subfigure}{0.32\linewidth}
\centering
\begin{tikzpicture}
\begin{axis}[
    width=\linewidth,
    height=0.7\linewidth,
    xlabel={Arc length},
    ylabel={$|t_x|$},
    grid=major,
    major grid style={line width=0.25pt, draw=gray!30},
    tick align=outside,
    tick label style={font=\normalsize},
    label style={font=\normalsize},
    every axis plot/.append style={line width=1.4pt},
]

\addplot[
    color=red
]
table[
    x=arc_length,
    y=tx,
    col sep=space,
    restrict expr to domain={\thisrow{dataset_id}}{2:2}
]{stress_norm_comparison.txt};

\addplot[
    color=black,
    dashed
]
table[
    x=arc_length,
    y=tx,
    col sep=space,
    restrict expr to domain={\thisrow{dataset_id}}{1:1}
]{stress_norm_comparison.txt};

\end{axis}
\end{tikzpicture}
\caption{$|t_x|$}
\end{subfigure}
\hfill
\begin{subfigure}{0.32\linewidth}
\centering
\begin{tikzpicture}
\begin{axis}[
    width=\linewidth,
    height=0.7\linewidth,
    xlabel={Arc length},
    ylabel={$|t_y|$},
    grid=major,
    major grid style={line width=0.25pt, draw=gray!30},
    tick align=outside,
    tick label style={font=\normalsize},
    label style={font=\normalsize},
    every axis plot/.append style={line width=1.4pt},
]

\addplot[color=red]
table[
    x=arc_length,
    y=ty,
    col sep=space,
    restrict expr to domain={\thisrow{dataset_id}}{2:2}
]{stress_norm_comparison.txt};

\addplot[color=black, dashed]
table[
    x=arc_length,
    y=ty,
    col sep=space,
    restrict expr to domain={\thisrow{dataset_id}}{1:1}
]{stress_norm_comparison.txt};

\end{axis}
\end{tikzpicture}
\caption{$|t_y|$}
\end{subfigure}
\hfill
\begin{subfigure}{0.32\linewidth}
\centering
\begin{tikzpicture}
\begin{axis}[
    width=\linewidth,
    height=0.7\linewidth,
    xlabel={Arc length},
    ylabel={$|t_z|$},
    grid=major,
    major grid style={line width=0.25pt, draw=gray!30},
    tick align=outside,
    tick label style={font=\normalsize},
    label style={font=\normalsize},
    every axis plot/.append style={line width=1.4pt},
]

\addplot[color=red]
table[
    x=arc_length,
    y=tz,
    col sep=space,
    restrict expr to domain={\thisrow{dataset_id}}{2:2}
]{stress_norm_comparison.txt};

\addplot[color=black, dashed]
table[
    x=arc_length,
    y=tz,
    col sep=space,
    restrict expr to domain={\thisrow{dataset_id}}{1:1}
]{stress_norm_comparison.txt};

\end{axis}
\end{tikzpicture}
\caption{$|t_z|$}
\end{subfigure}

\vspace{0.5em}

\begin{tikzpicture}
\begin{axis}[
    hide axis,
    xmin=0, xmax=1,
    ymin=0, ymax=1,
    legend columns=2,
    legend style={
        draw=none,
        fill=none,
        font=\footnotesize
    }
]

\addlegendimage{red, line width=1.4pt}
\addlegendentry{SCZM ($h=0.003125$)}

\addlegendimage{black, dashed, line width=1.4pt}
\addlegendentry{IFM ($h_{\Omega}=0.00292$)}

\end{axis}
\end{tikzpicture}

\caption{Comparison of the traction magnitude $|\boldsymbol{t}_i|$ along the diagonal of the domain in the reference configuration at $t=30$ (\secref{sec:multi_neml2_2D}). The non-zero out-of-plane traction component arises from the use of a fully three-dimensional crystal plasticity constitutive model, even though the problem geometry is two-dimensional.}
\label{fig:2D_grains_cauchy_stress}
\end{figure}
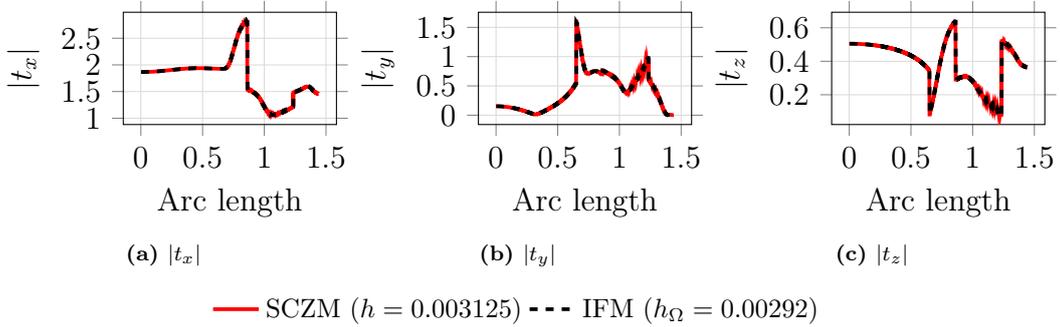

We further assess stress accuracy by evaluating traction magnitudes $|\bs{t}_i| = |\bs{\sigma}\cdot\bs{e}_i|$ along the domain diagonal at $t=30$ (\figref{fig:2D_grains_cauchy_stress}). The SCZM results are in close agreement with the interface-fitted solution, indicating accurate stress transmission across grain boundaries.

\begin{figure}[htb]
\centering

\begin{subfigure}{0.48\textwidth}
    \centering
    \includegraphics[width=\linewidth,trim = {60 200 370 200}, clip]{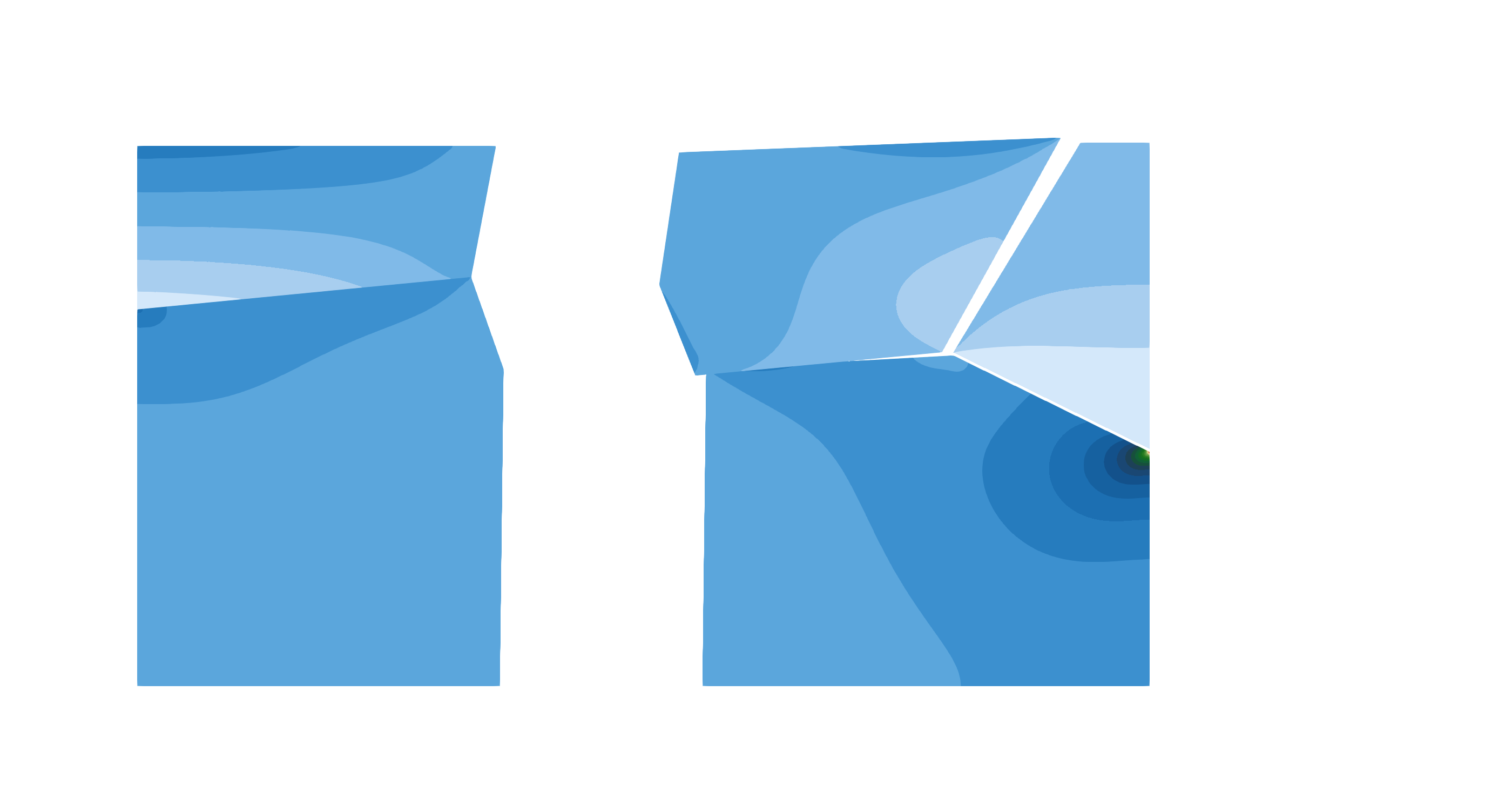}
    \caption{IFM}
\end{subfigure}
\hfill
\begin{subfigure}{0.48\textwidth}
    \centering
    \includegraphics[width=\linewidth,trim = {60 200 370 200}, clip]{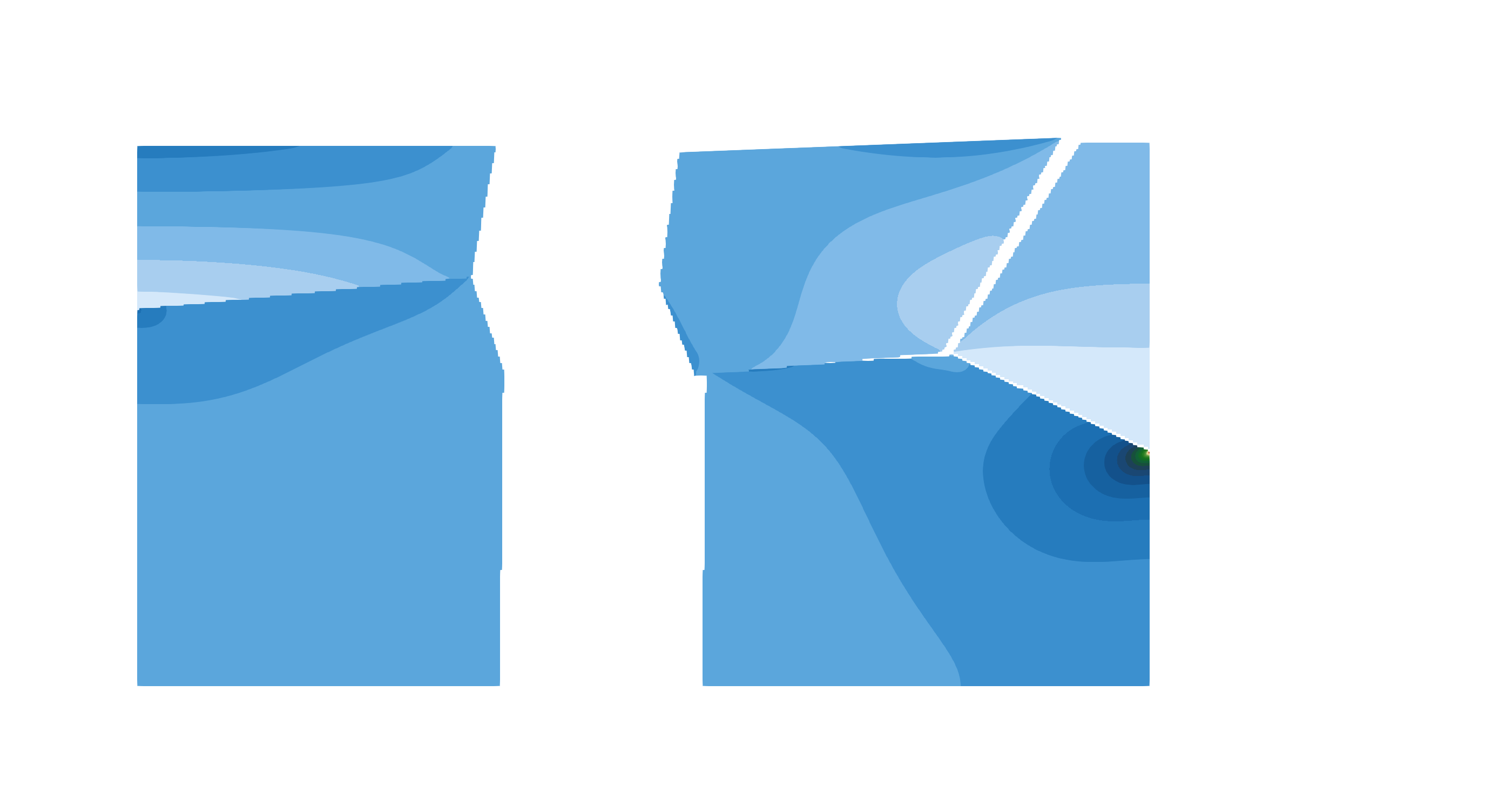}
    \caption{SCZM}
\end{subfigure}

\vspace{0.5em}

\begin{subfigure}{0.48\textwidth}
    \centering
    \includegraphics[width=\linewidth,trim = {0 60 250 60}, clip]{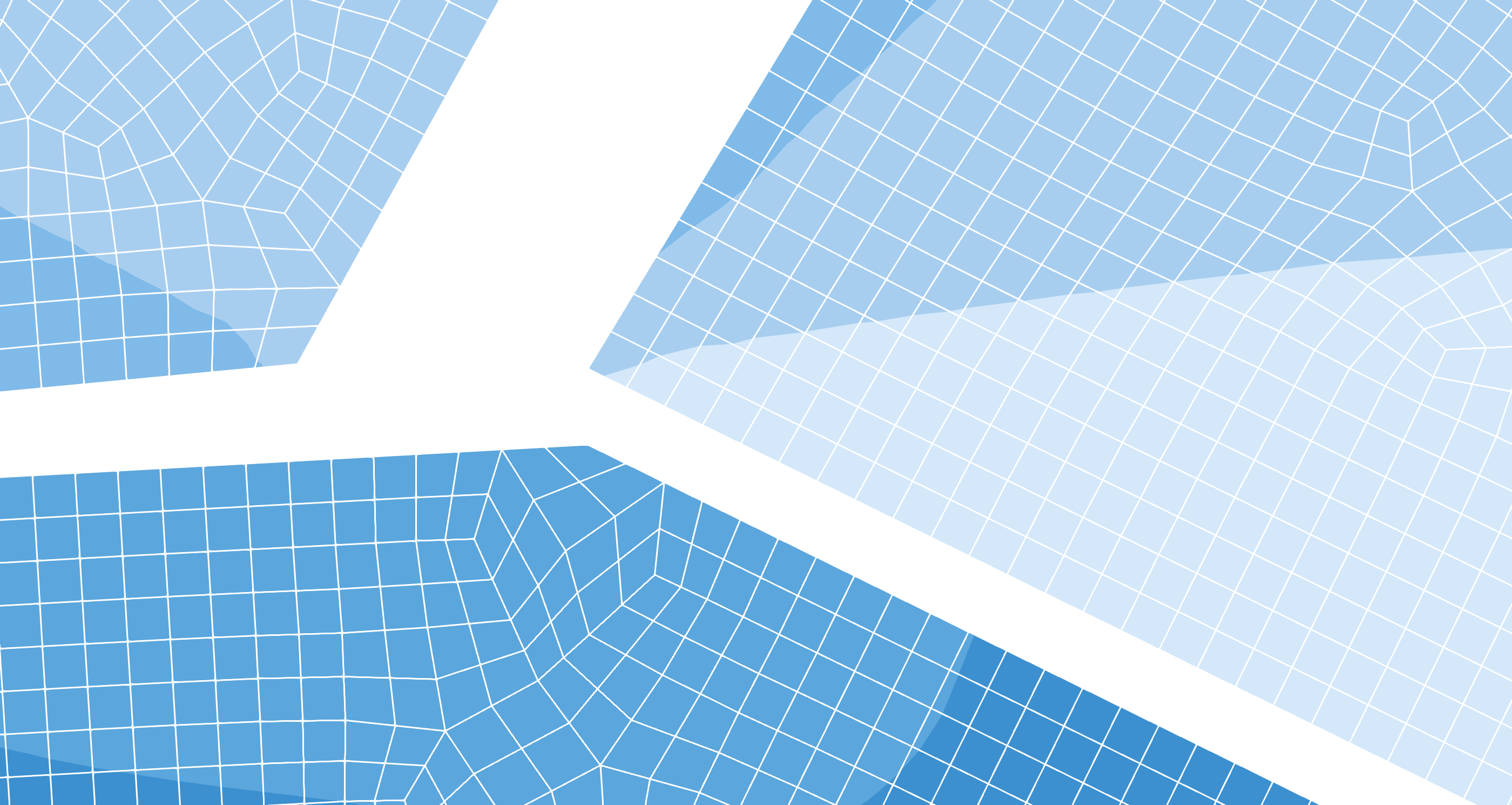}
    \caption{IFM (triple junction)}
\end{subfigure}
\hfill
\begin{subfigure}{0.48\textwidth}
    \centering
    \includegraphics[width=\linewidth,trim = {0 60 250 60}, clip]{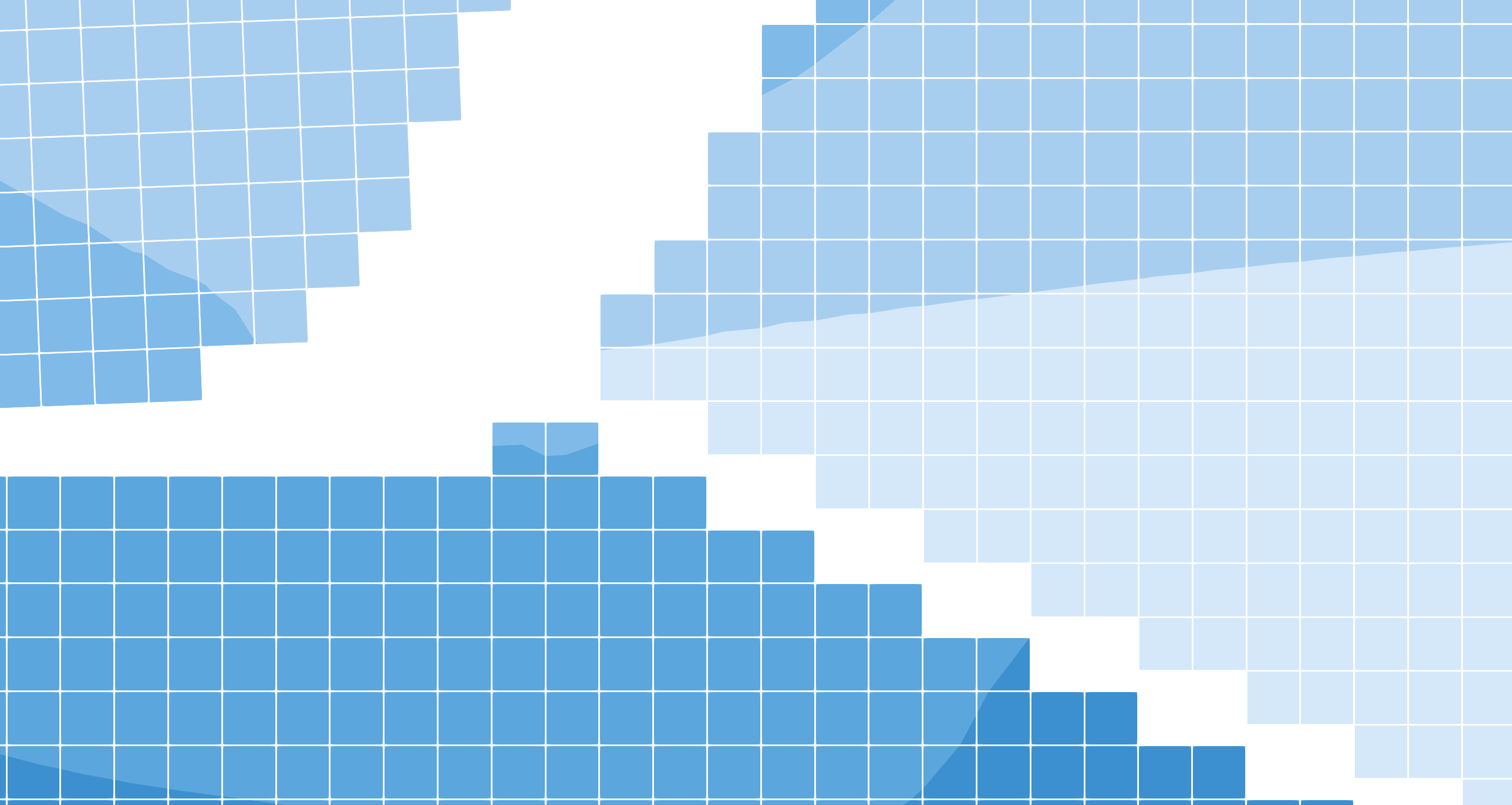}
    \caption{SCZM (triple junction)}
\end{subfigure}

\vspace{0.5em}

\begin{subfigure}{0.48\textwidth}
    \centering
    \includegraphics[width=\linewidth,trim = {0 60 0 60}, clip]{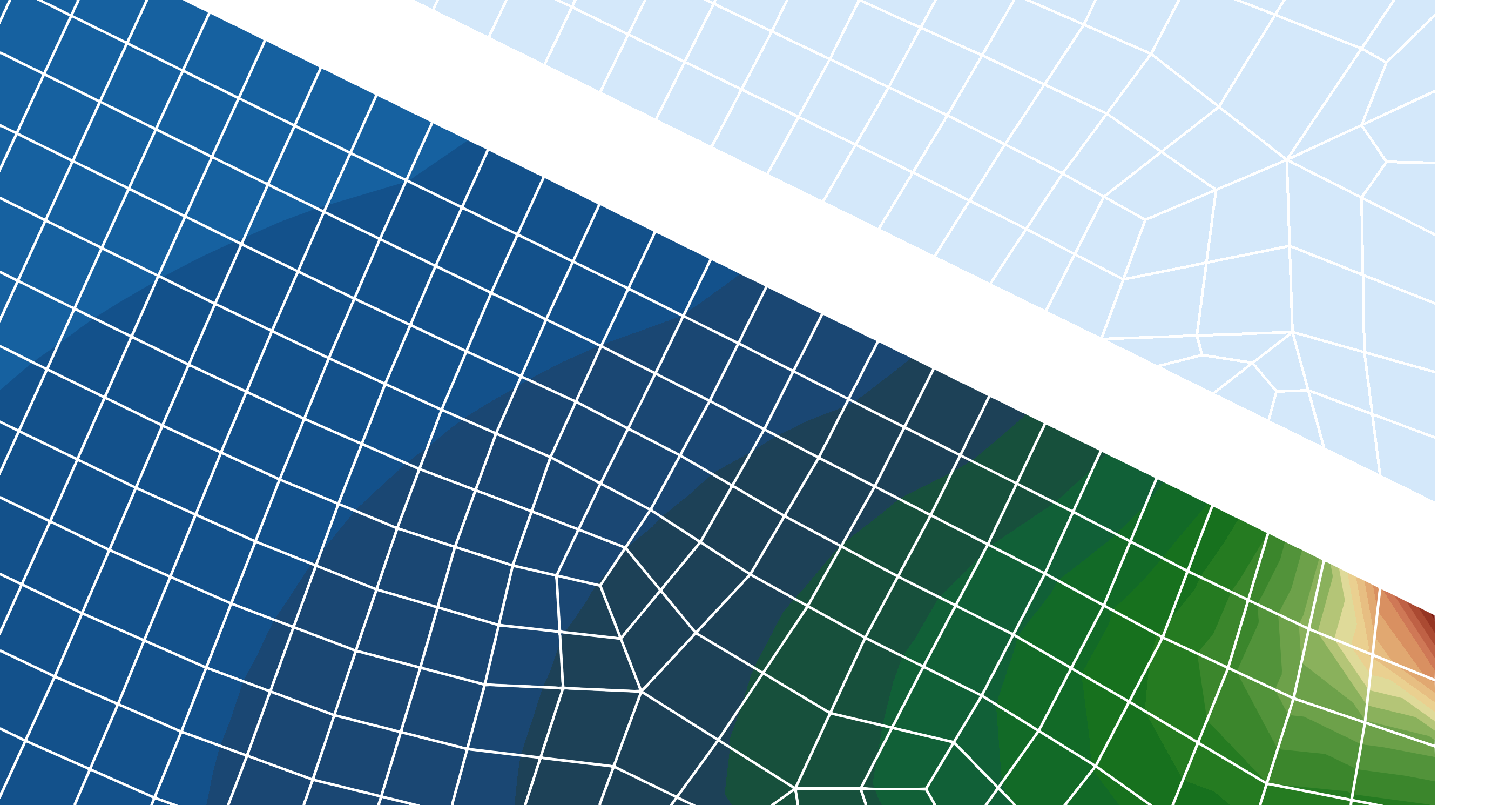}
    \caption{IFM (stress concentration)}
\end{subfigure}
\hfill
\begin{subfigure}{0.48\textwidth}
    \centering
    \includegraphics[width=\linewidth,trim = {0 60 0 60}, clip]{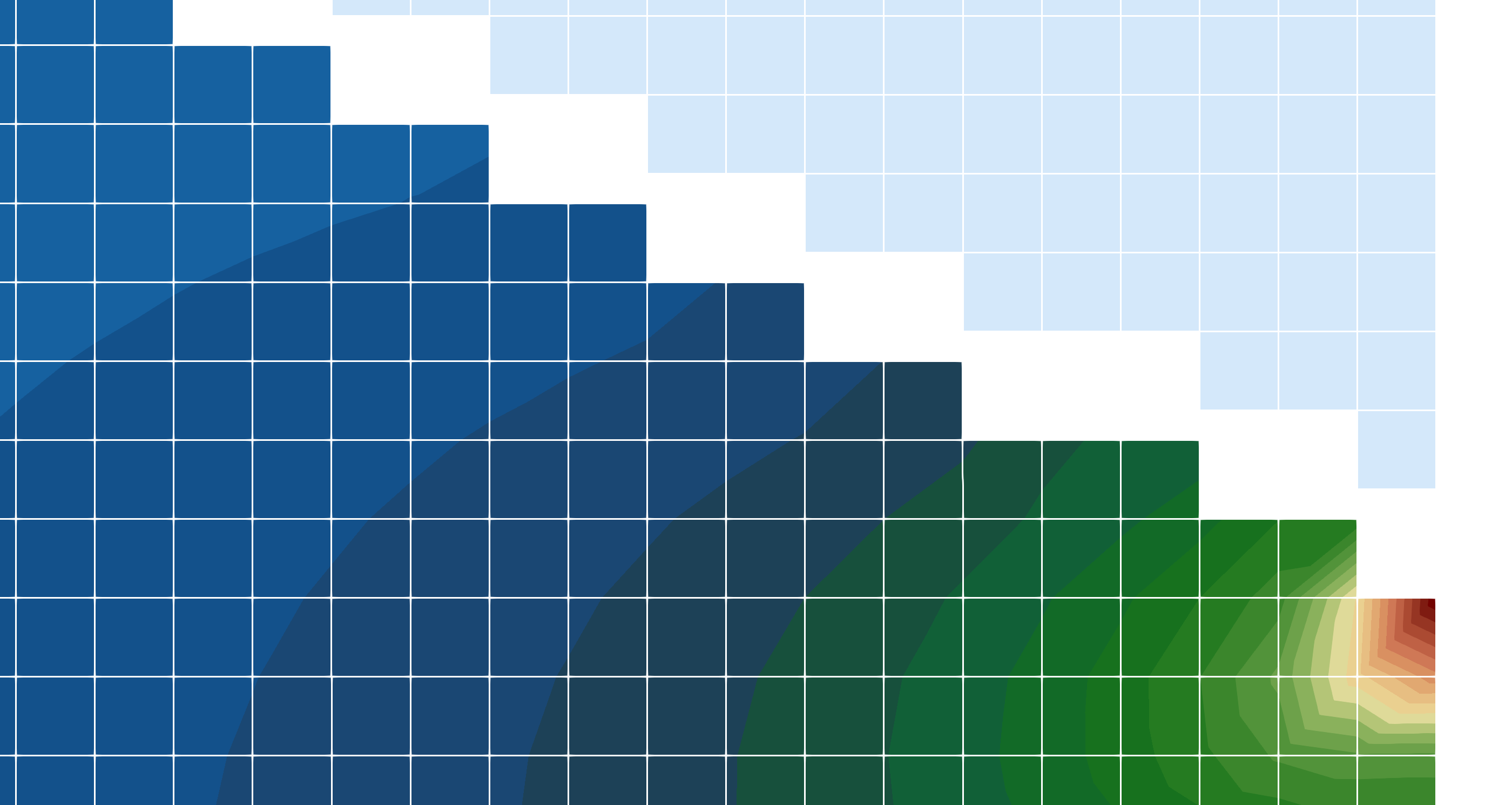}
    \caption{SCZM (stress concentration)}
\end{subfigure}

\includegraphics[width=0.15\linewidth]{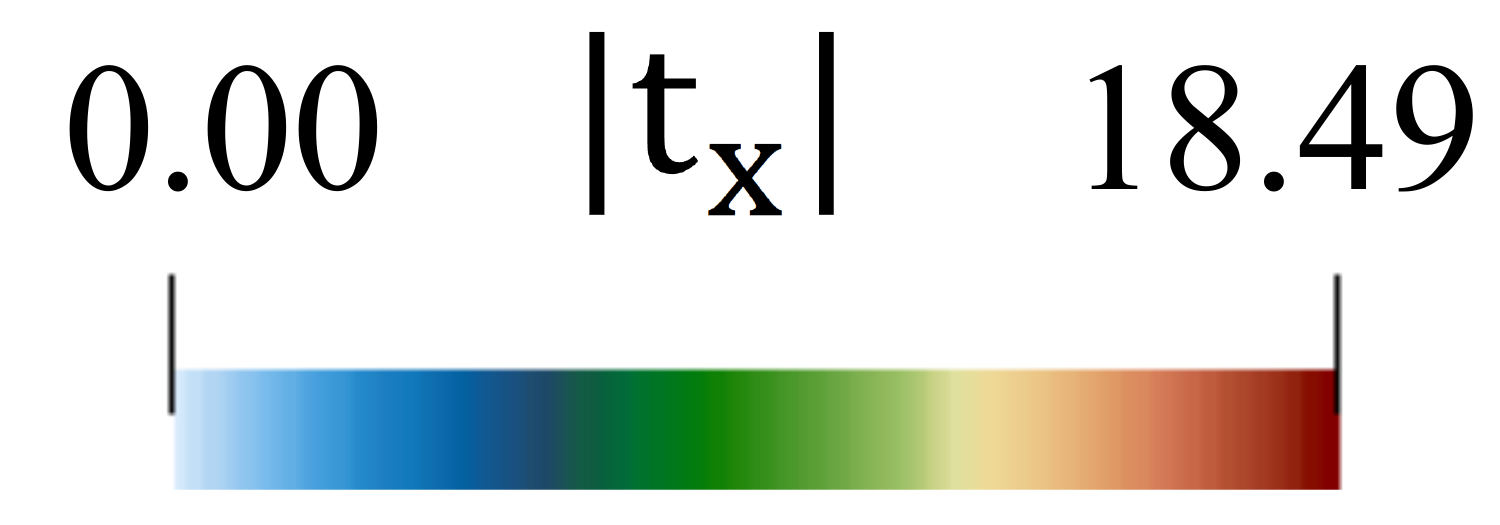}

\caption{
Comparison of traction magnitude $|\boldsymbol{t}_x|$ at $t=30$ between the IFM and SCZM (\secref{sec:multi_neml2_2D}).
The first row shows the full-field distribution, the second row presents a zoom-in near the triple-junction region, and the third row highlights the stress concentration region.
SCZM accurately reproduces both global and local features of the solution without requiring an interface-fitted mesh.
}
\label{fig:stress_comparison_t30}
\end{figure}

Full-field comparisons of $|\bs{t}_x|$ at $t=30$ (\figref{fig:stress_comparison_t30}) show that SCZM captures both the global stress distribution and localized features, including stress concentrations. In particular, zoomed-in views near triple-junction regions demonstrate that SCZM accurately resolves complex multi-interface interactions. The location and magnitude of peak stresses are also well reproduced, confirming that SCZM can capture localized stress concentrations without requiring interface-fitted meshes.

\subsection{Crystal plasticity on 3D polycrystalline RVE} \label{sec:multi_neml2_3D}

We consider a three-dimensional polycrystalline RVE to assess the robustness and scalability of SCZM in a fully realistic setting. The problem involves complex grain geometries, three-dimensional interface networks, and history-dependent crystal plasticity with the same parameters as in \secref{sec:single_neml2_2D}.

\begin{figure}[htb]
    \centering

    \begin{subfigure}{0.49\textwidth}
        \centering
        \includegraphics[width=\linewidth,trim = {800 200 800 200}, clip]{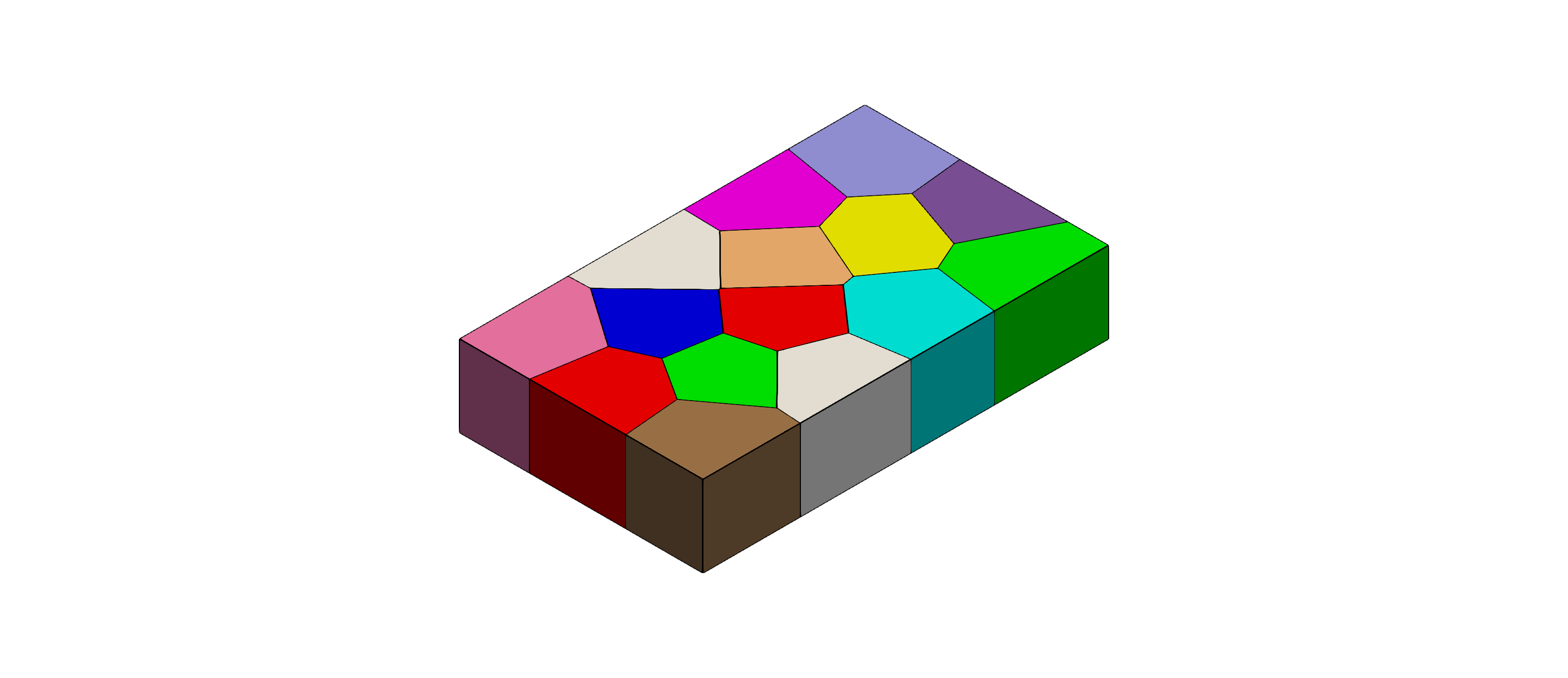}
        \caption{Grain geometry}
        \label{fig:3D_geometry}
    \end{subfigure}
    \hfill
    \begin{subfigure}{0.49\textwidth}
        \centering
        \includegraphics[width=\linewidth,trim = {800 200 800 200}, clip]{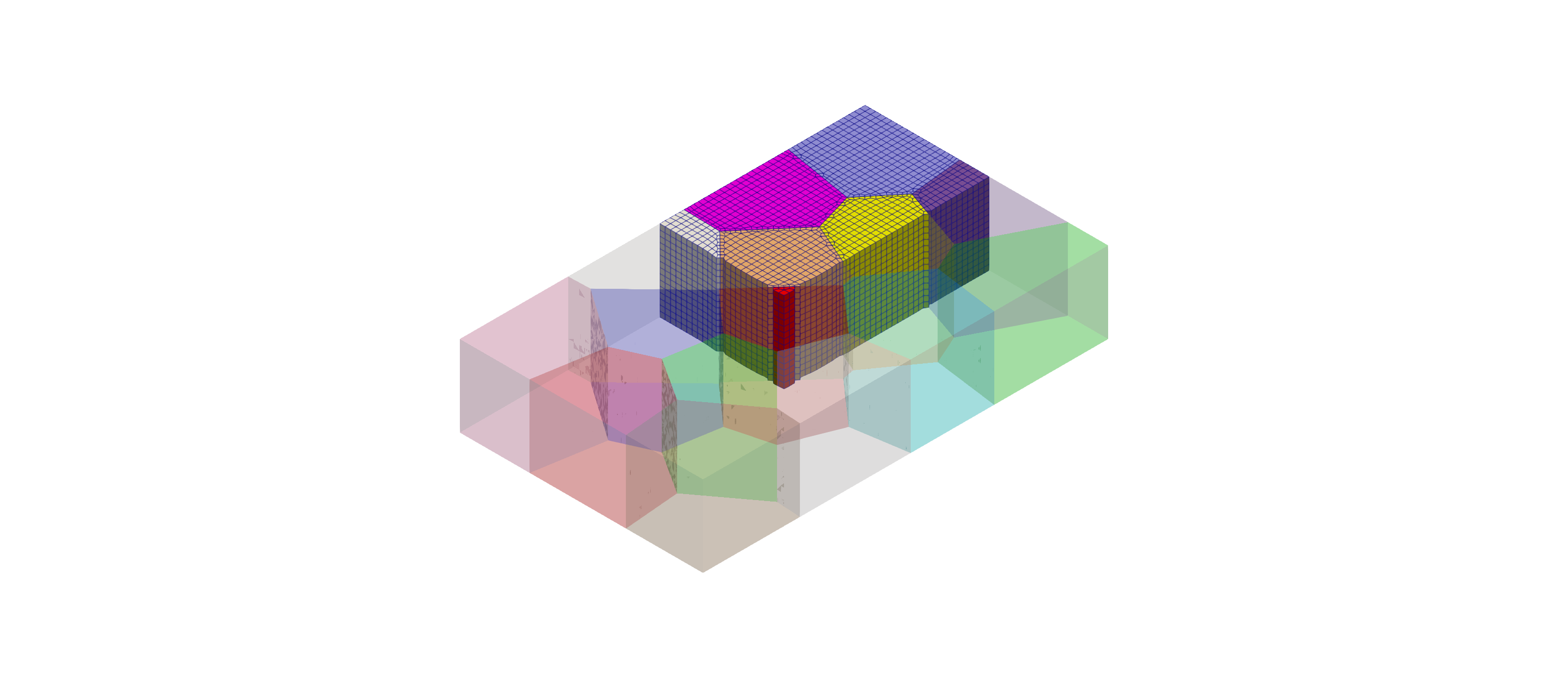}
        \caption{IFM: Interface-fitted mesh}
        \label{fig:3D_IFM}
    \end{subfigure}

    \vspace{0.5em}

    \begin{subfigure}{0.49\textwidth}
        \centering
        \includegraphics[width=\linewidth,trim = {800 200 800 200}, clip]{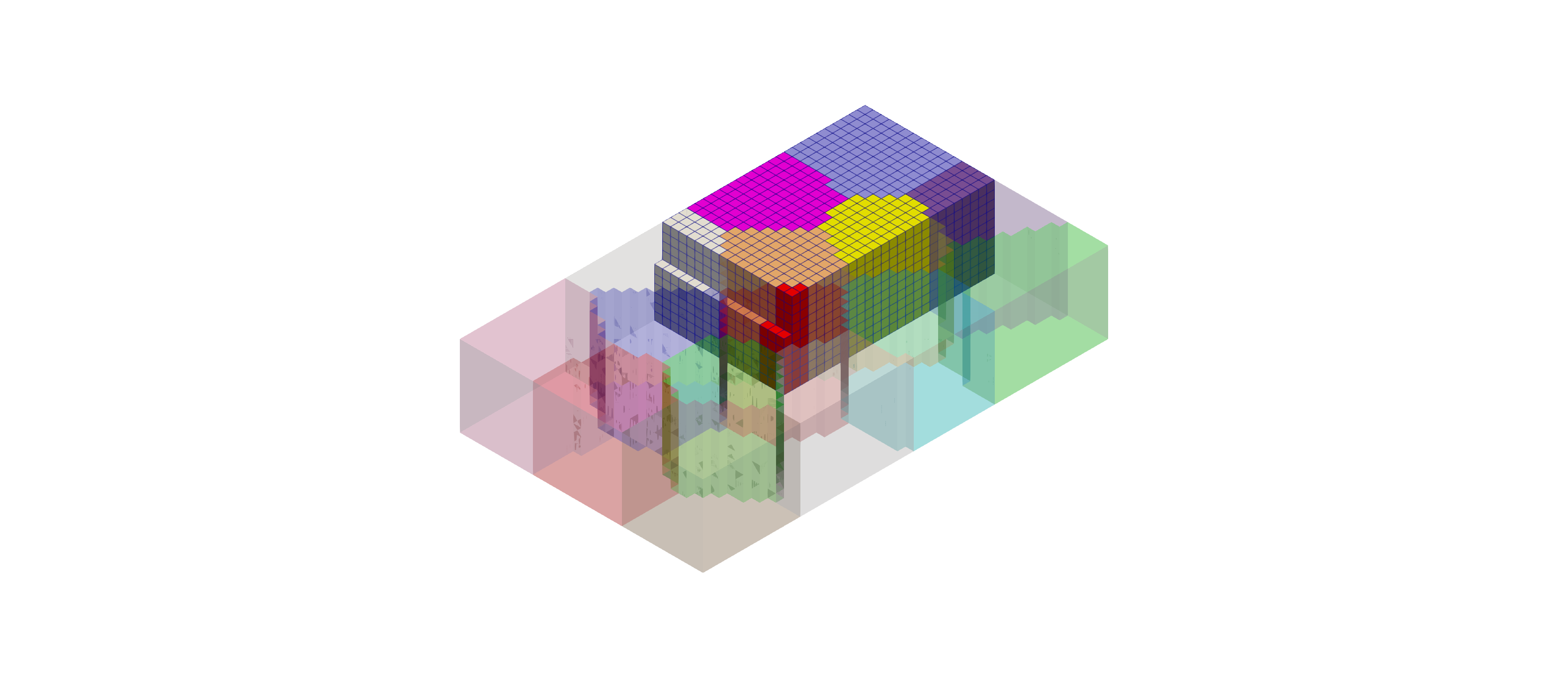}
        \caption{SCZM: Non-interface-fitted hexahedral mesh}
        \label{fig:3D_SCZM}
    \end{subfigure}

    \caption{Comparison of interface representations in the 3D polycrystalline domain (\secref{sec:multi_neml2_3D}).
    IFM resolves interfaces explicitly, while SCZM employs a surrogate non-interface-fitted discretization.}

    \label{fig:3D_grain_IFM_SCZM}
\end{figure}

The computational domain $[-750,750]\times[-150,150]\times[-450,450]$ is partitioned into 15 grains (\figref{fig:3D_geometry}). Results are compared between an interface-fitted hexahedral discretization (IFM, generated using \sculpt{}) and a non-interface-fitted hexahedral mesh using SCZM (\figref{fig:3D_IFM}, \figref{fig:3D_SCZM}). Boundary conditions fix displacements on the $x^{-}$, $y^{-}$, and $z^{-}$ faces, while a displacement $u_x = 0.5\,t$ is applied on the $x^{+}$ face.

\begin{figure}[htb]
    \centering

    \begin{subfigure}{0.32\textwidth}
        \includegraphics[width=\textwidth,trim={700 200 500 120},clip]{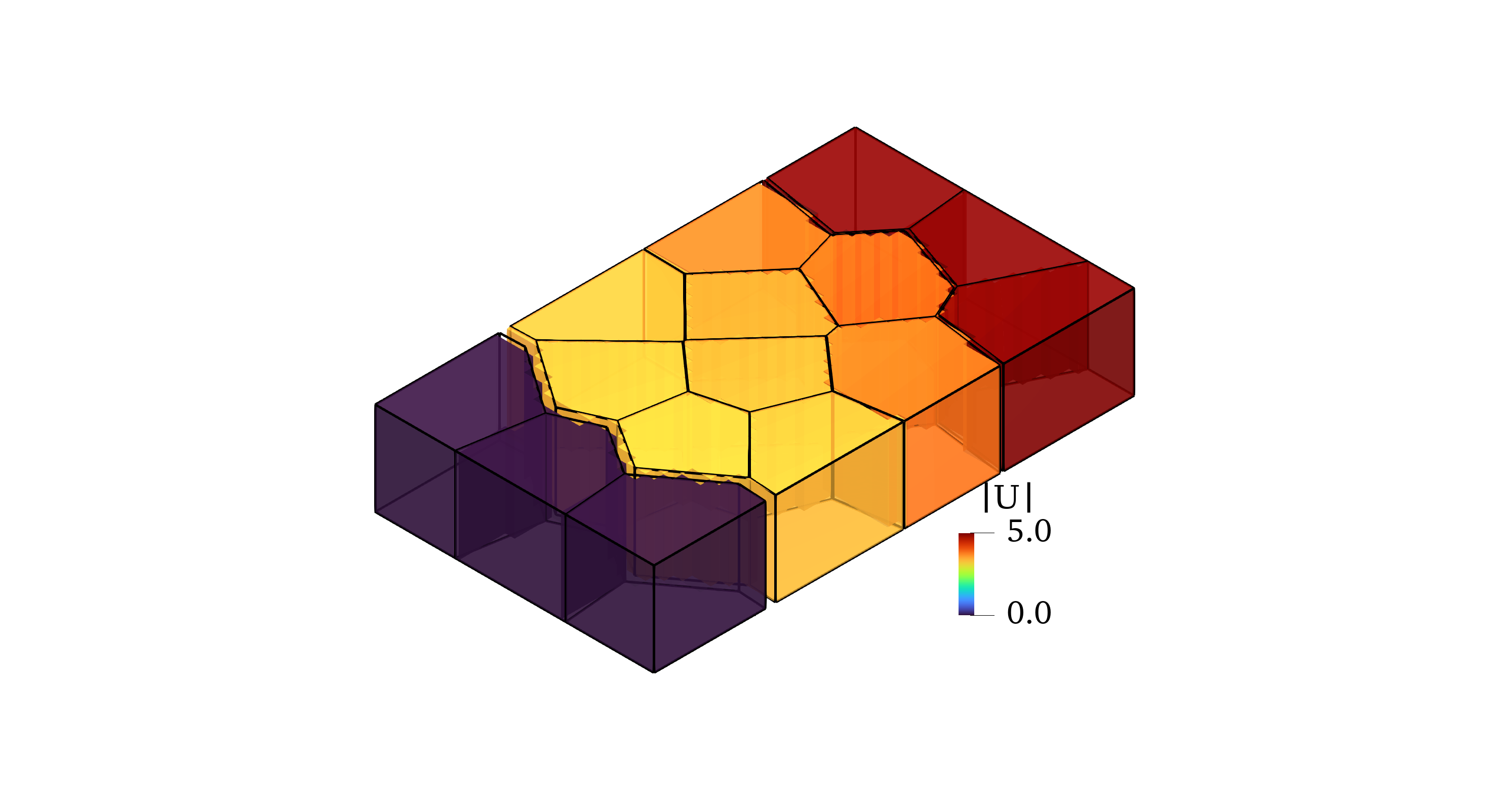}
        \caption{$t=10$}
        \label{fig:GI50_t10}
    \end{subfigure}
    \hfill
    \begin{subfigure}{0.32\textwidth}
        \includegraphics[width=\textwidth,trim={700 200 500 120},clip]{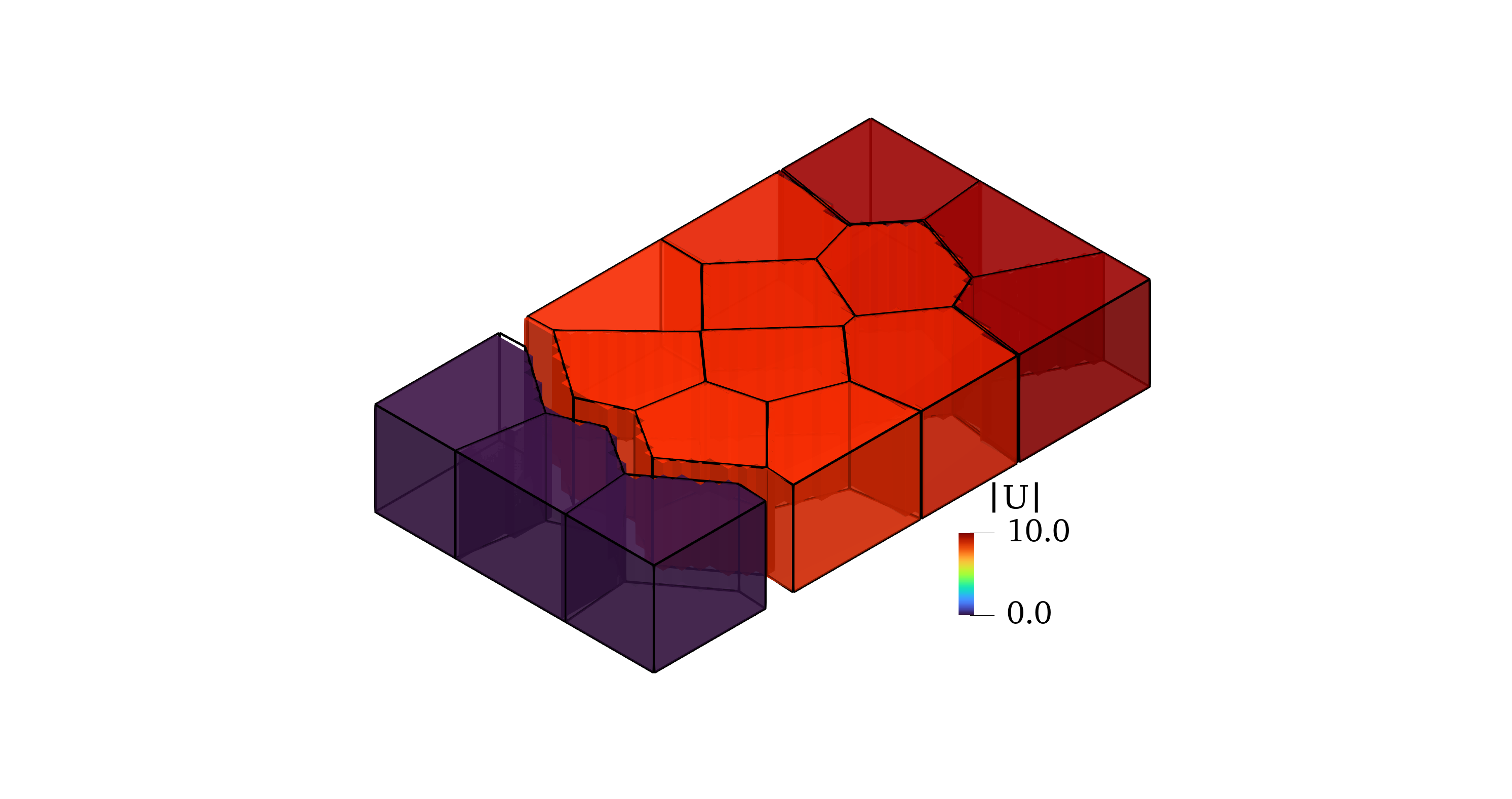}
        \caption{$t=20$}
        \label{fig:GI50_t20}
    \end{subfigure}
    \hfill
    \begin{subfigure}{0.32\textwidth}
        \includegraphics[width=\textwidth,trim={700 200 500 120},clip]{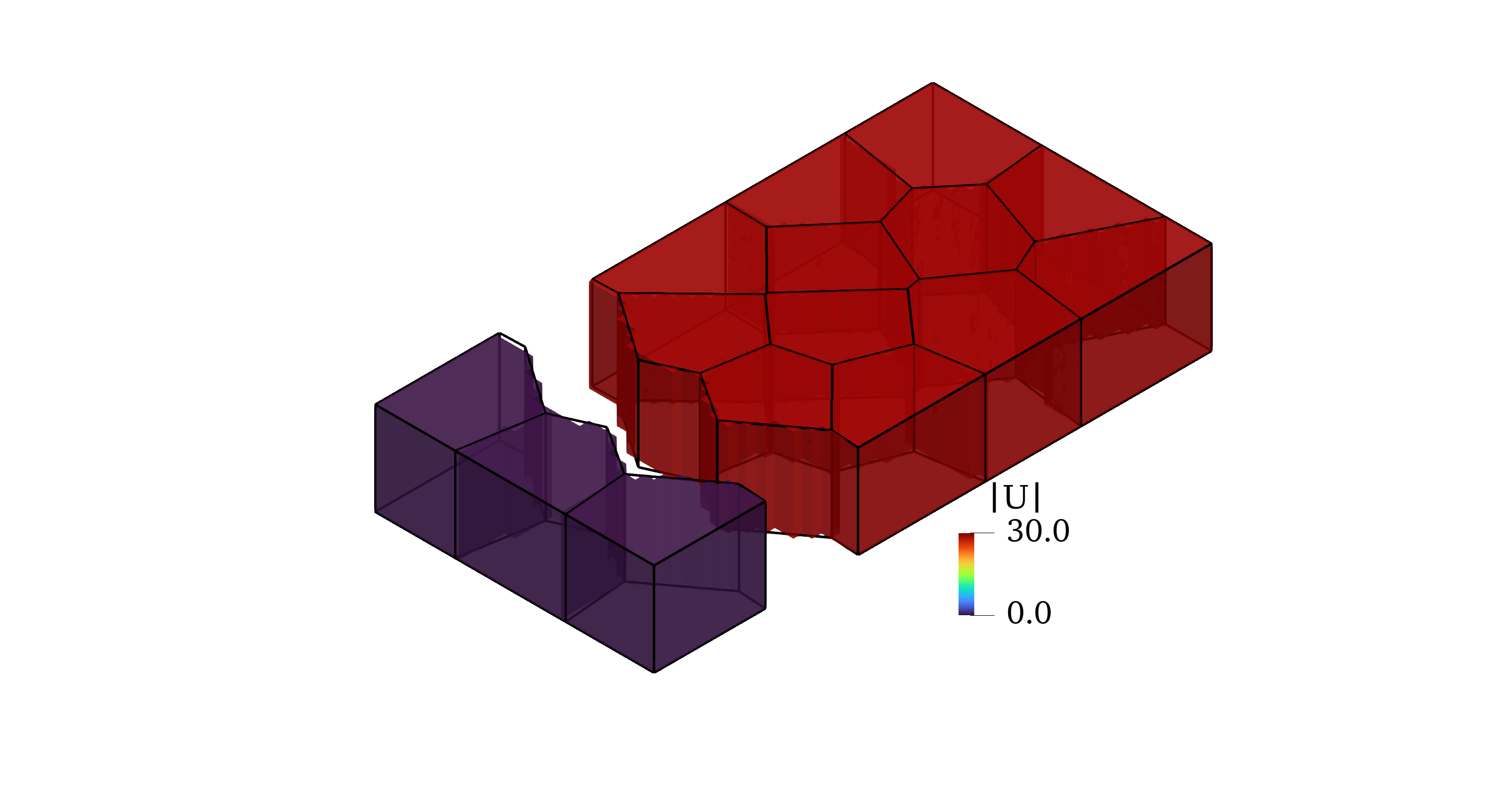}
        \caption{$t=60$}
        \label{fig:GI50_t60}
    \end{subfigure}

    \vspace{1em}

    \begin{subfigure}{0.32\textwidth}
        \includegraphics[width=\linewidth,trim = {700 200 500 120}, clip]{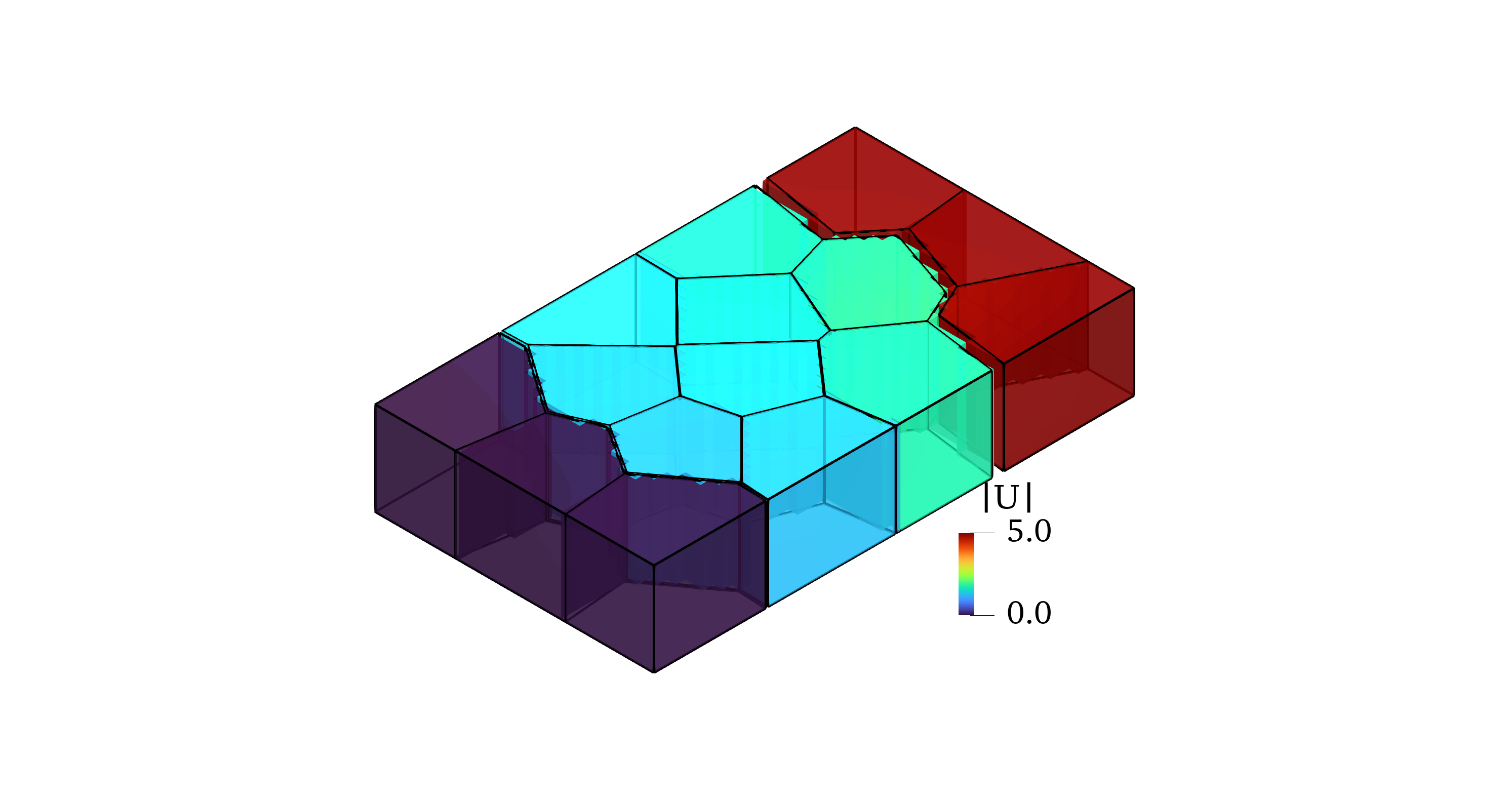}
        \caption{$t=10$}
        \label{fig:GI30_t10}
    \end{subfigure}
    \hfill
    \begin{subfigure}{0.32\textwidth}
        \includegraphics[width=\linewidth,trim = {700 200 500 120}, clip]{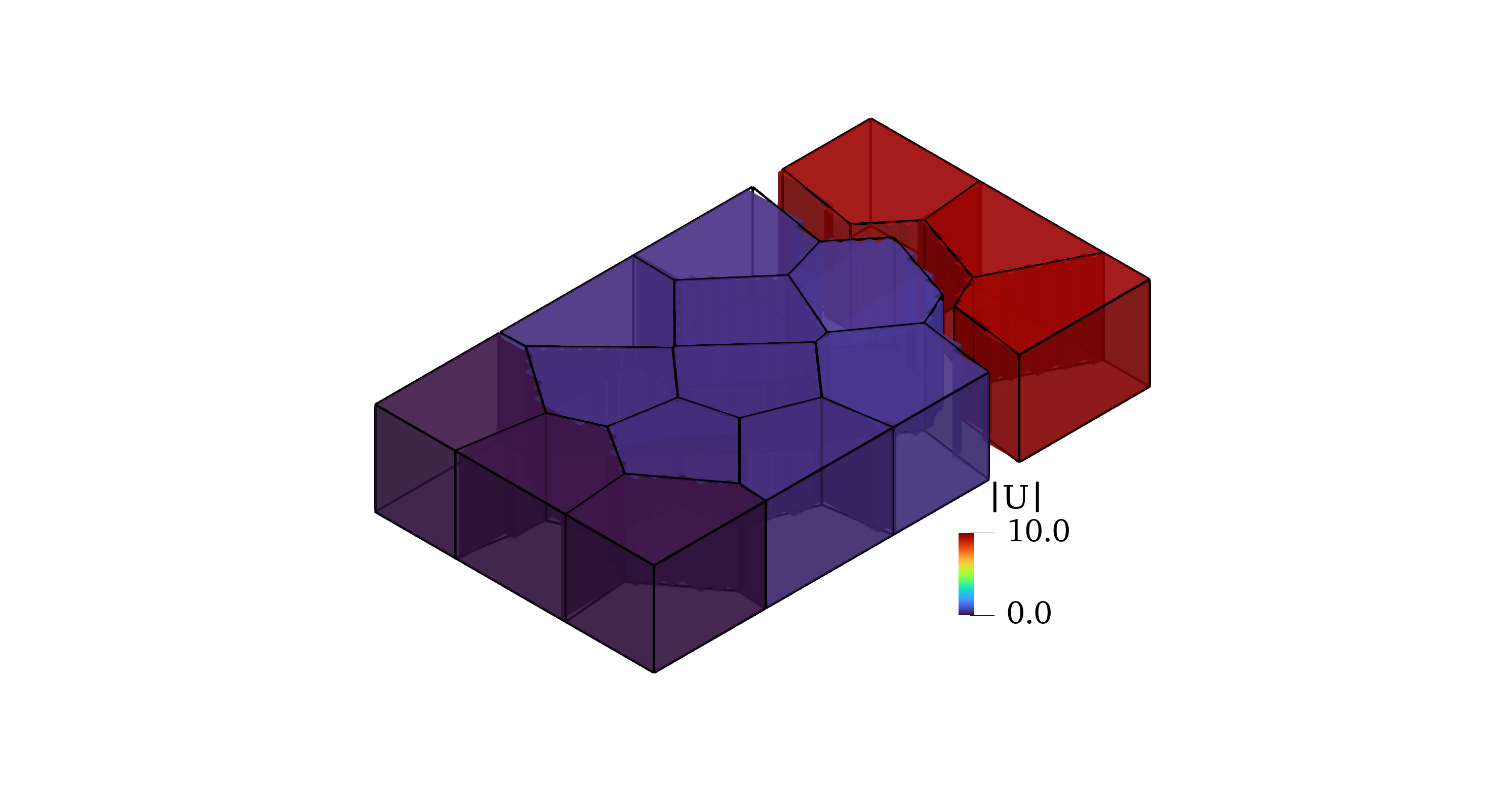}
        \caption{$t=20$}
        \label{fig:GI30_t20}
    \end{subfigure}
    \hfill
    \begin{subfigure}{0.32\textwidth}
        \includegraphics[width=\linewidth,trim = {700 200 500 120}, clip]{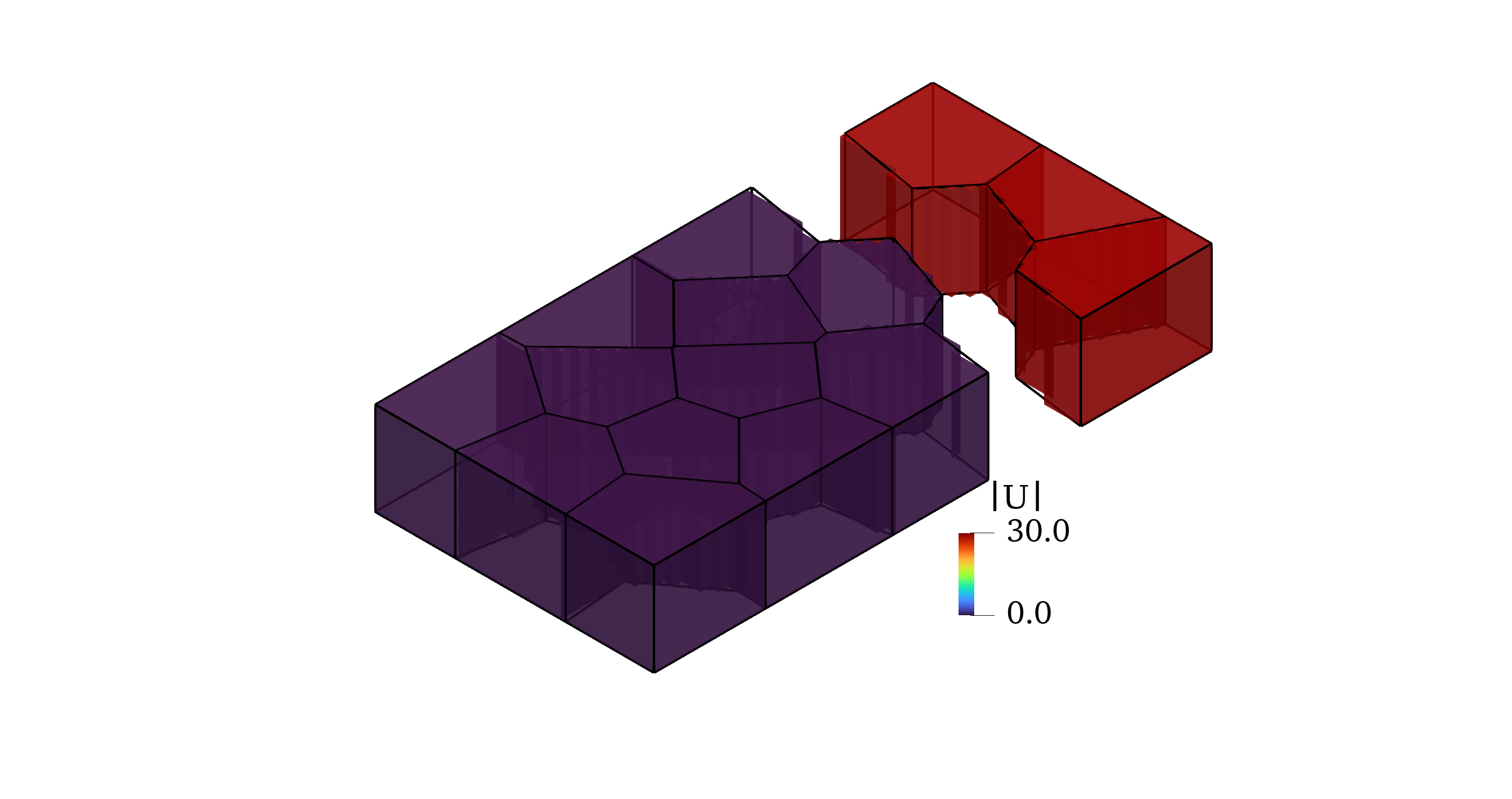}
        \caption{$t=60$}
        \label{fig:GI30_t60}
    \end{subfigure}

    \caption{Comparison of displacement magnitude obtained using SCZM for two sets of interface parameters based on a bilinear mixed-mode traction-separation law at three representative time instants (\secref{sec:multi_neml2_3D}).
    The deformation is visualized with a displacement scaling factor of 10 to enhance visibility. The color contours, however, correspond to the actual (unscaled) displacement magnitude.
    The top row corresponds to $K = 50, G_{I c} = 30$, $G_{II c} = 25$, $N = 3$, and $S = 5$, while the bottom row corresponds to $K = 100, G_{I c} = 30$, $G_{II c} = 200$, $N = 10$, and $S = 20$.
    The black lines indicate the reference interface-fitted configuration under the applied displacement.} \label{fig:sczm_comparison_3D_two_cases}
\end{figure}

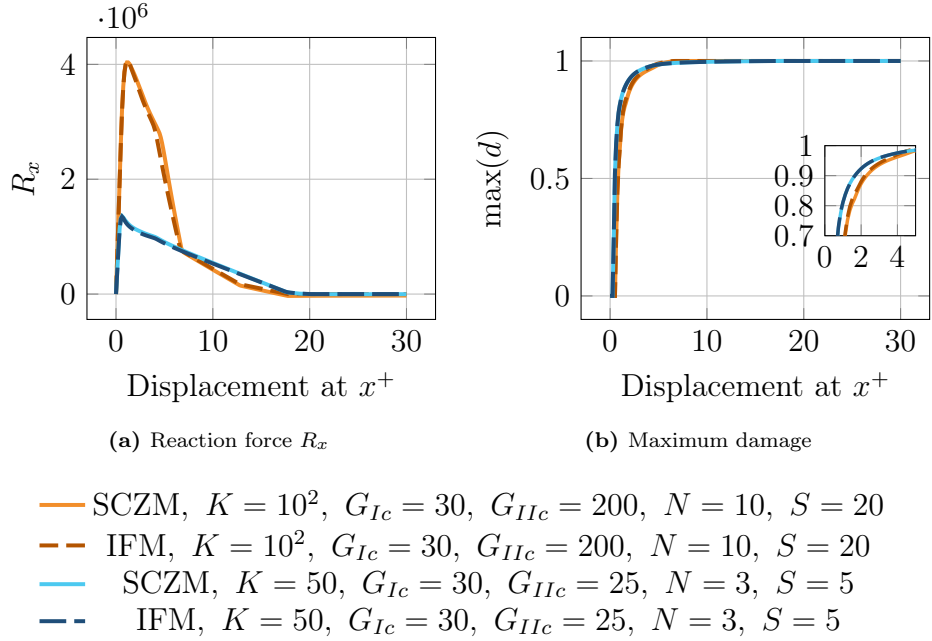
\begin{figure}[htb]
\centering

\begin{subfigure}{0.45\textwidth}
\centering
\begin{tikzpicture}
\begin{axis}[
    width=\linewidth,
    xlabel={Displacement at $x^+$},
    ylabel={$R_x$},
    grid=both
]

\addplot[color=sczmA, line width=1.5pt]
table[x expr=0.5*\thisrow{time}, y=react_y]
{sczm_react_K1e2_GI30_GII200_N10_S20.txt};

\addplot[color=ifmA, line width=1.5pt, dash pattern=on 7pt off 3pt]
table[x expr=0.5*\thisrow{time}, y=react_y]
{ifm_react_K1e2_GI30_GII200_N10_S20.txt};

\addplot[color=sczmB, line width=1.5pt]
table[x expr=0.5*\thisrow{time}, y=react_y]
{sczm_react_K50_GI30_GII25_N3_S5.txt};

\addplot[color=ifmB, line width=1.5pt, dash pattern=on 10pt off 3pt]
table[x expr=0.5*\thisrow{time}, y=react_y]
{ifm_react_K50_GI30_GII25_N3_S5.txt};

\end{axis}
\end{tikzpicture}
\caption{Reaction force $R_x$}
\end{subfigure}
%
\begin{subfigure}{0.45\textwidth}
\centering
\begin{tikzpicture}
\begin{axis}[
    width=\linewidth,
    xlabel={Displacement at $x^+$},
    ylabel={$\max(d)$},
    grid=both
]

\addplot[color=sczmA, line width=1.5pt]
table[
    x expr=0.5*\thisrow{time},
    y=max_damage
]{sczm_damage_K1e2_GI30_GII200_N10_S20.txt};

\addplot[
    color=ifmA,
    line width=1.5pt,
    dash pattern=on 7pt off 3pt
]
table[
    x expr=0.5*\thisrow{time},
    y=max_damage
]{ifm_damage_K1e2_GI30_GII200_N10_S20.txt};

\addplot[color=sczmB, line width=1.5pt]
table[
    x expr=0.5*\thisrow{time},
    y=max_damage
]{sczm_damage_K50_GI30_GII25_N3_S5.txt};

\addplot[
    color=ifmB,
    line width=1.5pt,
    dash pattern=on 10pt off 3pt
]
table[
    x expr=0.5*\thisrow{time},
    y=max_damage
]{ifm_damage_K50_GI30_GII25_N3_S5.txt};

\end{axis}

\begin{axis}[
    at={(0.52\linewidth,0.18\linewidth)},
    anchor=south west,
    width=0.45\linewidth,
    height=0.45\linewidth,
    xmin=0, xmax=5,
    ymin=0.7, ymax=1.0,
    grid=both,
    clip=true,   
    axis background/.style={fill=white}, 
]
\addplot[color=sczmA, line width=1.3pt]
table[
    x expr=0.5*\thisrow{time},
    y=max_damage
]{sczm_damage_K1e2_GI30_GII200_N10_S20.txt};

\addplot[
    color=ifmA,
    line width=1.3pt,
    dash pattern=on 7pt off 3pt
]
table[
    x expr=0.5*\thisrow{time},
    y=max_damage
]{ifm_damage_K1e2_GI30_GII200_N10_S20.txt};

\addplot[color=sczmB, line width=1.3pt]
table[
    x expr=0.5*\thisrow{time},
    y=max_damage
]{sczm_damage_K50_GI30_GII25_N3_S5.txt};

\addplot[
    color=ifmB,
    line width=1.3pt,
    dash pattern=on 10pt off 3pt
]
table[
    x expr=0.5*\thisrow{time},
    y=max_damage
]{ifm_damage_K50_GI30_GII25_N3_S5.txt};

\end{axis}
\end{tikzpicture}
\caption{Maximum damage}
\end{subfigure}
\vspace{0.3cm}

\begin{tikzpicture}
\begin{axis}[
    hide axis,
    xmin=0, xmax=1,
    ymin=0, ymax=1,
    legend columns=1, 
    legend style={
        draw=none,
        /tikz/every even column/.append style={column sep=0.5cm},
        at={(0.5,0.5)},
        anchor=center
    }
]

\addlegendimage{color=sczmA, line width=1.5pt}
\addlegendentry{$\mathrm{SCZM},\ K=10^{2},\ G_{I c}=30,\ G_{II c}=200,\ N=10,\ S=20$}

\addlegendimage{color=ifmA, line width=1.5pt, dash pattern=on 7pt off 3pt}
\addlegendentry{$\mathrm{IFM},\ K=10^{2},\ G_{I c}=30,\ G_{II c}=200,\ N=10,\ S=20$}

\addlegendimage{color=sczmB, line width=1.5pt}
\addlegendentry{$\mathrm{SCZM},\ K=50,\ G_{I c}=30,\ G_{II c}=25,\ N=3,\ S=5$}

\addlegendimage{color=ifmB, line width=1.5pt, dash pattern=on 10pt off 3pt}
\addlegendentry{$\mathrm{IFM},\ K=50,\ G_{I c}=30,\ G_{II c}=25,\ N=3,\ S=5$}

\end{axis}
\end{tikzpicture}

\caption{Reaction force in the loading direction ($R_x$) and maximum damage as functions of the imposed displacement.}
\label{fig:3D_grain_t_Rx_neml2}

\end{figure}

\figref{fig:sczm_comparison_3D_two_cases} shows displacement magnitude fields for two sets of traction-separation parameters. In both cases, SCZM reproduces the deformation patterns of the interface-fitted solution with high fidelity, demonstrating that the method captures global mechanical response despite the use of non-interface-fitted meshes. Quantitative comparisons of reaction force and maximum damage (\figref{fig:3D_grain_t_Rx_neml2}) show excellent agreement between SCZM and IFM throughout the loading history. This confirms that SCZM accurately captures both the global response and the evolution of interfacial damage in the presence of nonlinear, history-dependent constitutive behavior.

\begin{figure}[htbp]
\centering
\begin{subfigure}[b]{0.49\textwidth}
\centering
\includegraphics[width=\linewidth,trim = {20 0 60 150}, clip]{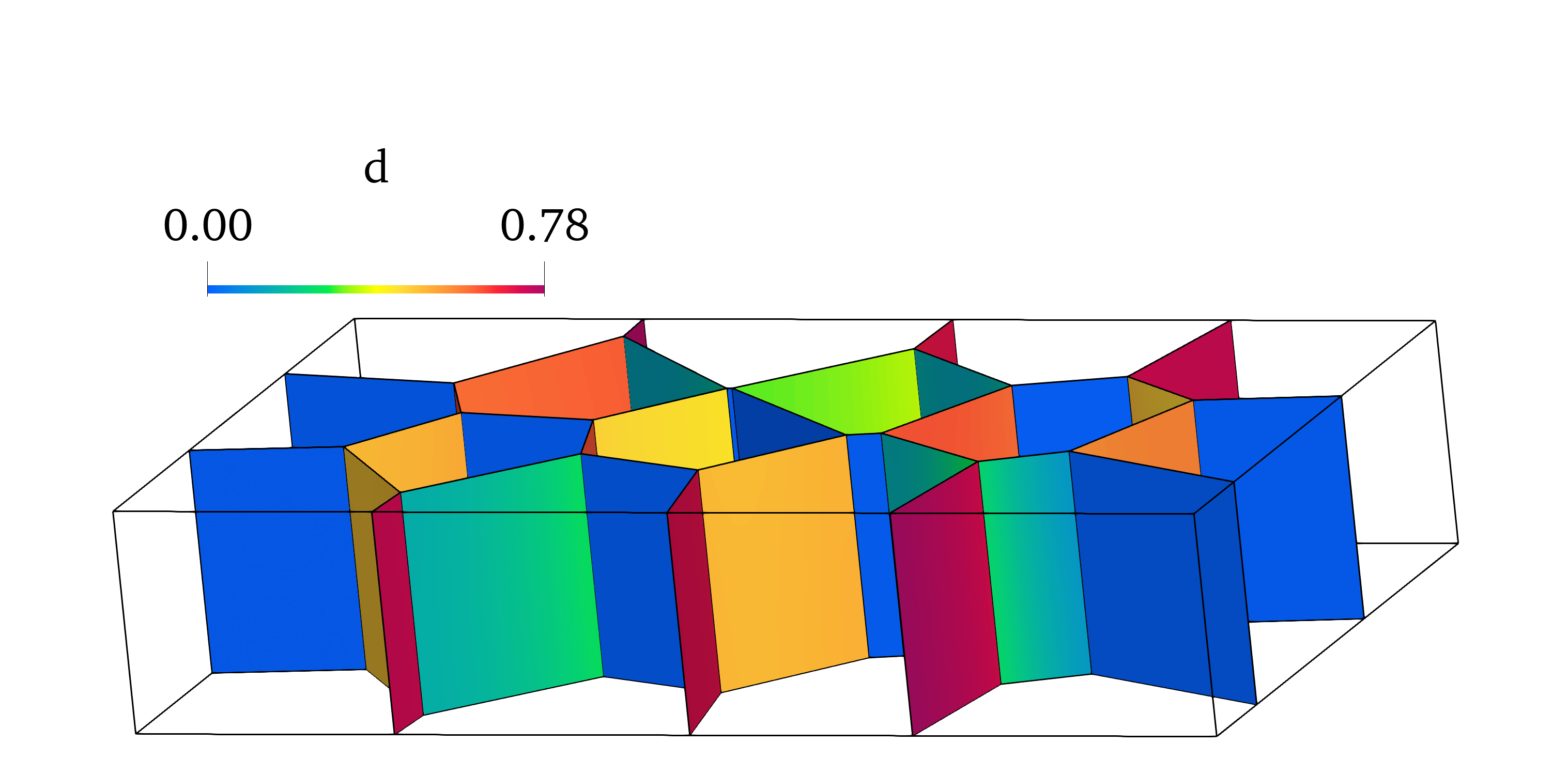}
\caption{IFM at $t=1.75$}
\label{fig:d_t1p75_bfm_surface_K1e2}
\end{subfigure}
\hfill
\begin{subfigure}[b]{0.49\textwidth}
\centering
\includegraphics[width=\linewidth,trim = {20 0 60 150}, clip]{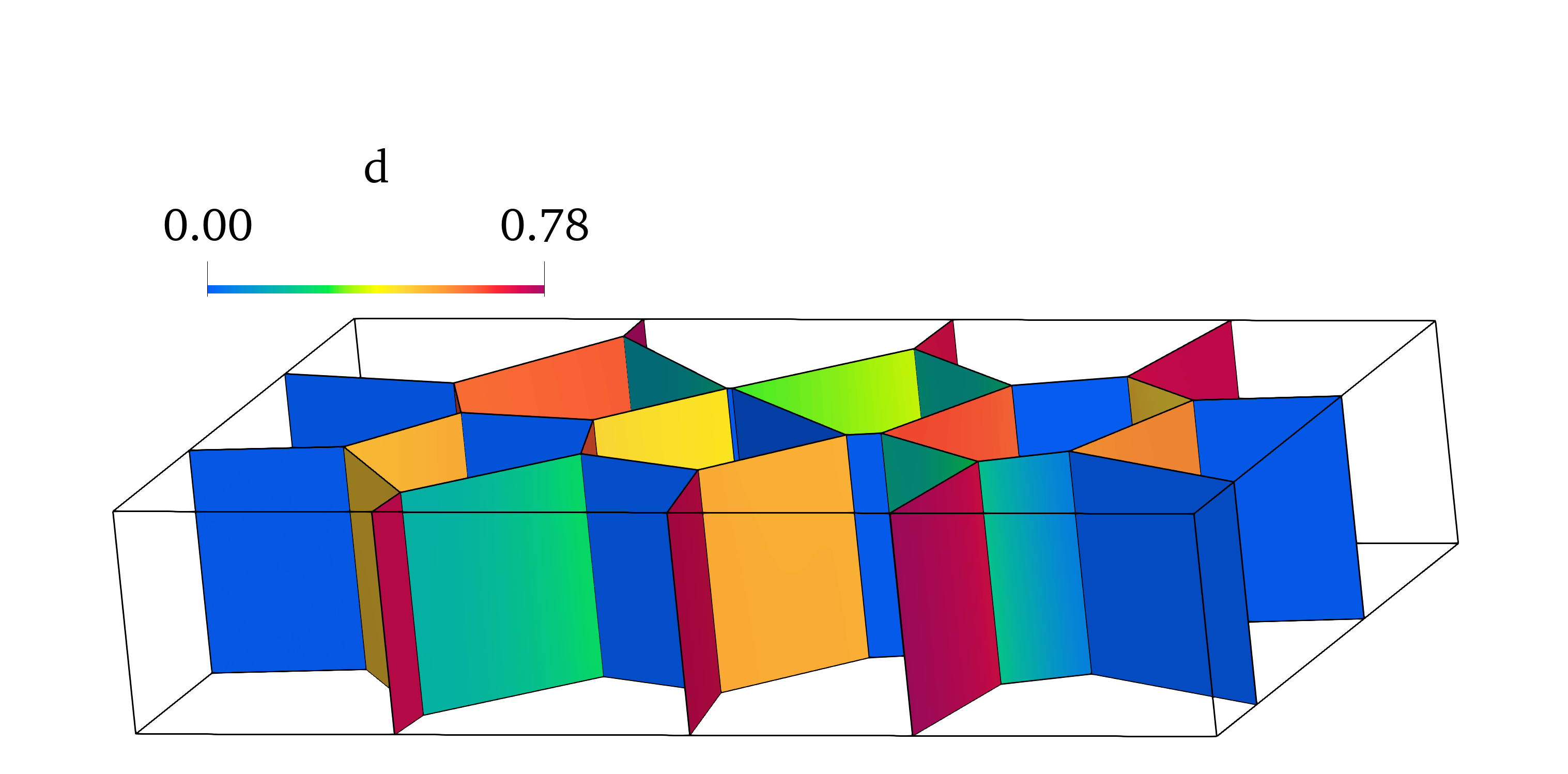}
\caption{SCZM at $t=1.75$}
\label{fig:d_t1p75_sczm_surface_K1e2}
\end{subfigure}

\begin{subfigure}[b]{0.49\textwidth}
\centering
\includegraphics[width=\linewidth,trim = {20 0 60 150}, clip]{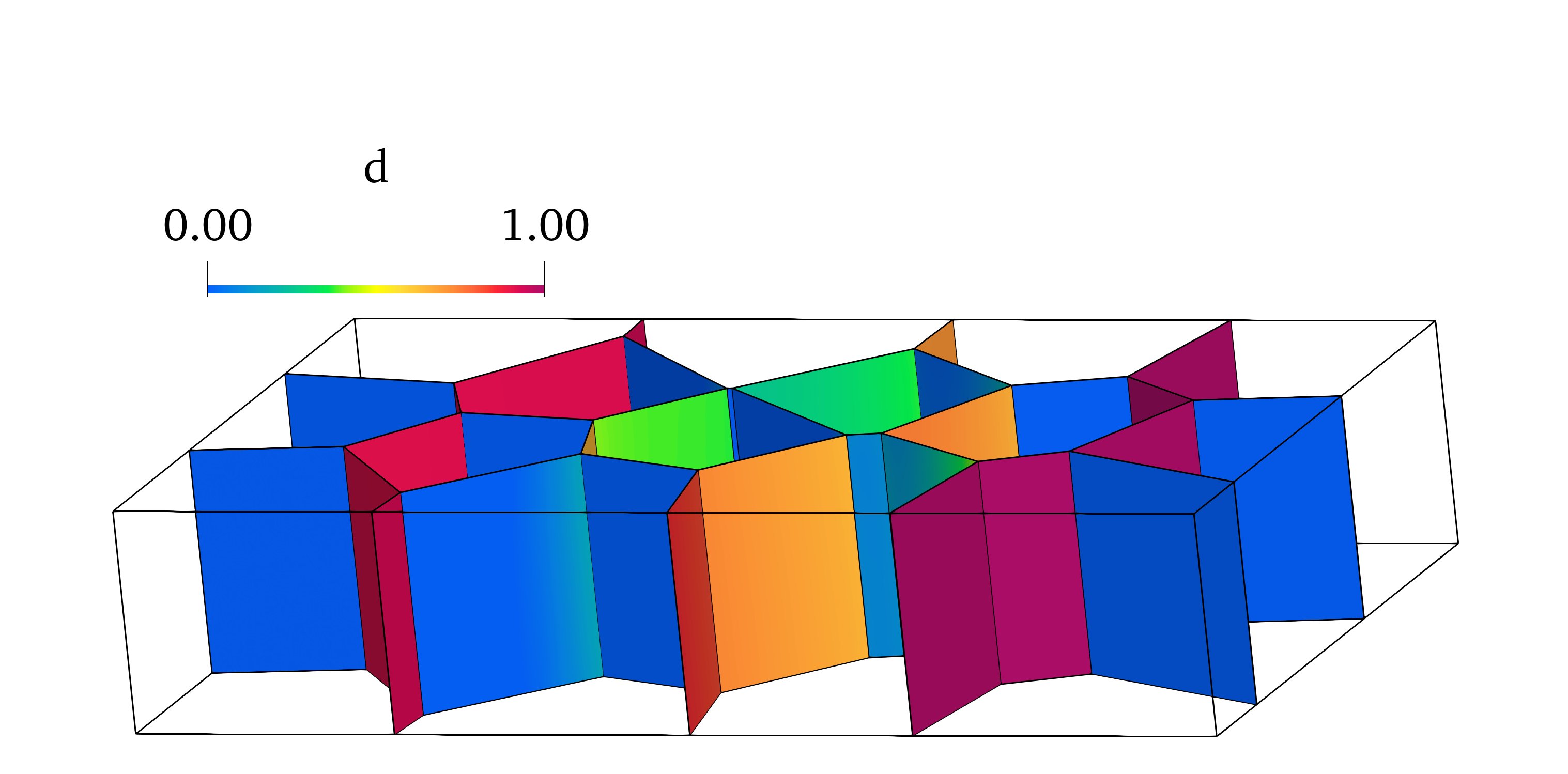}
\caption{IFM at $t=60$}
\label{fig:d_t60_bfm_surface_K1e2}
\end{subfigure}
\hfill
\begin{subfigure}[b]{0.49\textwidth}
\centering
\includegraphics[width=\linewidth,trim = {20 0 60 150}, clip]{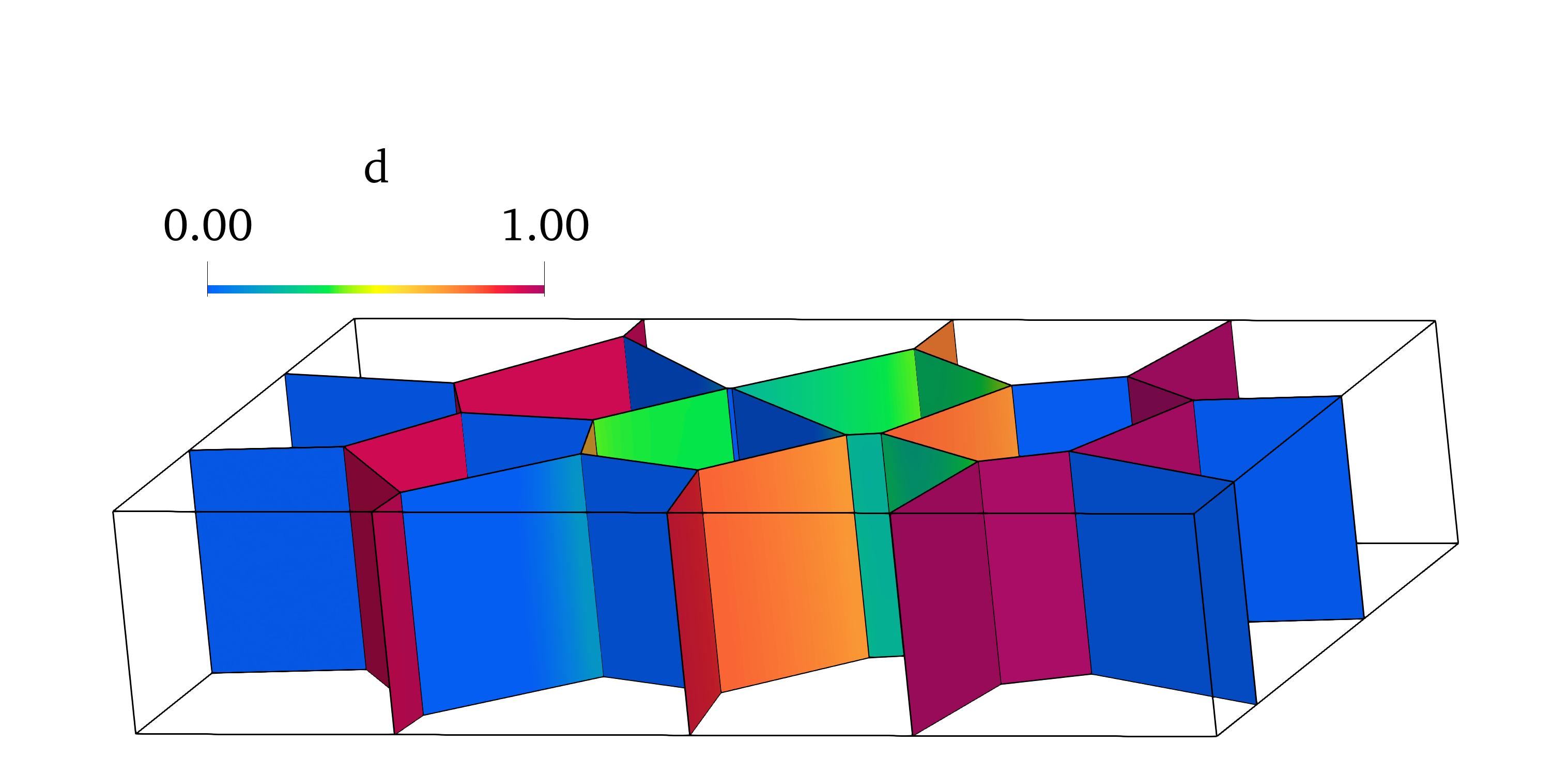}
\caption{SCZM at $t=60$}
\label{fig:d_t60_sczm_surface_K1e2}
\end{subfigure}
\caption{Comparison of the projected damage field on the interface surfaces at $t=1.75$ (top row) and $t=60$ (bottom row) for the IFM and SCZM cases with $(K, G_{I c}, G_{II c}, N, S) = (100, 30, 200, 10, 20)$. The two approaches yield highly consistent damage distributions across the interface network, demonstrating that the projected damage field obtained from SCZM provides a visualization in close agreement with the IFM results.}
\label{fig:d_t1p75_bfm_sczm_compare_K1e2}
\end{figure}

\begin{figure}[htbp]
\centering
\begin{subfigure}[b]{0.49\textwidth}
\centering
\includegraphics[width=\linewidth,trim = {20 0 60 150}, clip]{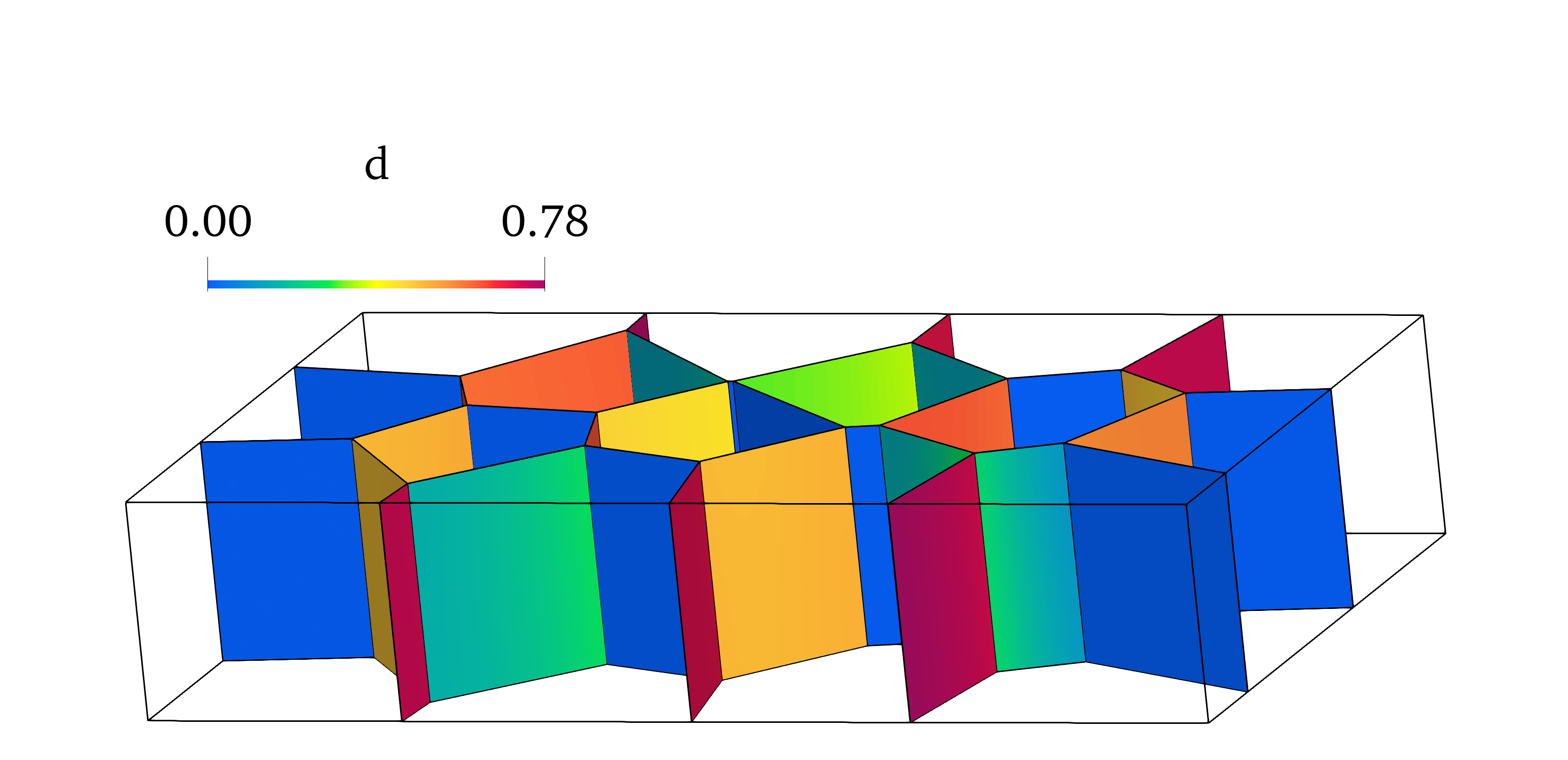}
\caption{IFM at $t=1.75$}
\label{fig:d_t1p75_bfm_surface_K5e1}
\end{subfigure}
\hfill
\begin{subfigure}[b]{0.49\textwidth}
\centering
\includegraphics[width=\linewidth,trim = {20 0 60 150}, clip]{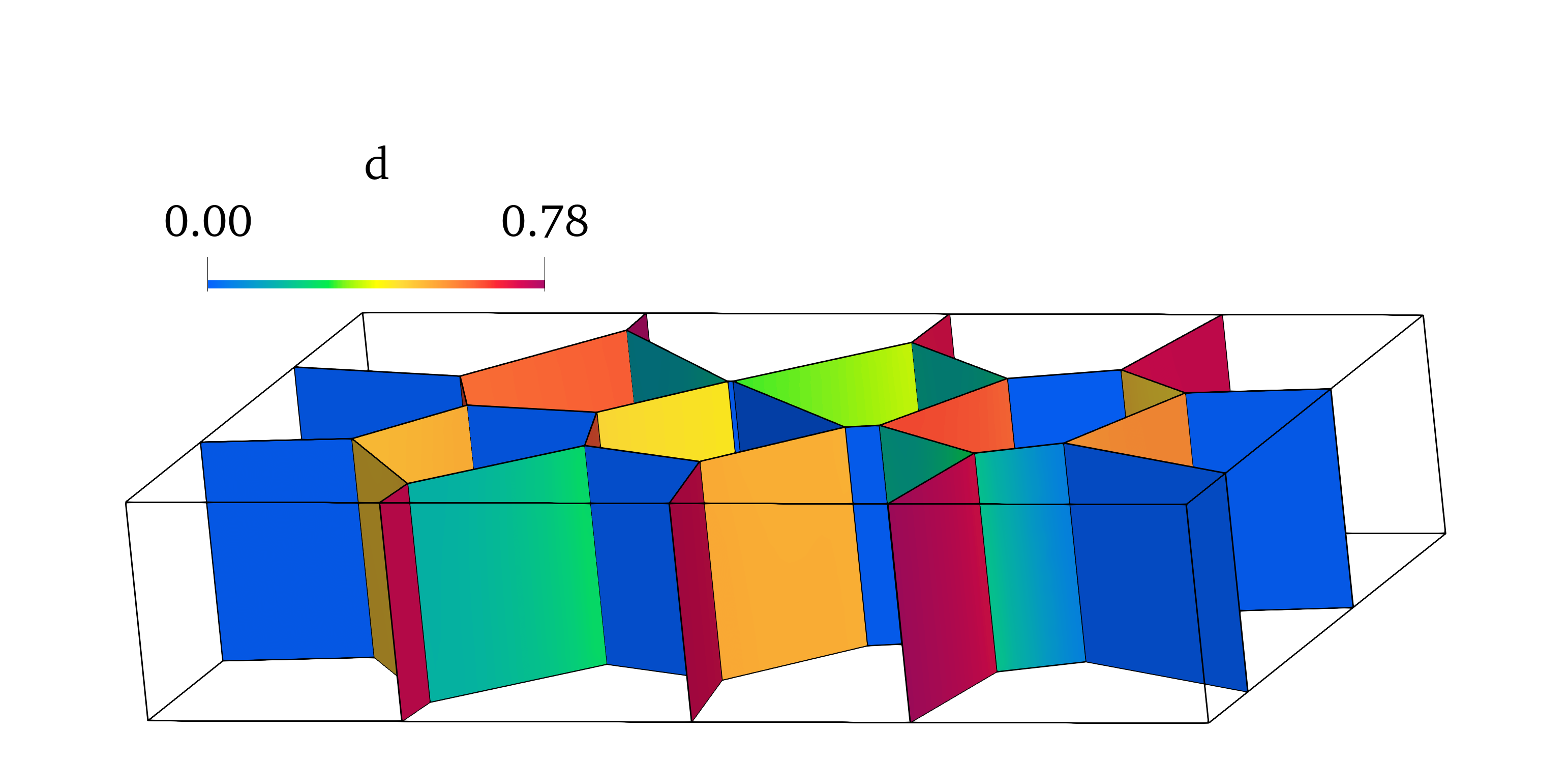}
\caption{SCZM at $t=1.75$}
\label{fig:d_t1p75_sczm_surface_K5e1}
\end{subfigure}

\begin{subfigure}[b]{0.49\textwidth}
\centering
\includegraphics[width=\linewidth,trim = {20 0 60 150}, clip]{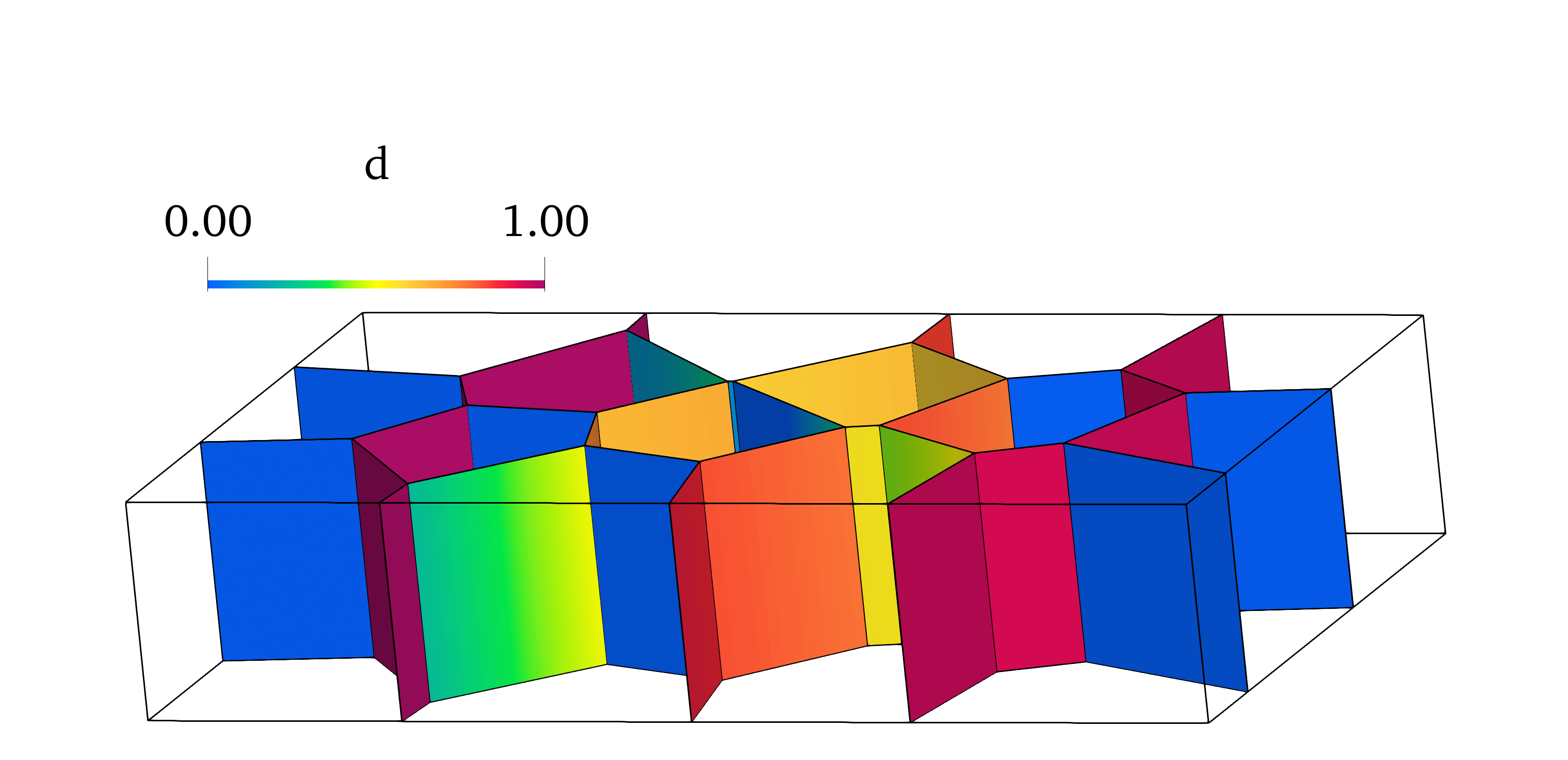}
\caption{IFM at $t=60$}
\label{fig:d_t60_bfm_surface_K5e1}
\end{subfigure}
\hfill
\begin{subfigure}[b]{0.49\textwidth}
\centering
\includegraphics[width=\linewidth,trim = {20 0 60 150}, clip]{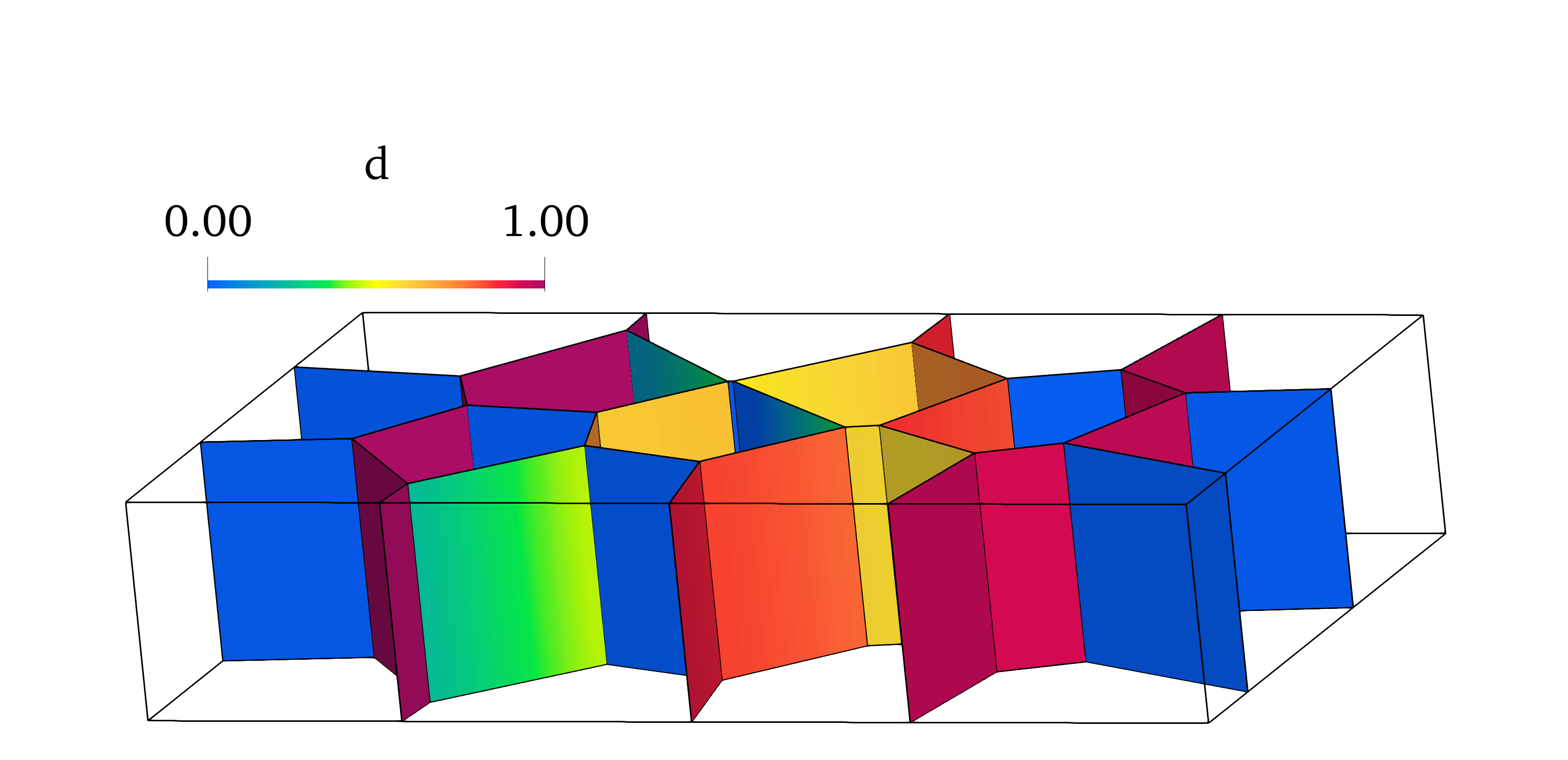}
\caption{SCZM at $t=60$}
\label{fig:d_t60_sczm_surface_K5e1}
\end{subfigure}

\caption{Comparison of the projected damage field on the interface surfaces at $t=1.75$ (top row) and $t=60$ (bottom row) for the IFM and SCZM cases with $(K, G_{I c}, G_{II c}, N, S) = (50, 30, 25, 3, 5)$. The two approaches produce highly consistent damage distributions across the interface network, demonstrating that the projected damage field obtained from SCZM provides a visualization in close agreement with the IFM results even under a softer interface response.}
\label{fig:d_t1p75_bfm_sczm_compare_K5e1}
\end{figure}

To assess interface-level accuracy, damage distributions are compared by projecting both SCZM and IFM results onto the true interface using a consistent projection procedure.
Specifically, for the IFM case, damage values evaluated at interface Gauss points are directly used, while for the SCZM case, surrogate Gauss points are first projected onto the true interface and associated with their corresponding field values.
The resulting pointwise data are then transferred onto an interface mesh via nearest-neighbor or inverse-distance weighting to obtain a continuous representation of the damage field.
The resulting damage fields (\figref{fig:d_t1p75_bfm_sczm_compare_K1e2}, \figref{fig:d_t1p75_bfm_sczm_compare_K5e1}) show strong agreement with the IFM solution for both parameter sets. Minor discrepancies are observed in localized regions, which can be attributed to the omission of higher-order correction terms in the shifted formulation. Specifically, the current SCZM formulation neglects higher-order contributions arising from the geometric gap between the surrogate and true interfaces, which involve derivatives of the traction field. These terms are difficult to recover consistently with first-order finite elements. Recent developments in Gap-SBM~\cite{collins2025gap,antonelli2026isogeometric} suggest a pathway to incorporate such corrections without explicit higher-order reconstruction, representing a promising direction for future work.

Overall, these results demonstrate that SCZM extends naturally to three-dimensional polycrystalline systems, maintaining accuracy, robustness, and scalability without requiring interface-fitted meshes.

\section{Conclusions and Future Work}
\label{Sec:Conclusions}

The shifted cohesive zone method (SCZM) enables accurate enforcement of cohesive interface behavior on non-interface-fitted meshes by evaluating the interface contribution on a surrogate interface. In this work, we extend the shifted fracture framework to cohesive-zone formulations and demonstrate its applicability to problems involving crystal plasticity. The key developments include a simplified weak formulation defined on the surrogate interface, an implementation within the \moose{} framework integrated with \neml{}, and a geometry-aware, PCA-enhanced point classification algorithm that reduces preprocessing cost.

The proposed formulation is verified through MMS-based studies, which confirm first-order convergence consistent with the geometric approximation of the interface, while standard non-interface-fitted approaches fail to converge. Across a range of benchmark problems, including isotropic elasticity, multiple traction-separation laws, and history-dependent crystal plasticity, SCZM reproduces interface-fitted results with high fidelity. The results also highlight the importance of the directional correction term, which is essential for accurately resolving local stress fields and avoiding oscillatory artifacts near misaligned interfaces. Applications to 2D and 3D polycrystalline RVEs demonstrate that SCZM captures interfacial responses in complex microstructures without the need for interface-fitted discretizations.

Several directions for future work remain. Incorporating ideas from the Gap-SBM framework offers a promising path toward higher-order consistency by accounting for the geometric gap between surrogate and true interfaces. Coupling SCZM with adaptive or octree-based discretizations could further improve efficiency for large-scale three-dimensional simulations. Extensions to multiphysics interface problems, such as thermomechanical coupling, and to evolving interfaces would significantly broaden the applicability of the framework.

\section*{Acknowledgements}
The work of the Argonne authors was sponsored by the U.S. Department of Energy, under Contract No. DE-AC02-06CH11357 with Argonne National Laboratory, managed and operated by UChicago Argonne LLC.

Programmatic support was provided by the Nuclear Energy Advanced Modeling \& Simulation Program of the Department of Nuclear Energy. The authors gratefully acknowledge the support and direction of Ben Spencer, at Idaho National Laboratory, program director David Andersson, and Federal Manager Dave Henderson.

\clearpage
\appendix

\section{Supplementary implementation details}
\label{app:supplementary_implementation}

\subsection{Detailed algorithms for morphology-aware point classification}
\label{app:inout_algorithms}

This subsection presents the complete pseudo-code underlying the morphology-aware point-in-polygon framework summarized in \secref{sec:implementation}. While the main text emphasizes the methodological ideas and principal contributions of the approach, the present appendix provides the full algorithmic workflow, including PCA construction, OBB and projected-diagonal evaluation, k-d tree-based candidate filtering, and the multi-direction ray-tracing consistency checks used for robust point classification.

We begin by invoking \Algref{alg:ComputeBoundaryPCA} to compute the principal axes, including the directions of maximum and minimum variance. Subsequently, \Algref{alg:BuildOBBAndKDTree} constructs an oriented bounding box (OBB), the maximum projected diagonal, and a k-d tree~\cite{blanco2014nanoflann} to accelerate the filtering of \SBMElem{}s during ray-element intersection checks.

The main routine for determining the point classification (\textit{sideness}) is described in \Algref{alg:sideness}. This algorithm calls \Algref{alg:TraceRayAgainstGeometry}, which counts the number of intersections between a ray and the \SBMElem{}s. If the number of intersections is even, the point is classified as $\texttt{OUT}$ (outside the geometry); if the number is odd, the point is classified as $\texttt{IN}$ (inside the geometry).

To enhance the robustness of our framework, we shoot an additional ray in the opposite direction -- still aligned with the minimum variance direction -- to verify whether the $\texttt{IN}$ and $\texttt{OUT}$ results from the first two rays are consistent. If one of these rays does not intersect the geometry, the point is immediately classified as $\texttt{OUT}$, and no further consistency check is necessary, as shown in \Algref{alg:sideness}.

If the results of the two initial rays are inconsistent (one ray suggests $\texttt{IN}$, while the other suggests $\texttt{OUT}$), we perform an additional verification by shooting pairs of rays in the maximum variance and second variance directions (along opposite orientations). If any of these additional rays does not intersect the geometry, the point is classified as $\texttt{OUT}$.

If all rays consistently intersect the geometry, we return a warning indicating inconsistent ray results. However, such cases are rare, even when tested on complex 3D geometries.

\begin{algorithm}[htb]
\footnotesize
\caption{\textsc{ComputeBoundaryPCA}}
\begin{algorithmic}[1]
\Require $\{ \mathbf{x}_i \}_{i=1}^N$: coordinates of all boundary nodes
\Ensure Principal axes: $\vec{u}, \vec{v}, \vec{w}$

\State $\bar{\mathbf{x}} \gets \frac{1}{N} \sum_{i=1}^N \mathbf{x}_i$
\State $X \gets \left[ \mathbf{x}_1 - \bar{\mathbf{x}}, \dots, \mathbf{x}_N - \bar{\mathbf{x}} \right]^T$
\State $X = U \Sigma V^T$ \Comment{Singular value decomposition}
\State $\vec{u} \gets V^T_{0,:},\quad \vec{v} \gets V^T_{1,:},\quad \vec{w} \gets V^T_{2,:}$
\State Normalize $\vec{u}, \vec{v}, \vec{w}$
\end{algorithmic}
\label{alg:ComputeBoundaryPCA}
\end{algorithm}

\begin{algorithm}[htb]
\footnotesize
\caption{\textsc{BuildOBBAndKDTreeAndMaxProjectedDiagonal}}
\begin{algorithmic}[1]
\Require PCA basis vectors $\vec{u}, \vec{v}, \vec{w}$, \SBMElem{} boundary elements $\{ e_i \}$
\Ensure Centroids of \SBMElem{} boundary elements $\{ \mathbf{c}_i \}$, the oriented bounding box (OBB), a k-d tree built from projected centroids, and the maximum projected length $L_{\max,\text{projected}}$

\State Initialize bounding intervals: $[u_{\min}, u_{\max}], [v_{\min}, v_{\max}], [w_{\min}, w_{\max}]$
\For{$e_i$ in boundary elements}
  \State $\mathbf{c}_i \gets$ centroid of $e_i$
  \State Project $\mathbf{c}_i$ onto the plane orthogonal to $\vec{r}$ (where $\vec{r} = \vec{v}$ in 2D, $\vec{w}$ in 3D): $\mathbf{p}_i = \mathbf{c}_i - (\mathbf{c}_i \cdot \vec{r}) \vec{r}$
  \For{each node $\mathbf{x}$ in $e_i$}
    \State $u \gets (\mathbf{x} - \bar{\mathbf{x}}) \cdot \vec{u}$, similarly for $v$ and $w$ \Comment{$\bar{\mathbf{x}}$ is obtained from~\Algref{alg:ComputeBoundaryPCA}}
    \State Update extrema for $[u_{\min}, u_{\max}], [v_{\min}, v_{\max}], [w_{\min}, w_{\max}]$
  \EndFor
  \State box $\gets$ AABB (axis-aligned bounding box) of $e_i$
  \State $\vec{\text{diagonal}} \gets \vec{\max(\text{box})} - \vec{\min(\text{box})}$
  \State $\vec{\text{projected\_diagonal}} \gets \vec{\text{diagonal}} - (\vec{\text{diagonal}} \cdot \vec{r}) \vec{r}$
  \If {$L_{\max,\text{projected}} < |\vec{\text{projected\_diagonal}}|$}
    \State $L_{\max,\text{projected}} \gets |\vec{\text{projected\_diagonal}}|$
  \EndIf
\EndFor

\State Define the OBB as: $\bar{\mathbf{x}} + [u_{\min}, u_{\max}]\vec{u} + [v_{\min}, v_{\max}]\vec{v} + [w_{\min}, w_{\max}]\vec{w}$
\State Optionally expand the OBB to account for numerical tolerance
\State Build a k-d tree using the projected centroids $\{\mathbf{p}_i\}$
\end{algorithmic}
\label{alg:BuildOBBAndKDTree}
\end{algorithm}

\begin{algorithm}[htb]
\footnotesize
\caption{\textsc{Point Classification using Multiple Ray Directions}}
\begin{algorithmic}[1]
\Require PCA basis vectors $\vec{u}, \vec{v}, \vec{w}$; query point $\mathbf{p}$; global bounding box; k-d tree; and all \SBMElem{}s on the true boundary $\G$
\Ensure Classification of $\mathbf{p}$ as $\texttt{IN}$, $\texttt{OUT}$, or $\texttt{ON}$

\If{$\mathbf{p}$ lies outside the global bounding box}
  \State \Return $\texttt{OUT}$
\EndIf

\Comment{\textit{Trace rays along the $\vec{w}$ (3D) or $\vec{v}$ (2D) direction first}}
\For{$i \in \{0, 1\}$}
  \State $\mathbf{r}_i^{\mathrm{start}} \gets \texttt{generateRayStart}(\mathbf{p}, \texttt{inverted}= (i == 1), \vec{w} \ \mathrm{in\ 3D}\ /\ \vec{v} \ \mathrm{in\ 2D})$ \Comment{See~\Algref{alg:generateRayStart}}
  \State $c_i \gets \textsc{TraceRayAgainstGeometry}(\mathbf{r}_i^{\mathrm{start}}, \mathbf{p}, L_{\max,\mathrm{proj}}, \mathrm{k\text{-}d\ tree}, \{\SBMElem\})$ \Comment{See~\Algref{alg:TraceRayAgainstGeometry}}
  \If{$c_i == 0$}
    \State \Return $\texttt{OUT}$
  \EndIf
\EndFor

\If{$c_0 \bmod 2 == c_1 \bmod 2$}
  \State \Return $\texttt{IN}$ if $c_0$ is odd; otherwise $\texttt{OUT}$
\EndIf

\Comment{\textit{Use additional ray directions ($\vec{u}, \vec{v}$ in 3D; $\vec{u}$ in 2D) for verification}}
\For{$\mathbf{e} \in \begin{cases}
\{ \vec{u}, \vec{v} \}, & \mathrm{in\ 3D} \\
\{ \vec{u} \}, & \mathrm{in\ 2D}
\end{cases}$}
  \For{$i \in \{0, 1\}$}
    \State $\mathbf{r}_i^{\mathrm{start}} \gets \texttt{generateRayStart}(\mathbf{p}, \texttt{inverted}= (i == 1), \mathbf{e})$ \Comment{See~\Algref{alg:generateRayStart}}
    \State $c_i \gets \textsc{TraceRayAgainstGeometry}(\mathbf{r}_i^{\mathrm{start}}, \mathbf{p}, \{\SBMElem\})$ \Comment{See~\Algref{alg:TraceRayAgainstGeometry}}
    \If{$c_i == 0$}
      \State \Return $\texttt{OUT}$
    \EndIf
  \EndFor
\EndFor

\State \Return Warning: inconsistent ray results, possible geometry error.
\end{algorithmic}
\label{alg:sideness}
\end{algorithm}

\begin{algorithm}[htb]
\footnotesize
\caption{\textsc{TraceRayAgainstGeometry}: Count ray-geometry intersections}
\begin{algorithmic}[1]
\Require Ray origin $\mathbf{ray}_{start}$, endpoint $\mathbf{ray}_{end}$, k-d tree, search radius $L_{\max,\text{projected}}$, and all \SBMElem{}s on the true boundary $\G$ (\{\SBMElem\})
\Ensure Classification of $\mathbf{ray}_{end}$ ($\texttt{IN, OUT, ON}$)

\State Initialize intersection counter: $I \gets 0$
\If{k-d tree and $L_{\max,\text{projected}}$ are provided}
  \State candidate \SBMElem{}s $\gets$ \textsc{CollectCandidateElementIDs}($\mathbf{ray}_{start}$, $L_{\max,\text{projected}}$) \Comment{See~\Algref{alg:CollectCandidateElementIDs}}
\Else
  \State candidate \SBMElem{}s $\gets$ all \SBMElem{}s on the geometry
\EndIf
\For{each \SBMElem{} in candidate list}
  \State $\mathbf{c}$, r $\gets \text{Bounding ball of \SBMElem{}}$
  \If{\textsc{RayRegionFastReject}($\mathbf{ray}_{start}, \mathbf{ray}_{end} - \mathbf{ray}_{start}, \mathbf{c}, r$)}  \Comment{See~\Algref{alg:bboxcheck}}
    \State \textbf{continue}
  \EndIf
  \If{point lies on \SBMElem{}}
    \State \Return $\texttt{ON}$
  \EndIf
  \State $T \gets$ \textsc{IfLineIntersects\SBMElem}($\mathbf{ray}_{start}, \mathbf{ray}_{end}, \SBMElem$)
  \If{$T = \texttt{INTERSECT}$}
    \State $I \gets I + 1$
  \EndIf
\EndFor
\If{$I$ is odd}
  \State \Return $\texttt{IN}$
\Else
  \State \Return $\texttt{OUT}$
\EndIf
\end{algorithmic}
\label{alg:TraceRayAgainstGeometry}
\end{algorithm}

\begin{algorithm}[htb]
\footnotesize
\caption{\textsc{CollectCandidateElementIDs}}
\begin{algorithmic}[1]
\Require Query point $\mathbf{q}$, k-d tree, search radius $L_{\max,\text{projected}}$, and all \SBMElem{}s on the true boundary $\G$ (\{\SBMElem\})
\Ensure List of selected elements

\State Project $\mathbf{q}$ onto the plane orthogonal to $\vec{r}$ ($\vec{v}$ in 2D, $\vec{w}$ in 3D), yielding $\mathbf{q}'$
\State Perform k-d tree radius search: $\texttt{radiusSearch}(\mathbf{q}', L_{\max,\text{projected}})$
\State IDs $\gets \text{list of matched element IDs}$
\State \Return \{\SBMElem\}[IDs]
\end{algorithmic}
\label{alg:CollectCandidateElementIDs}
\end{algorithm}

\begin{algorithm}[htb]
\footnotesize
\caption{\textsc{generateRayStart}: Generate the starting point for ray tracing}
\begin{algorithmic}[1]
\Require Query point $\mathbf{p}$, inversion flag $\texttt{inverted}$, ray direction (one of the PCA basis vectors $\{\vec{u}, \vec{v}, \vec{w}\}$)
\Ensure Ray starting point $\mathbf{r}^{\mathrm{start}}$

\State $\texttt{SAFE\_FACTOR} \gets 1.1$
\State $\texttt{axisLength} \gets$ length of the OBB axis aligned with the chosen PCA basis vector
\State $\texttt{halfAxisLength} \gets \texttt{axisLength} / 2.0$

\If{the projection of $\mathbf{p}$ along the chosen basis vector $<$ $\texttt{halfAxisLength}$}
    \State $\mathbf{corner} \gets$ \texttt{getMaximalCorner()} if $\texttt{inverted}$, else \texttt{getMinimalCorner()}
    \State $\texttt{multiplier} \gets -1.0$ if $\texttt{inverted}$, else $1.0$
\Else
    \State $\mathbf{corner} \gets$ \texttt{getMinimalCorner()} if $\texttt{inverted}$, else \texttt{getMaximalCorner()}
    \State $\texttt{multiplier} \gets 1.0$ if $\texttt{inverted}$, else $-1.0$
\EndIf

\State $\mathbf{projPoint} \gets$ projection of $\mathbf{p}$ onto the plane orthogonal to the chosen basis vector and passing through $\mathbf{corner}$
\State \Return $\mathbf{projPoint} - \texttt{SAFE\_FACTOR} \times \texttt{axisLength} \times \texttt{multiplier} \times \mathbf{rayDirection}$

\end{algorithmic}
\label{alg:generateRayStart}
\end{algorithm}

\begin{algorithm}[htb]
\footnotesize
\caption{\textsc{RayRegionFastReject}: Quick exclusion of far elements}
\begin{algorithmic}[1]
\Require Ray origin $\mathbf{o}$, direction $\vec{d}$, target center $\mathbf{c}$, radius $r$
\Ensure \texttt{true} if rejection is valid

\State Compute AABB: $[\mathbf{l}, \mathbf{u}] = \text{min/max}(\mathbf{o}, \mathbf{o} + \vec{d}) \pm r$
\If{$\mathbf{c} \notin [\mathbf{l}, \mathbf{u}]$}
    \State \Return \texttt{true}
\EndIf
\State Project: $\mathbf{P}_b = \mathbf{o} + \frac{(\mathbf{c} - \mathbf{o}) \cdot \vec{d}}{\vec{d} \cdot \vec{d}} \vec{d}$
\If{$\|\mathbf{P}_b - \mathbf{c}\| > r$}
\State \Return \texttt{true}
\Else
\State \Return \texttt{false}
\EndIf
\end{algorithmic}
\label{alg:bboxcheck}
\end{algorithm}

\subsection{Element-wise grain assignment for surrogate RVEs}
\label{app:rve_generation}

To generate surrogate RVEs, we first obtain the water-tight boundary representation of each grain, i.e., $\Sigma = \sum_{i=1}^N \G_i$. The assignment principle is to identify the grain that occupies the largest fraction of an element's volume and assign that grain's ID to the element. This avoids the time-consuming process of generating conformal RVEs and enables efficient crystal plasticity simulations with surrogate RVEs.

\begin{algorithm}[htb]
\footnotesize
\caption{\textsc{ElementGrainIDMarkers}: Assigning grain ID to an element}
\begin{algorithmic}[1]
\Require Element $\mathcal{E}$, set of grain boundaries $\Sigma = \{\G_i\}_{i=1}^N$
\Ensure Grain ID $g^\ast$ associated with $\mathcal{E}$
\State Initialize $\phi^\ast \gets 0$ \Comment{Maximum volume fraction among grains}
\State Initialize $g^\ast \gets \varnothing$ \Comment{Index of the grain with $\phi^\ast$}
\For{each grain boundary $\G_i \in \Sigma$}
    \State Determine how many nodes of $\mathcal{E}$ lie inside $\G_i$
    \If{all nodes of $\mathcal{E}$ are inside $\G_i$}
        \State \Return $i$ \Comment{Element fully belongs to grain $i$}
    \Else
        \State Compute volume fraction $\phi_i = \dfrac{|\mathcal{E} \cap \Omega_i|}{|\mathcal{E}|}$
        \If{$\phi_i > \phi^\ast$}
            \State Update $\phi^\ast \gets \phi_i$
            \State Update $g^\ast \gets i$
        \EndIf
    \EndIf
\EndFor
\State \Return $g^\ast$ \Comment{Assign to grain with maximum occupancy}
\end{algorithmic}
\label{alg:surrogate_RVE}
\end{algorithm}

\section{Algorithmic Details and Visualization of SCZM-to-IFM Transfer}
\label{app:sczm_to_ifm_algorithms}

This section provides detailed algorithmic descriptions and visual validation for \secref{sec:SCZM2IFM}.

The algorithms for \secref{sec:SCZM2IFM} are summarized in \Algref{alg:ConformalizePixelatedMeshWithBoundaryRegions}, \Algref{alg:BuildSCZMSearchStructure}, \Algref{alg:ProjectSolutionToIFM}, and \Algref{alg:FindSourceElem}. The central geometric operation involves computing the intersection between an SCZM element $K$ and a material region $\Omega_N$, followed by a simplex decomposition of the resulting clipped region. In two dimensions, this corresponds to polygon clipping and triangulation, whereas in three dimensions, it corresponds to polyhedral clipping and tetrahedralization.

\begin{figure}[htb]
\centering

    \centering
    \includegraphics[width=\linewidth,trim = {0 0 0 0}, clip]{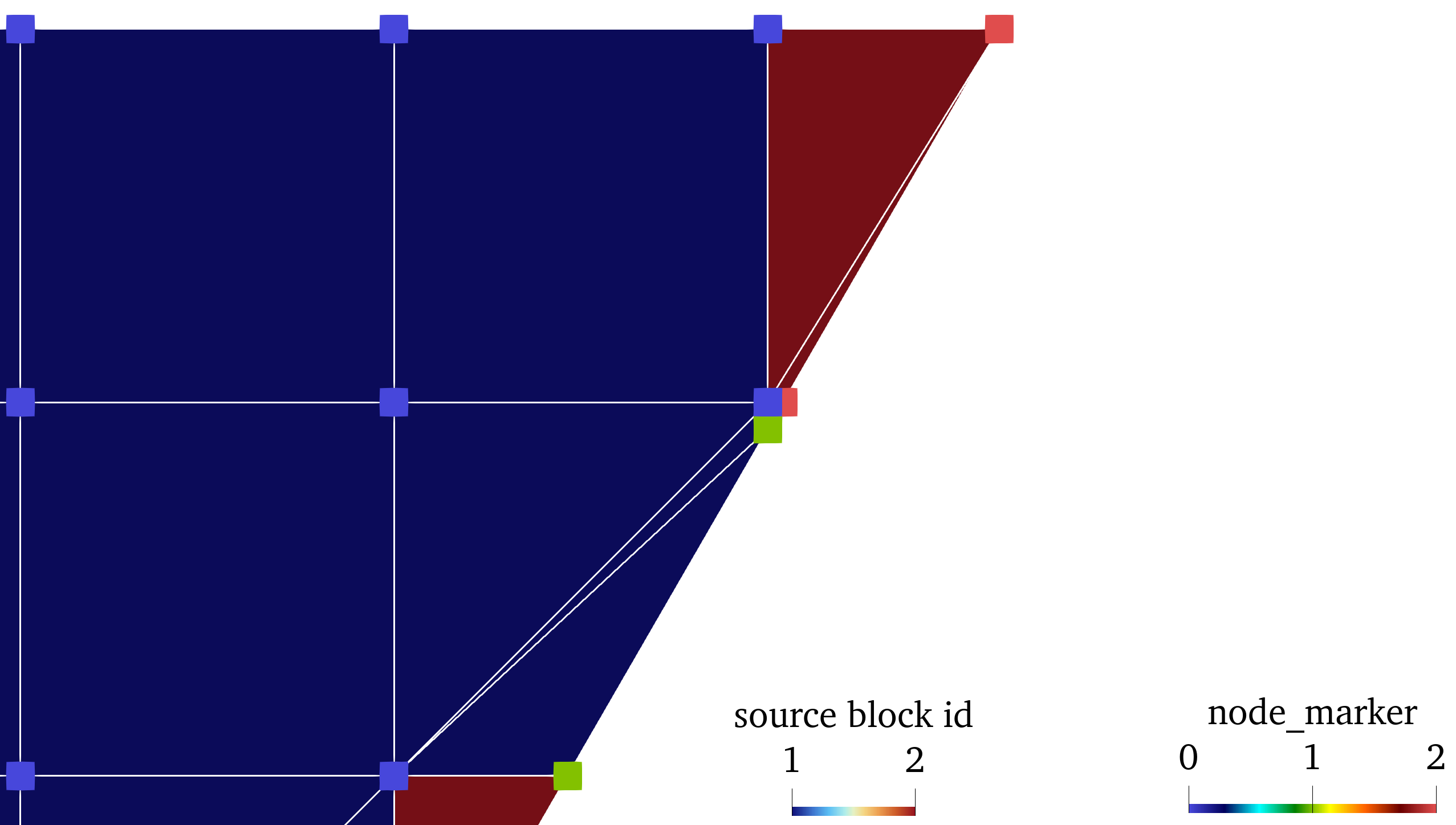}
    \caption{
    Visualization of the \texttt{node\_marker} classification on the generated IFM mesh.
    Rectangle-shaped points denote mesh nodes, while element colors indicate the corresponding source block or subdomain.
    Because the SCZM mesh is non-interface-fitted, the conformalization process locally removes or adds portions of the mesh topology to construct an interface-fitted mesh.
    Nodes directly inherited from the original SCZM mesh are assigned \texttt{node\_marker}$=0$.
    Newly generated nodes whose block ID matches the source element block ID are assigned \texttt{node\_marker}$=1$.
    Newly generated nodes whose block ID differs from the source element block ID are assigned \texttt{node\_marker}$=2$; these nodes arise from locally extended regions introduced during the IFM mesh-generation process.
    }
    \label{fig:node_marker}
\end{figure}

The \texttt{node\_marker} encodes the origin of each output node in \Algref{alg:ConformalizePixelatedMeshWithBoundaryRegions}, where \texttt{node\_marker}$=0$ denotes nodes coincident with original SCZM nodes, \texttt{node\_marker}$=1$ denotes newly generated nodes within the same material block, and \texttt{node\_marker}$=2$ denotes nodes associated with regions outside the original block assignment; these markers determine whether the solution is copied, interpolated, or recovered via a local Taylor expansion in \Algref{alg:ProjectSolutionToIFM}. A visual illustration and detailed explanation of the \texttt{node\_marker} classification are provided in \figref{fig:node_marker}.

\begin{figure}[htb]
\centering

\begin{subfigure}[b]{0.32\textwidth}
    \centering
    \includegraphics[width=\linewidth,trim = {300 10 1200 100}, clip]{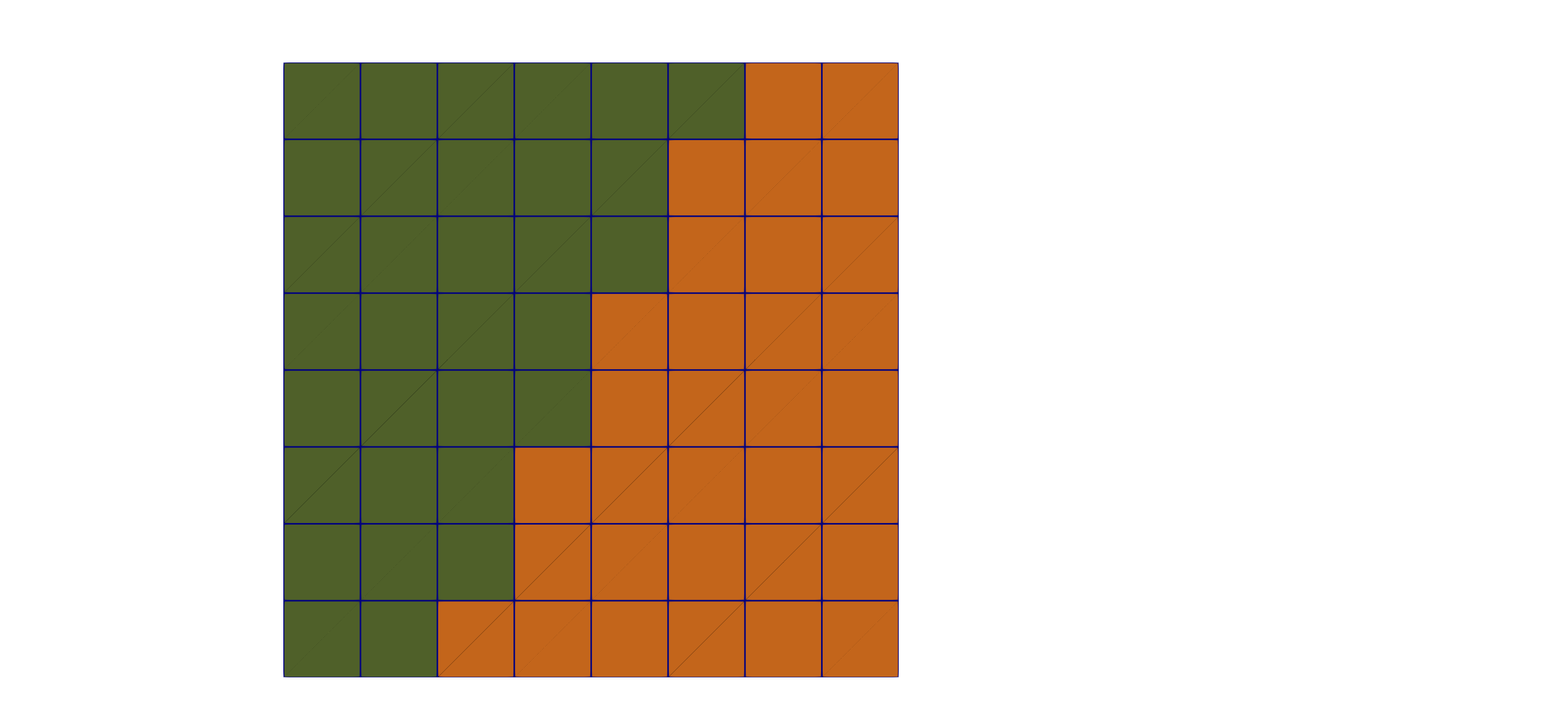}
    \caption{SCZM ($n=8$)}
\end{subfigure}
\hfill
\begin{subfigure}[b]{0.32\textwidth}
    \centering
    \includegraphics[width=\linewidth,trim = {300 10 1200 100}, clip]{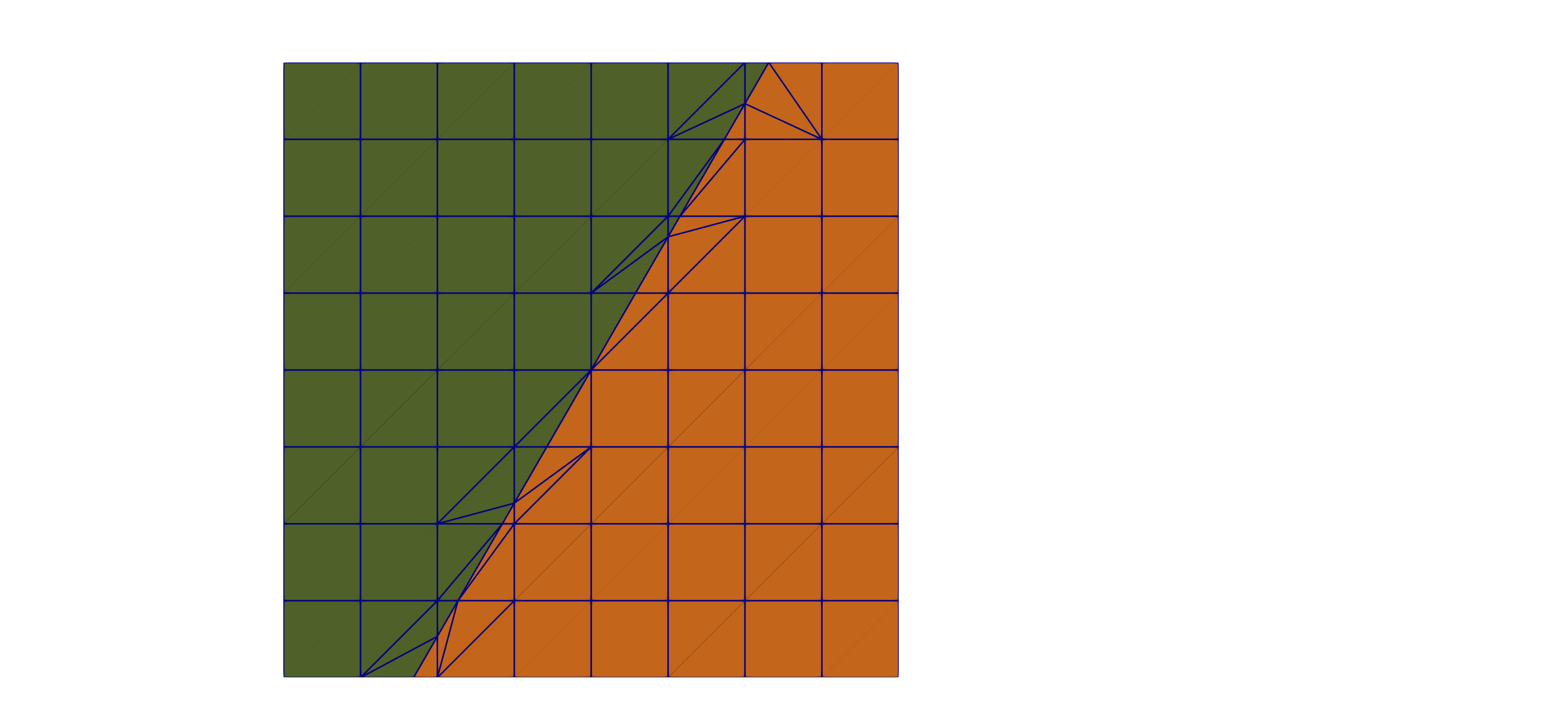}
    \caption{Interface-fitted mesh generated from the SCZM mesh ($n=8$)}
\end{subfigure}

\vspace{0.4em}

\begin{subfigure}[b]{0.32\textwidth}
    \centering
    \includegraphics[width=\linewidth,trim = {300 10 1200 100}, clip]{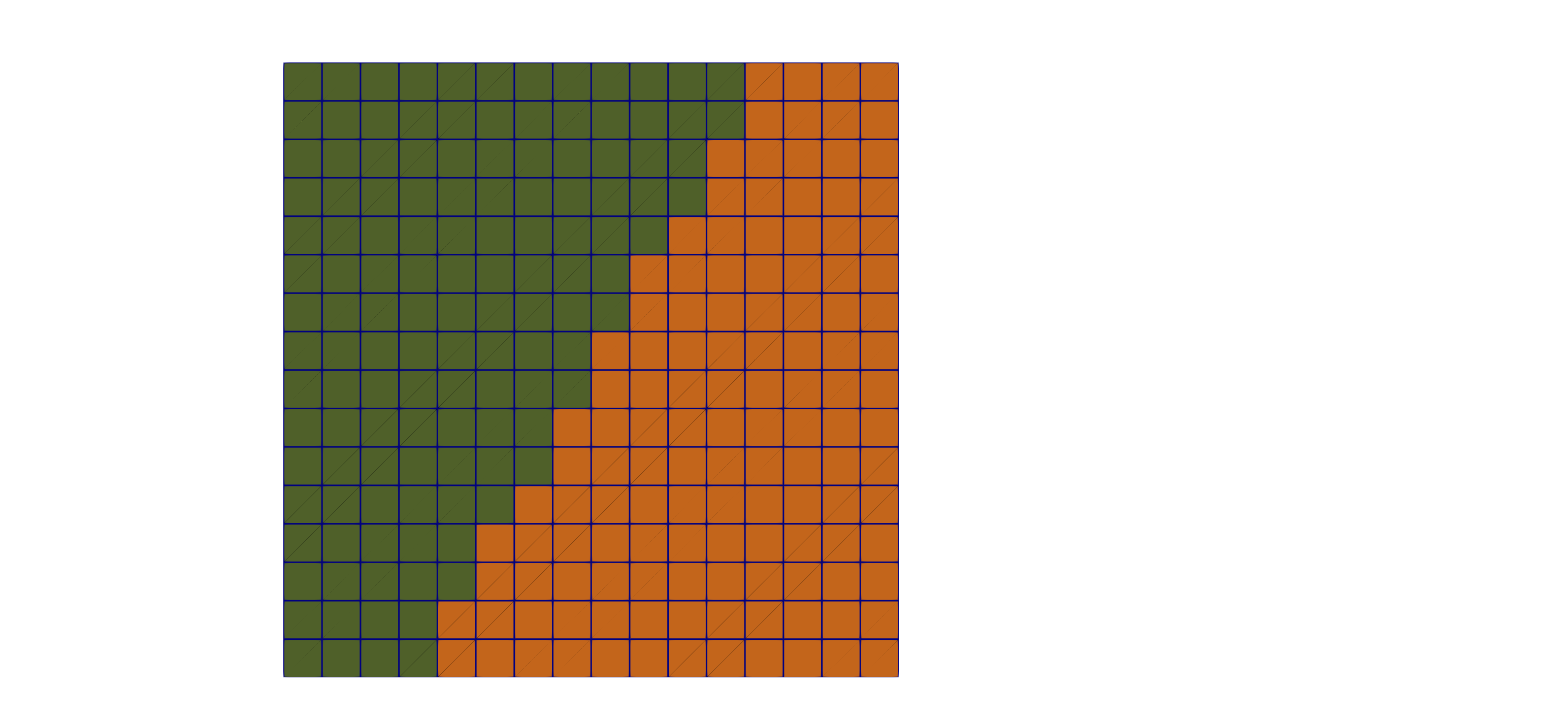}
    \caption{SCZM ($n=16$)}
\end{subfigure}
\hfill
\begin{subfigure}[b]{0.32\textwidth}
    \centering
    \includegraphics[width=\linewidth,trim = {300 10 1200 100}, clip]{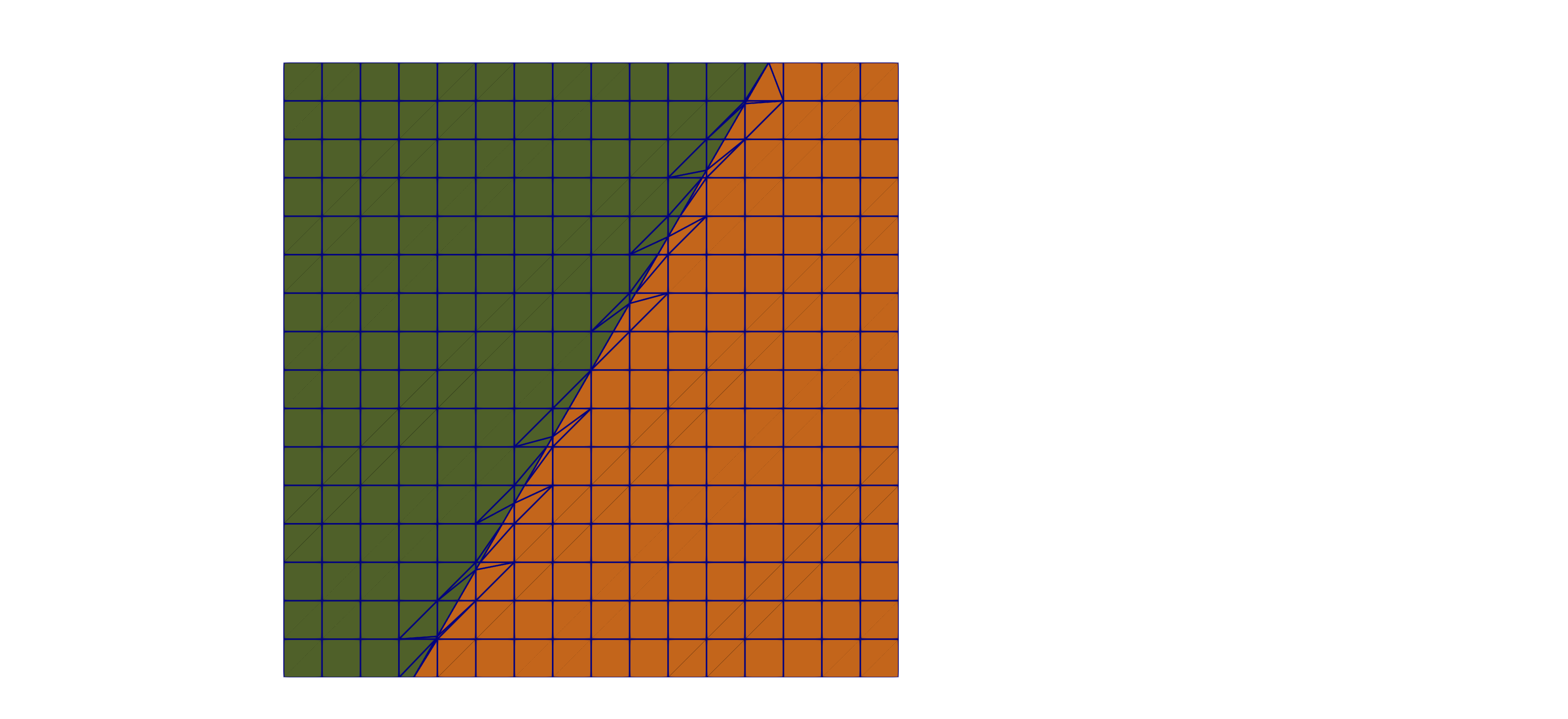}
    \caption{Interface-fitted mesh generated from the SCZM mesh ($n=16$)}
\end{subfigure}

\vspace{0.4em}

\begin{subfigure}[b]{0.32\textwidth}
    \centering
    \includegraphics[width=\linewidth,trim = {300 10 1200 100}, clip]{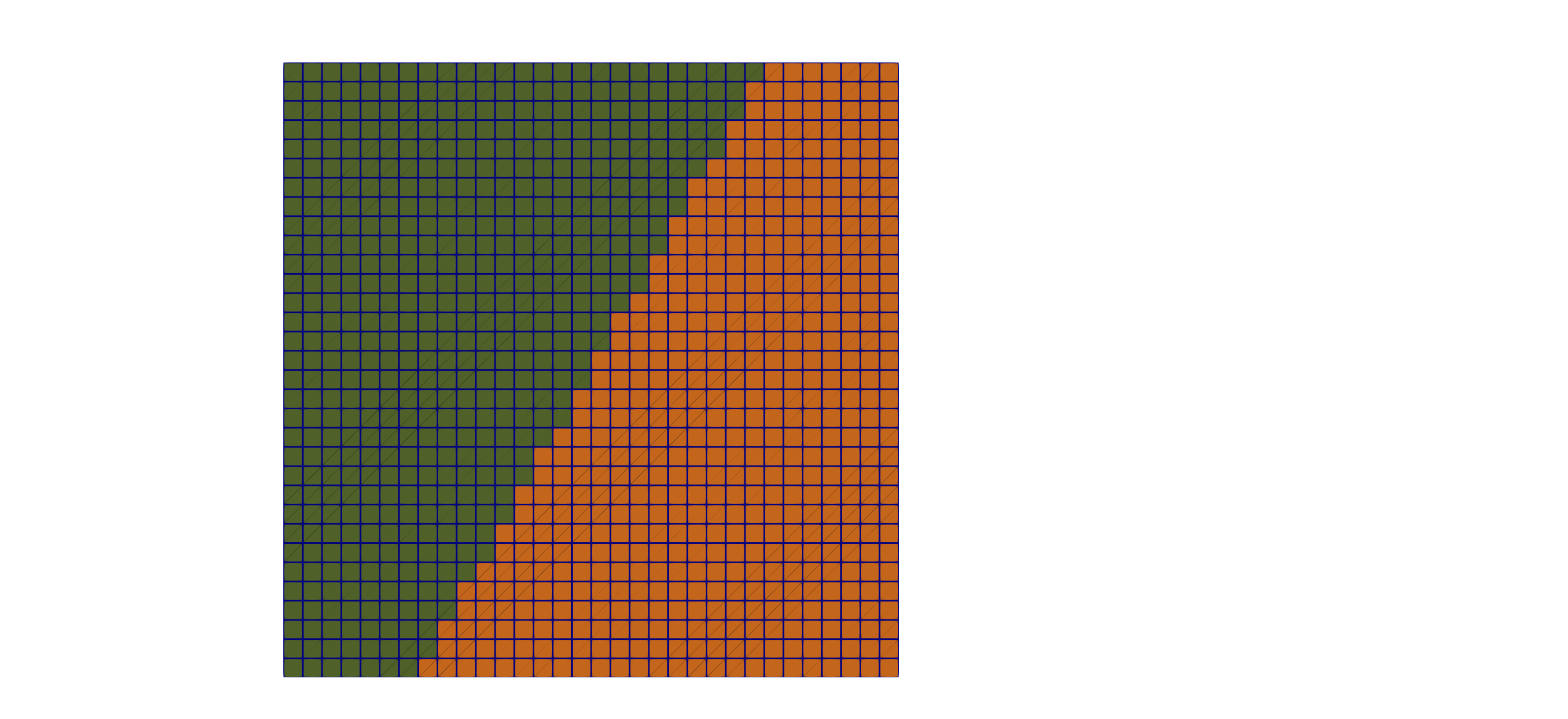}
    \caption{SCZM ($n=32$)}
\end{subfigure}
\hfill
\begin{subfigure}[b]{0.32\textwidth}
    \centering
    \includegraphics[width=\linewidth,trim = {300 10 1200 100}, clip]{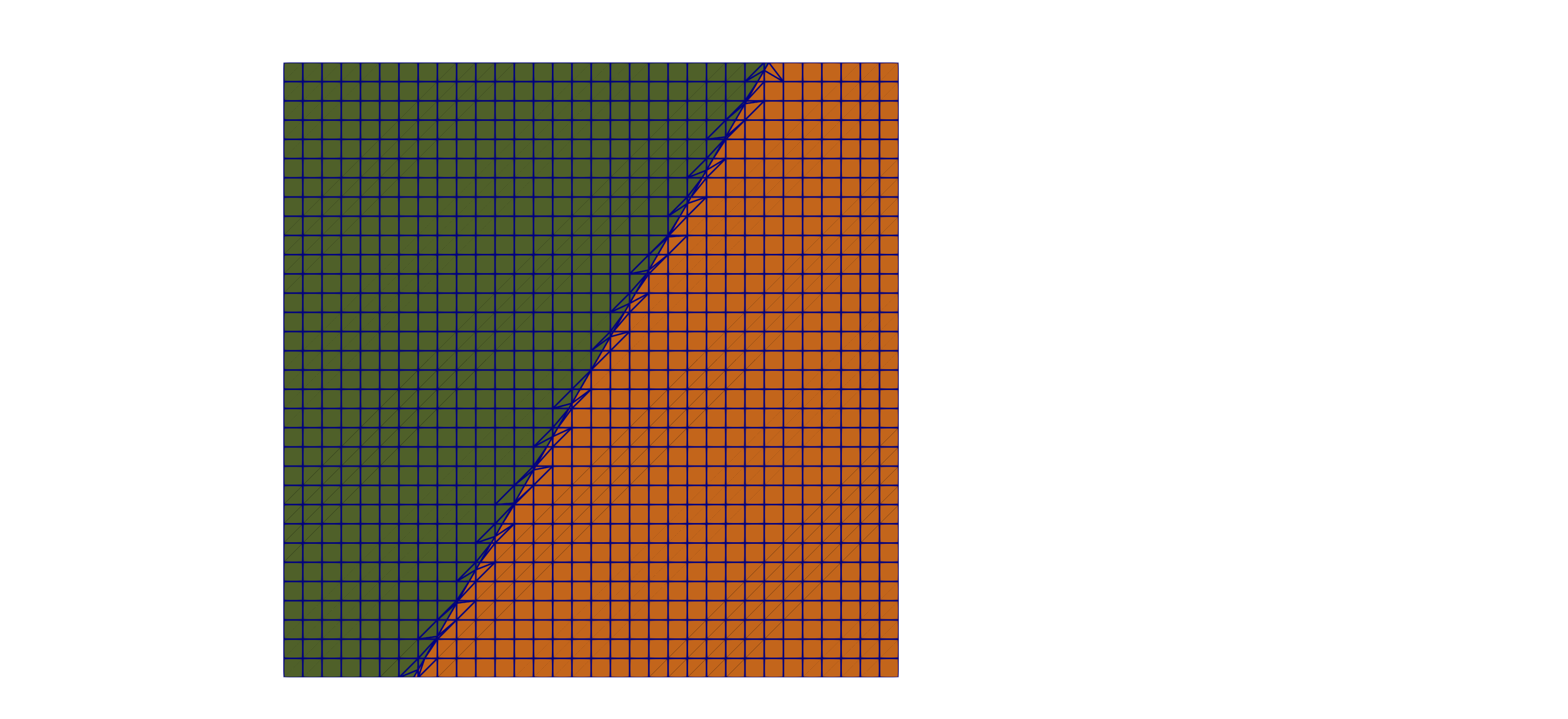}
    \caption{Interface-fitted mesh generated from the SCZM mesh ($n=32$)}
\end{subfigure}

\caption{
Illustration of the conformalization procedure from SCZM to IFM meshes.
The left column shows the original non-interface-fitted SCZM meshes, while the right column shows the corresponding interface-fitted meshes generated by the proposed algorithm.
Only elements intersected by the interface are locally repartitioned, while the remaining elements are preserved.
As the mesh resolution is increased from $n=8$ to $n=32$, the interface representation becomes progressively more accurate.
}
\label{fig:sczm_to_ifm_mesh}
\end{figure}

The resulting interface-fitted meshes are shown in \figref{fig:sczm_to_ifm_mesh}, where only elements intersected by the interface are locally modified while the majority of the SCZM mesh is retained, yielding an accurate interface representation without global remeshing.

\begin{figure}[htb]
\centering

\begin{minipage}{0.78\textwidth}
\centering

\begin{subfigure}[b]{0.48\textwidth}
    \centering
    \includegraphics[width=\linewidth,trim = {500 100 500 100}, clip]
    {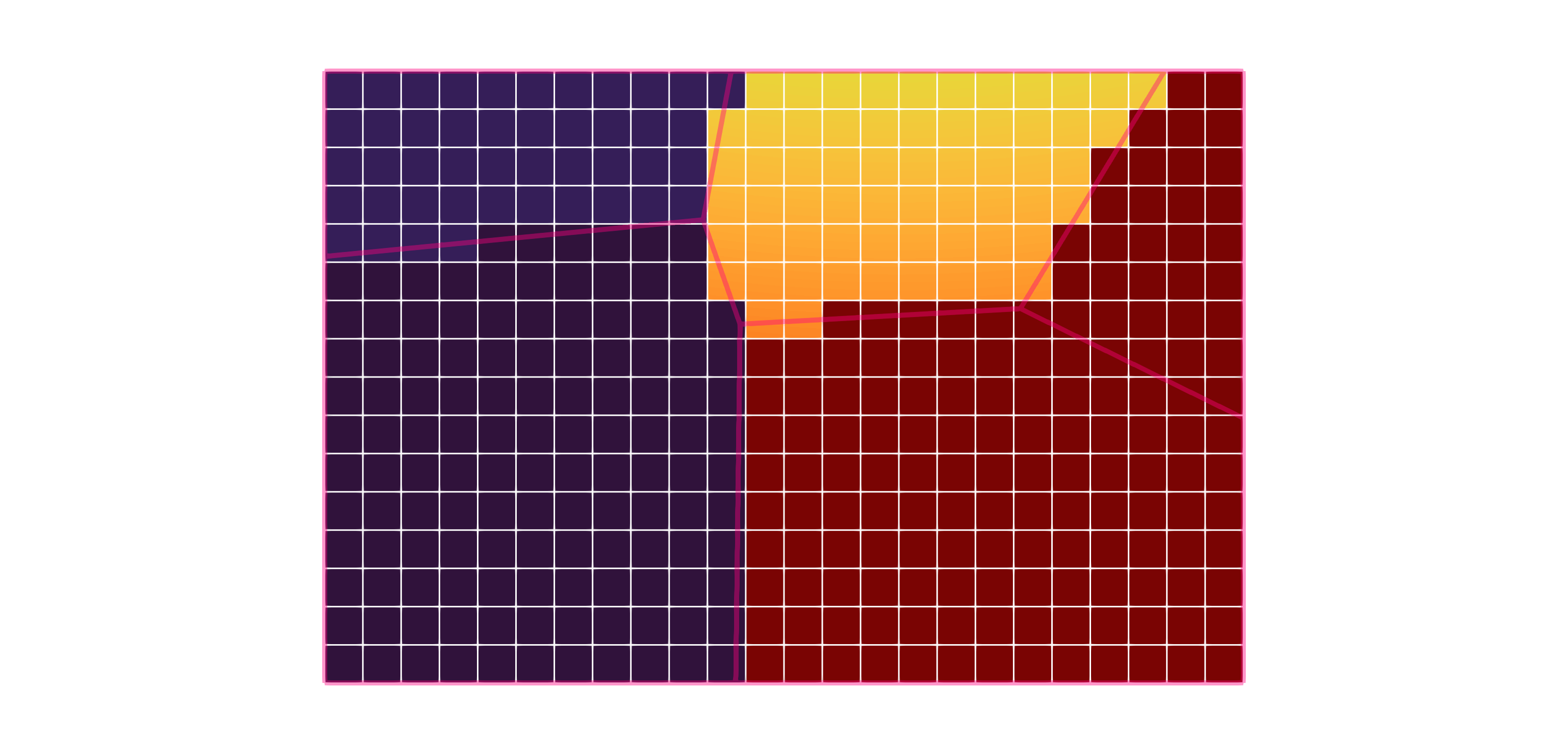}
    \caption{Block-wise partition on the original SCZM mesh.}
    \label{fig:sczm_block}
\end{subfigure}
\hfill
\begin{subfigure}[b]{0.48\textwidth}
    \centering
    \includegraphics[width=\linewidth,trim = {500 100 500 100}, clip]
    {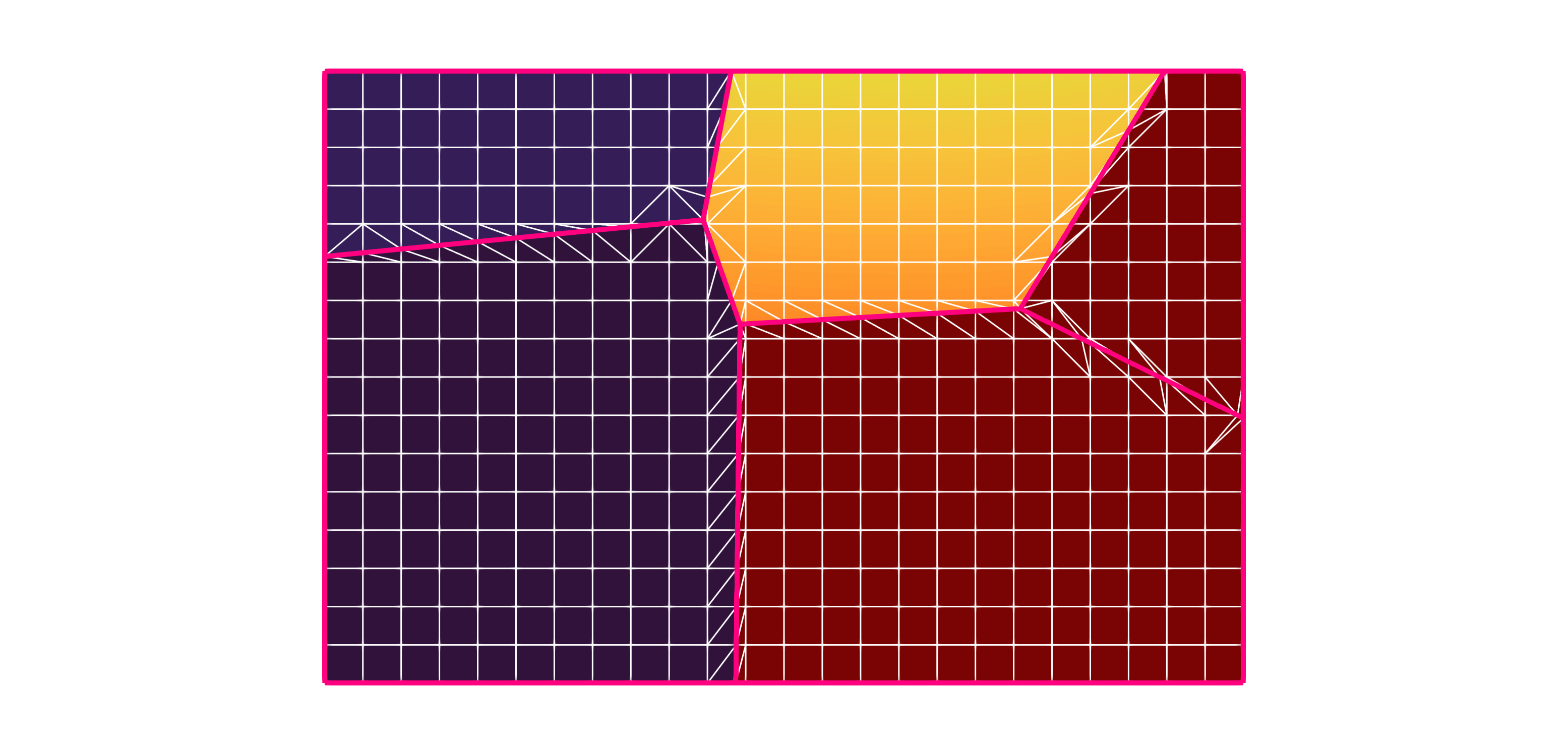}
    \caption{Reconstructed partition on the generated IFM mesh.}
    \label{fig:ifm_block}
\end{subfigure}

\end{minipage}
\hfill
\begin{minipage}{0.12\textwidth}
\centering
\includegraphics[width=\linewidth]{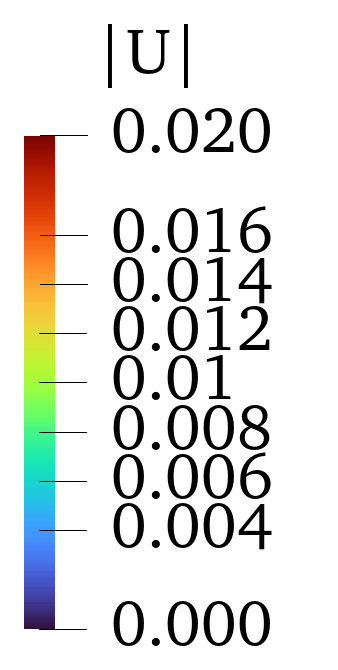}
\end{minipage}

\caption{
Comparison of block-wise partitions before and after conformalization.
The SCZM mesh provides a pixelated approximation of the material regions, while the generated IFM mesh recovers an interface-aligned representation.
}
\label{fig:sczm_block_to_ifm_block}
\end{figure}

This reconstruction also produces a consistent block-wise partition, as shown in \figref{fig:sczm_block_to_ifm_block}, where the pixelated SCZM representation is replaced by an interface-aligned geometry.

\begin{figure}[htb]
\centering

\begin{minipage}{0.82\textwidth}
\centering

\begin{subfigure}[b]{0.48\textwidth}
    \centering
    \includegraphics[width=\linewidth,trim = {300 10 1200 100}, clip]{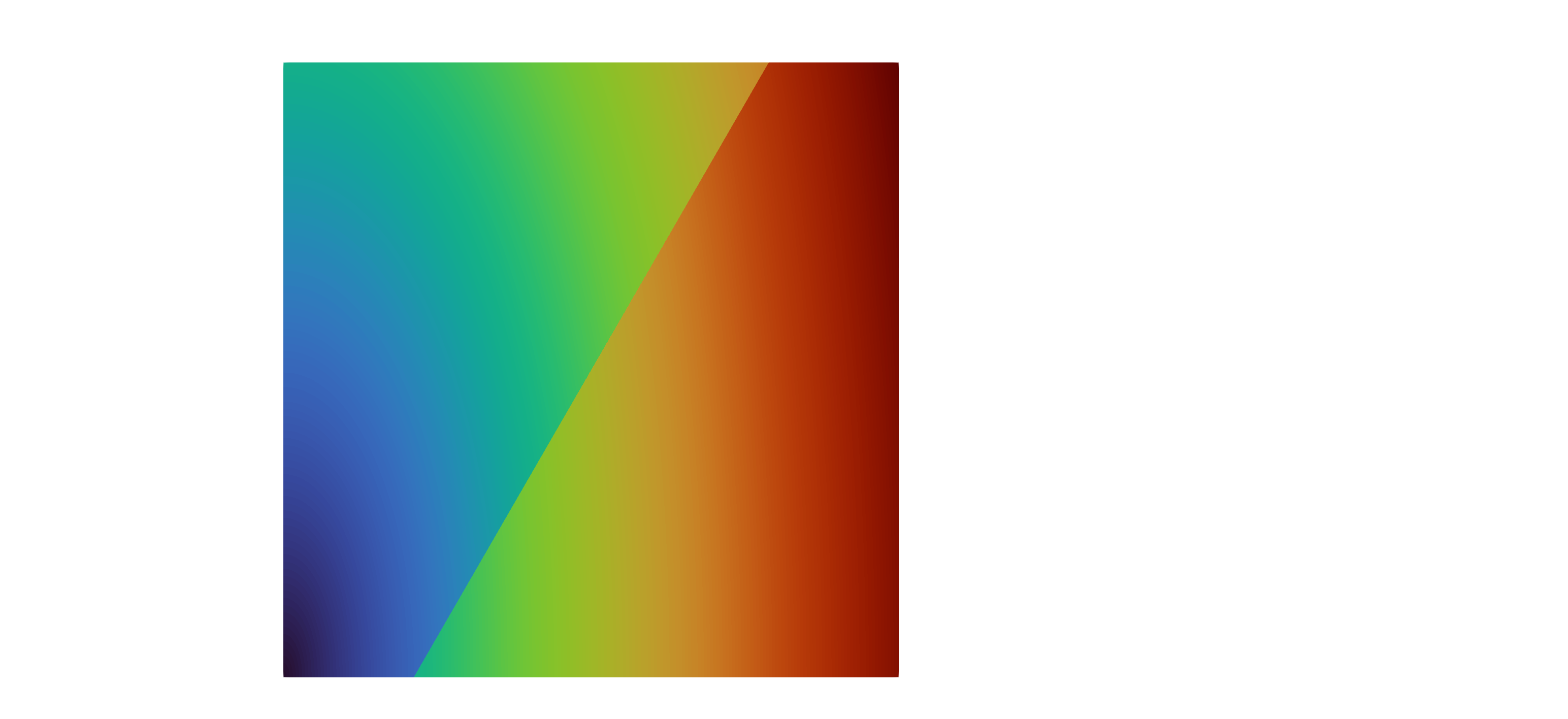}
    \caption{Reference interface-fitted solution.}
    \label{fig:bfm_ref}
\end{subfigure}
\hfill
\begin{subfigure}[b]{0.48\textwidth}
    \centering
    \includegraphics[width=\linewidth,trim = {300 10 1200 100}, clip]{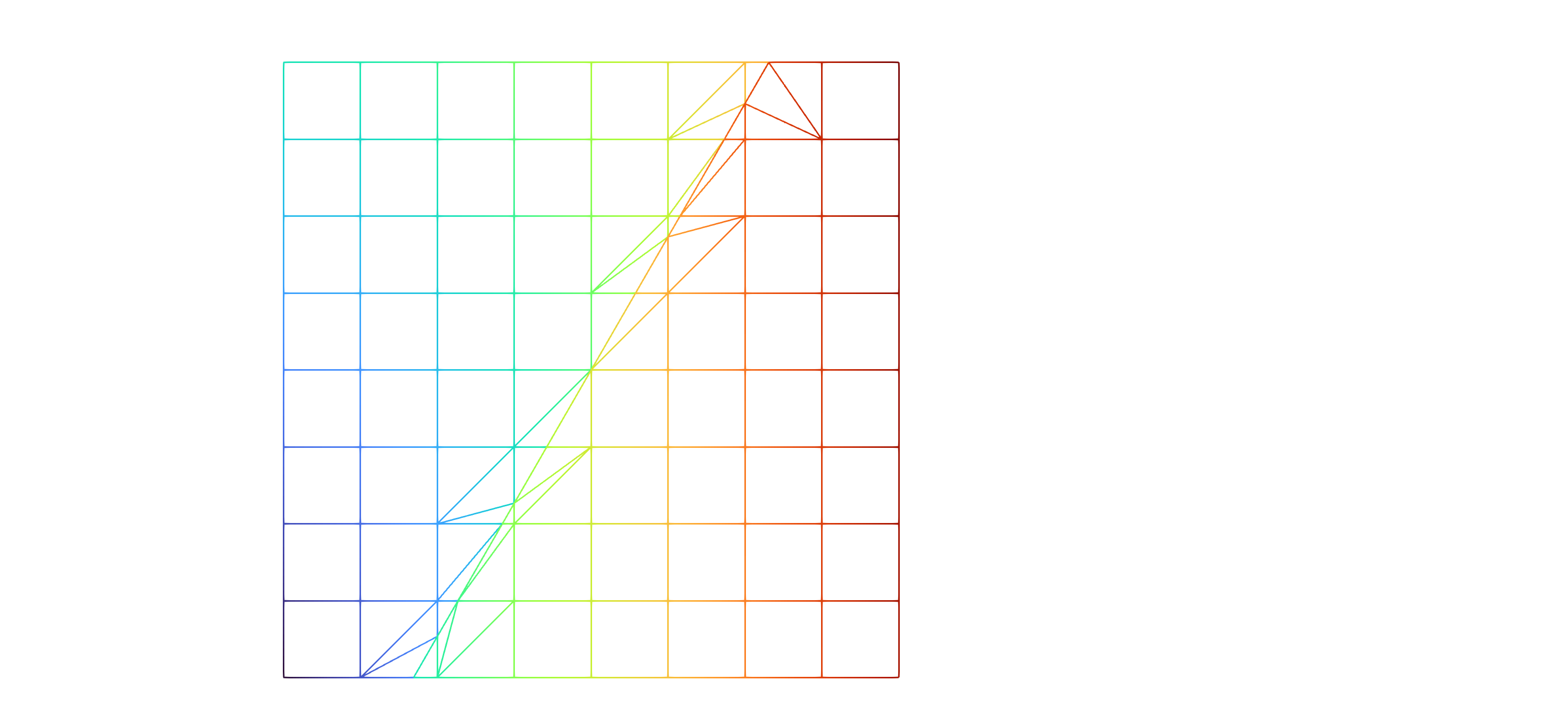}
    \caption{Projected IFM solution ($n=8$).}
    \label{fig:sczm_n8}
\end{subfigure}

\vspace{0.5em}

\begin{subfigure}[b]{0.48\textwidth}
    \centering
    \includegraphics[width=\linewidth,trim = {300 10 1200 100}, clip]{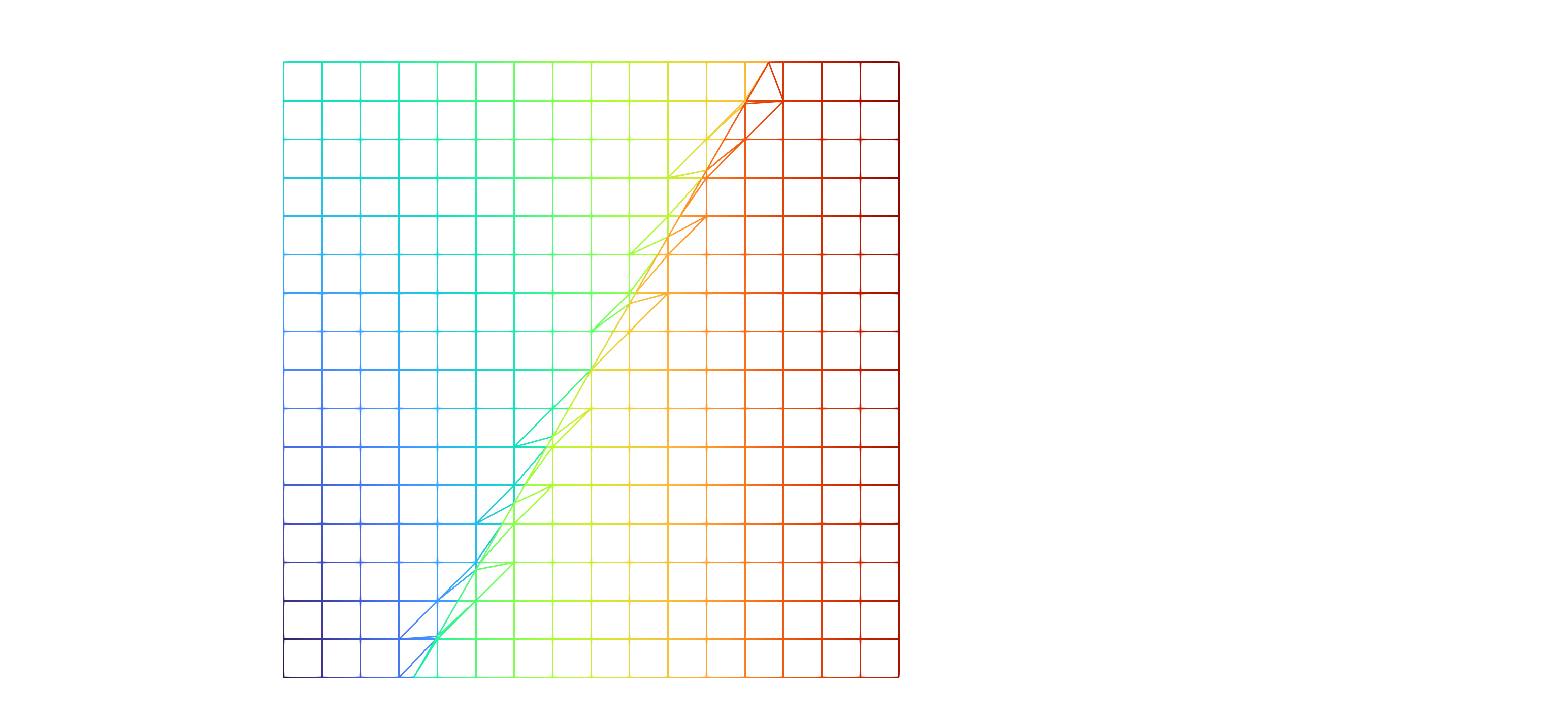}
    \caption{Projected IFM solution ($n=16$).}
    \label{fig:sczm_n16}
\end{subfigure}
\hfill
\begin{subfigure}[b]{0.48\textwidth}
    \centering
    \includegraphics[width=\linewidth,trim = {300 10 1200 100}, clip]{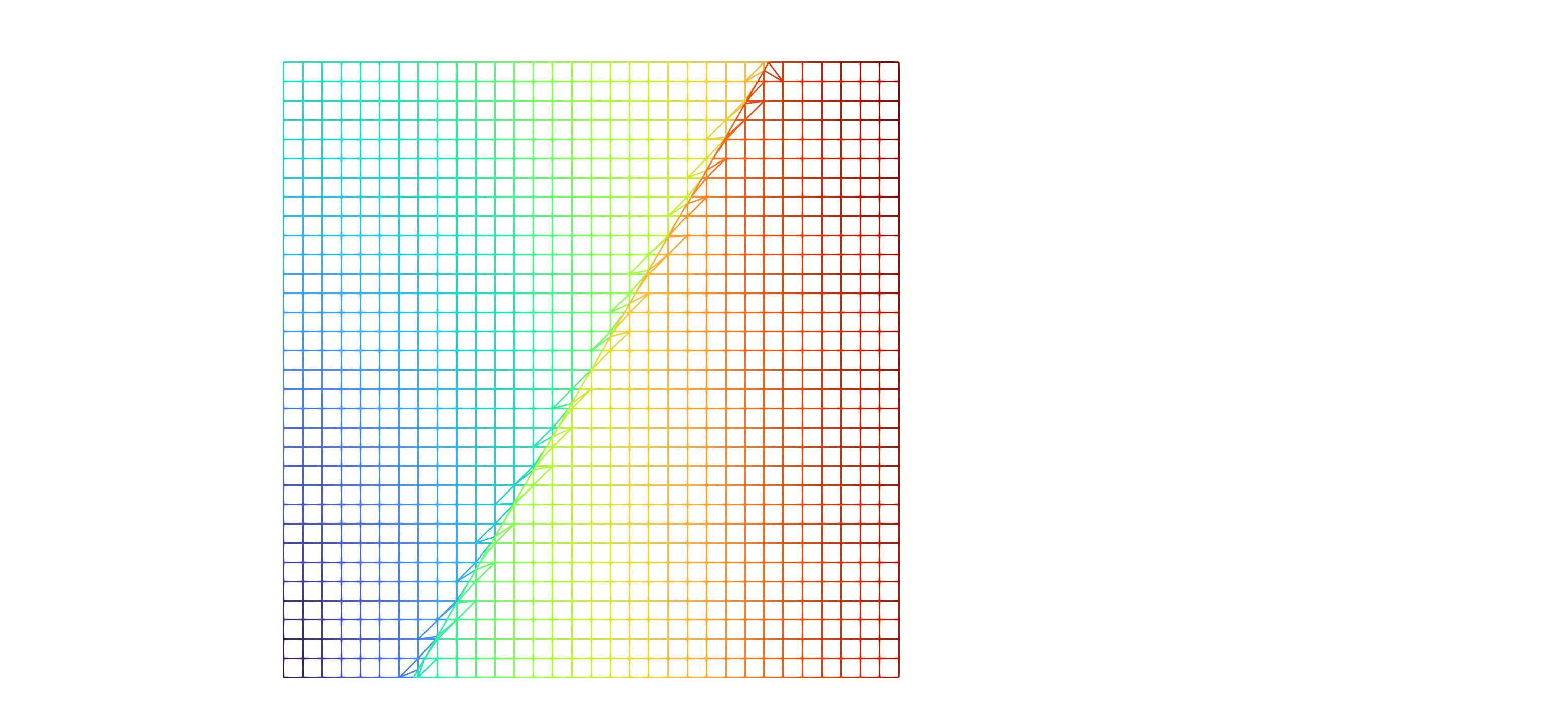}
    \caption{Projected IFM solution ($n=32$).}
    \label{fig:sczm_n32}
\end{subfigure}

\end{minipage}
\hfill
\begin{minipage}{0.13\textwidth}
\centering
\includegraphics[width=\linewidth]{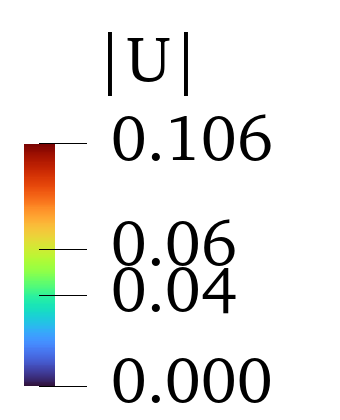}
\end{minipage}

\caption{
Comparison of projected IFM displacement magnitude solutions obtained from SCZM meshes with increasing mesh resolution ($h=\tfrac{1}{n}$).
The reference interface-fitted solution is shown in \subref{fig:bfm_ref}.
}
\label{fig:projection_refinement}

\end{figure}

\begin{figure}[htb]
\centering

\begin{subfigure}[b]{0.48\textwidth}
    \centering
    \includegraphics[width=\linewidth,trim = {300 10 1200 100}, clip]{bfm_t10.png}
    \caption{Reference IFM solution (surface).}
    \label{fig:bfm_surface}
\end{subfigure}
\hfill
\begin{subfigure}[b]{0.48\textwidth}
    \centering
    \includegraphics[width=\linewidth,trim = {300 10 1200 100}, clip]{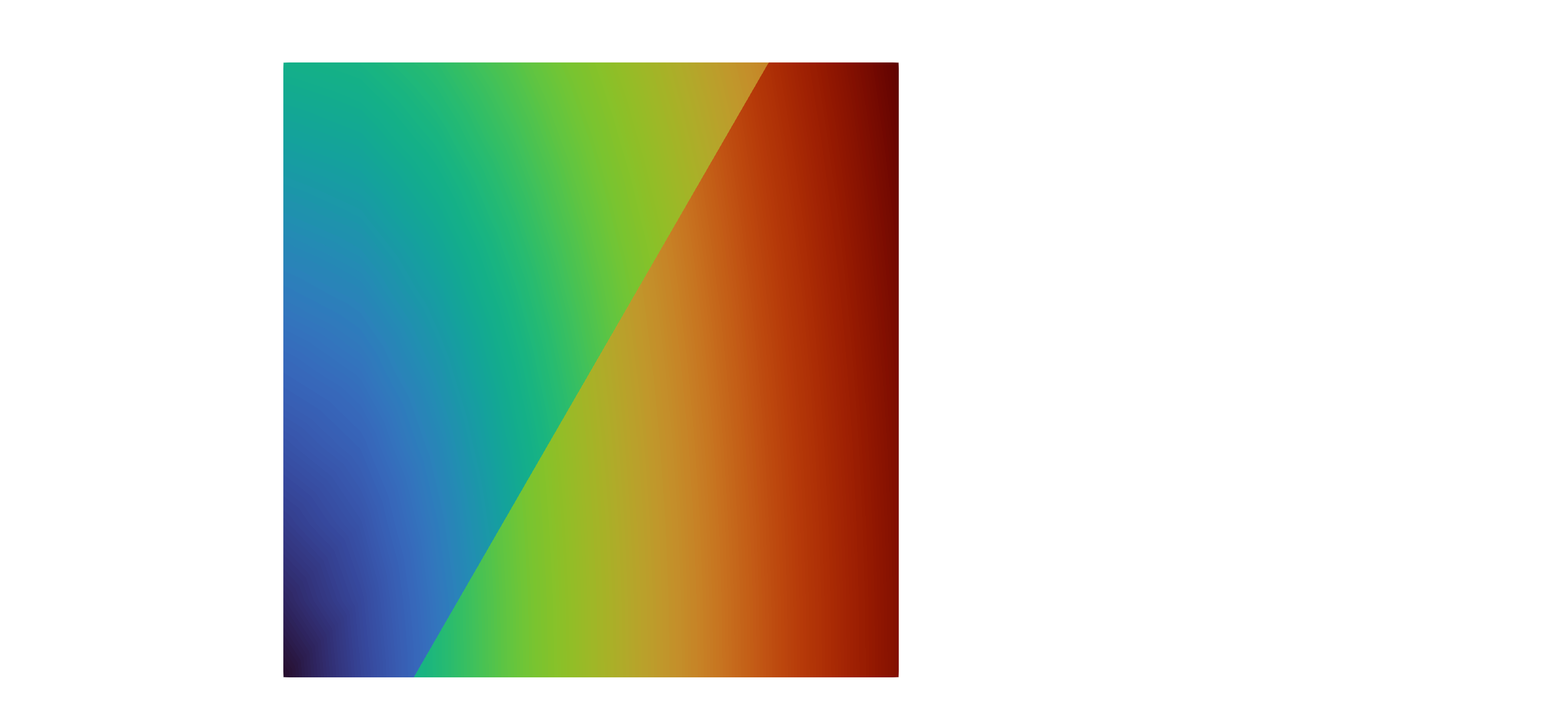}
    \caption{Projected SCZM solution ($n=8$, surface).}
    \label{fig:sczm_surface_n8}
\end{subfigure}

\caption{
Surface comparison between the reference interface-fitted solution and the projected SCZM solution on a coarse mesh ($n=8$).
The two solutions are visually indistinguishable, indicating that the projection preserves the solution field accurately.
}
\label{fig:surface_comparison}
\end{figure}

The corresponding projected solutions are presented in \figref{fig:projection_refinement}, where increasing the mesh resolution from $n=8$ to $n=32$ leads to progressively improved agreement with the reference solution computed on an interface-fitted mesh (\figref{fig:bfm_ref}); a surface comparison for the coarse mesh further confirms that the projected solution is visually indistinguishable from the reference (\figref{fig:surface_comparison}).

\begin{figure}[htb]
  \centering

  \begin{subfigure}[b]{0.48\textwidth}
    \centering
    \includegraphics[width=\linewidth,trim = {600 10 600 100}, clip]{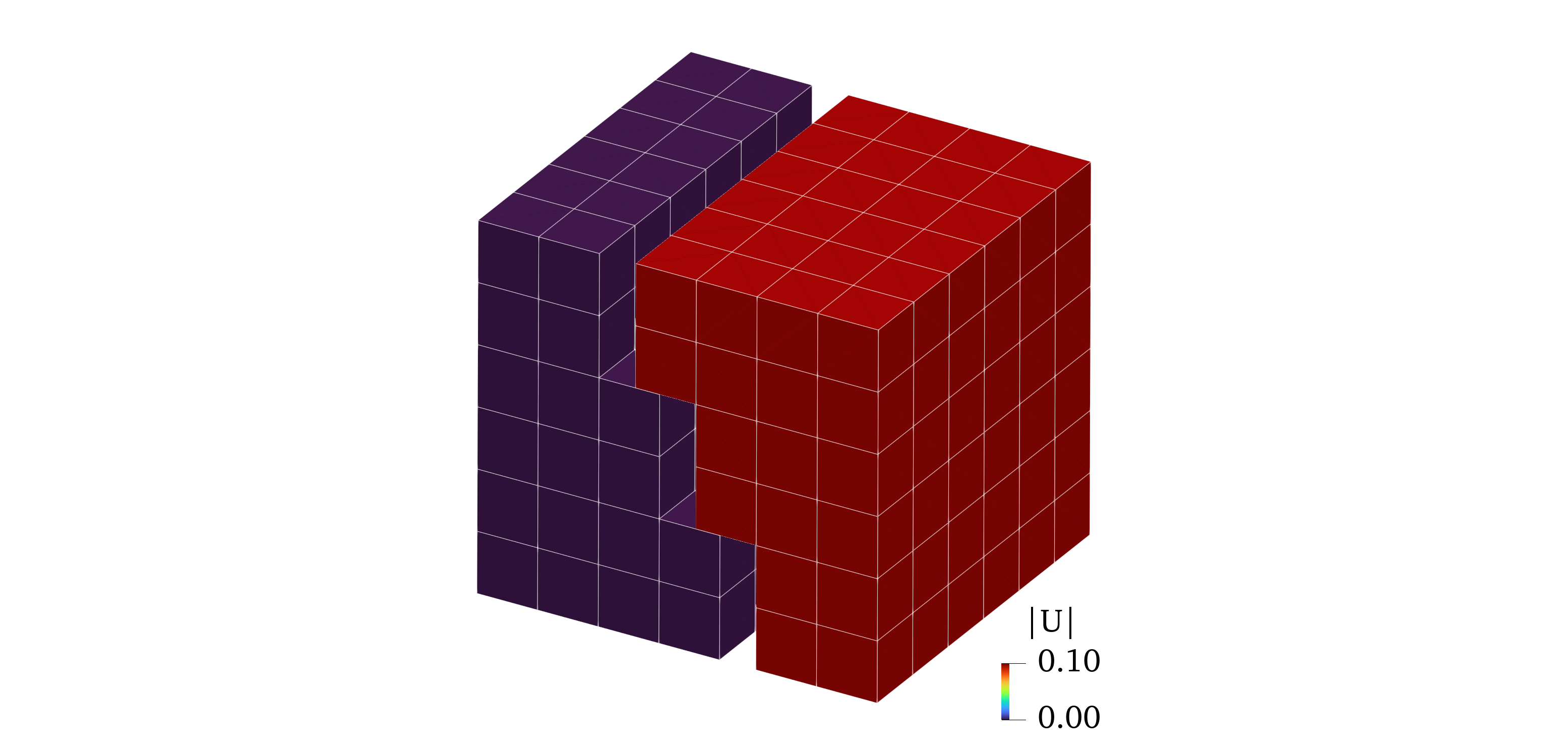}
    \caption{SCZM solution on the original non-interface-fitted mesh.}
    \label{fig:sczm_original_mesh}
  \end{subfigure}
  \hfill
  \begin{subfigure}[b]{0.48\textwidth}
    \centering
    \includegraphics[width=\textwidth,trim={600 10 600 100},clip]{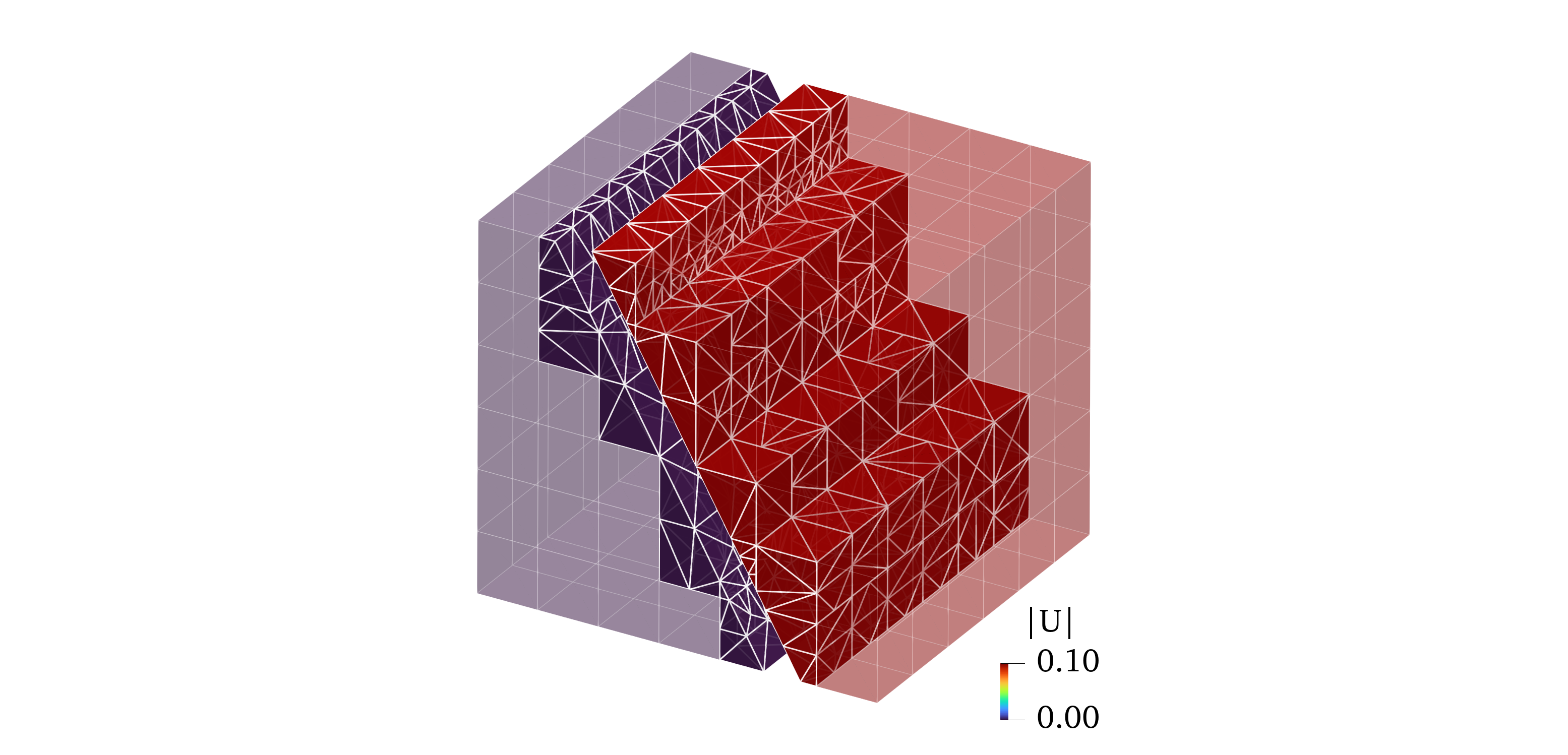}
    \caption{Projected IFM solution after local conformalization of the SCZM mesh.}
    \label{fig:ifm_conformalized_mesh}
  \end{subfigure}

  \caption{
  Illustration of the projection from the original SCZM mesh to the generated IFM mesh.
  The SCZM solution is computed on the original non-interface-fitted pixelated mesh.
  The IFM mesh is then obtained by locally conformalizing the SCZM mesh near the material interface, followed by projecting the SCZM solution onto the generated interface-fitted nodes.
  }
  \label{fig:sczm_to_ifm_projection}
\end{figure}

Finally, the overall workflow is illustrated in three dimensions in \figref{fig:sczm_to_ifm_projection}, where the SCZM solution is first computed on the non-interface-fitted mesh and then transferred to the interface-fitted mesh generated by the conformalization procedure.

\begin{algorithm}[htb]
\footnotesize
\caption{\textsc{ConformalizePixelatedMeshWithBoundaryRegions}}
\begin{algorithmic}[1]
\Require Block-wise closed material regions $\{\Omega_N\}$ and pixelated SCZM mesh $\mathcal{T}_{\mathrm{SCZM}}$
\Ensure Conformalized interface-fitted mesh $\mathcal{T}_{\mathrm{IFM}}$, block IDs, and node markers

\State Construct the block-wise material regions $\{\Omega_N\}$
\State Read the pixelated SCZM mesh $\mathcal{T}_{\mathrm{SCZM}}$
\State Initialize the output element list and the block-wise node registry

\For{each element $K \in \mathcal{T}_{\mathrm{SCZM}}$}
    \State Compute all non-empty clipped pieces $(K_N,N)$, where $K_N = K \cap \Omega_N$
    \If{no clipped piece exists}
        \State Discard $K$
    \ElsIf{exactly one clipped piece $(K_N,N)$ exists and $K_N \approx K$}
        \State Copy $K$ directly to the output mesh
        \State Assign block ID $N$ to the copied element
    \Else
        \For{each clipped piece $(K_N,N)$}
            \State Decompose $K_N$ into simplices
            \Comment{triangles in 2D; tetrahedra in 3D}
            \State Add the generated simplices to the output mesh with block ID $N$
        \EndFor
    \EndIf
\EndFor

\State Merge coincident output nodes within each block
\State Duplicate coincident nodes across different blocks
\State Classify each output node and assign its \texttt{node\_marker}
\State Assemble all output elements into $\mathcal{T}_{\mathrm{IFM}}$
\State Write $\mathcal{T}_{\mathrm{IFM}}$, block IDs, and node markers to file
\end{algorithmic}
\label{alg:ConformalizePixelatedMeshWithBoundaryRegions}
\end{algorithm}

\begin{algorithm}[htb]
\footnotesize
\caption{\textsc{BuildSCZMSearchStructure}}
\begin{algorithmic}[1]
\Require SCZM mesh $\mathcal{T}_{\mathrm{SCZM}}$
\Ensure Centroid-based k-d tree $T$, direct lookup structure $S_{\mathrm{lookup}}$, and flag \texttt{is\_uniform}

\State Initialize an empty centroid list $\mathcal{C}$

\For{each SCZM element $e \in \mathcal{T}_{\mathrm{SCZM}}$}
    \State Compute the element centroid $\mathbf{c}_e$
    \State Store $(\mathbf{c}_e,e)$ in $\mathcal{C}$
\EndFor

\State Build a k-d tree $T$ from the centroid list $\mathcal{C}$

\State Check whether $\mathcal{T}_{\mathrm{SCZM}}$ is a structured uniform mesh

\If{$\mathcal{T}_{\mathrm{SCZM}}$ is structured and uniform}
    \State Build the direct lookup structure $S_{\mathrm{lookup}}$
    \State Set \texttt{is\_uniform} $=$ \texttt{true}
\Else
    \State Set $S_{\mathrm{lookup}} =$ \texttt{null}
    \State Set \texttt{is\_uniform} $=$ \texttt{false}
\EndIf

\State \Return $T$, $S_{\mathrm{lookup}}$, \texttt{is\_uniform}
\end{algorithmic}
\label{alg:BuildSCZMSearchStructure}
\end{algorithm}

\begin{algorithm}[htb]
\footnotesize
\caption{\textsc{ProjectSolutionToIFM}}
\begin{algorithmic}[1]
\Require Target IFM nodes $\{\mathbf{x}_t\}$, node markers, \texttt{source\_node\_id}, SCZM mesh $\mathcal{T}_{\mathrm{SCZM}}$, SCZM solution $\mathbf{u}$, block IDs, search structure $T$, $S_{\mathrm{lookup}}$, and flag \texttt{is\_uniform}
\Ensure Projected solution $\mathbf{u}_t$ on the IFM mesh

\For{each target node $\mathbf{x}_t$ in target block $N$}

    \If{\texttt{node\_marker}$(\mathbf{x}_t)=0$}
        \State Obtain the corresponding SCZM node from \texttt{source\_node\_id}
        \State Verify that the source-node block is compatible with the target block $N$
        \State Copy the SCZM nodal solution directly to $\mathbf{u}_t(\mathbf{x}_t)$
        \State \textbf{continue}
    \EndIf

    \If{\texttt{node\_marker}$(\mathbf{x}_t)=1$}
        \State Find a containing SCZM element $e$ in block $N$ using \textsc{FindSourceElem}
        \If{$e$ is found}
            \State Evaluate $\mathbf{u}_t(\mathbf{x}_t)$ by shape-function interpolation on $e$
            \State \textbf{continue}
        \EndIf
    \EndIf

    \State Find the closest SCZM element $e^\ast$ whose block ID is $N$
    \State Select a reference point $\mathbf{x}^\ast$ associated with $e^\ast$
    \Comment{e.g., closest node, closest quadrature point, or element centroid}
    \State Evaluate $\mathbf{u}(\mathbf{x}^\ast)$ from the SCZM solution on $e^\ast$
    \State Evaluate the solution gradient $\nabla \mathbf{u}(\mathbf{x}^\ast)$ from the shape functions of $e^\ast$
    \State Set
    \[
        \mathbf{u}_t(\mathbf{x}_t)
        =
        \mathbf{u}(\mathbf{x}^\ast)
        +
        \nabla \mathbf{u}(\mathbf{x}^\ast)
        \left(\mathbf{x}_t-\mathbf{x}^\ast\right)
    \]

\EndFor

\State \Return $\mathbf{u}_t$
\end{algorithmic}
\label{alg:ProjectSolutionToIFM}
\end{algorithm}

\begin{algorithm}[htb]
\footnotesize
\caption{\textsc{FindSourceElem}}
\begin{algorithmic}[1]
\Require Target point $\mathbf{x}_t \in \mathbb{R}^d$, target block $N$, SCZM mesh $\mathcal{T}_{\mathrm{SCZM}}$, centroid-based k-d tree $T$, direct lookup structure $S_{\mathrm{lookup}}$, and flag \texttt{is\_uniform}
\Ensure Containing source element $e$ in block $N$, or \texttt{null} if no containing element is found

\If{\texttt{is\_uniform} = \texttt{true}}
    \State Compute the structured grid multi-index $\mathbf{i}=(i_1,\ldots,i_d)$ associated with $\mathbf{x}_t$
    \State Query $S_{\mathrm{lookup}}$ using $\mathbf{i}$ to obtain a candidate source element $e$
    \If{$e$ exists and the block ID of $e$ is $N$}
        \State Check whether $\mathbf{x}_t$ lies inside $e$
        \If{yes}
            \State \Return $e$
        \EndIf
    \EndIf
\EndIf

\State Query the centroid-based k-d tree $T$ to obtain nearby candidate elements of $\mathbf{x}_t$

\For{each candidate element $e$}
    \If{the block ID of $e$ is not $N$}
        \State \textbf{continue}
    \EndIf

    \If{$\mathbf{x}_t$ lies inside $e$}
        \State \Return $e$
    \EndIf
\EndFor

\State \Return \texttt{null}
\end{algorithmic}
\label{alg:FindSourceElem}
\end{algorithm}

\end{document}